\documentclass{article}

\usepackage{amsthm, mathrsfs, url} 
\usepackage[margin = 25mm]{geometry} 

\newtheorem{theorem}{Theorem}[section]
\newtheorem{proposition}[theorem]{Proposition}

\newtheorem{remark}[theorem]{Remark}
\newtheorem{definition}[theorem]{Definition}
\newtheorem{lemma}[theorem]{Lemma}

\usepackage{amsmath, amsfonts, amscd, amssymb, bm, mathtools}
\usepackage{enumerate}
\usepackage{tikz}
\usepackage{array}
\usepackage{booktabs}
\newcommand{\N}{\mathbb{N}}

\newcommand{\cA}{{\mathcal{A}}}

\newcommand{\cR}{\mathcal{R}}
\newcommand{\cI}{\mathcal{I}}
\newcommand{\cS}{\mathcal{S}}

\newcommand{\cD}{\mathscr{D}}

\newcommand{\R}{\mathbb{R}}

\newcommand{\bu}{\mathbf u}
\newcommand{\bU}{\mathbf U}

\newcommand{\bv}{\mathbf v}

\newcommand{\bF}{\mathbf F}
\newcommand{\bD}{\mathbf D}
\newcommand{\be}{\mathbf e}
\newcommand{\bfeta}{{\boldsymbol{\eta}}}
\newcommand{\bfETA}{\boldsymbol{H}}
\newcommand{\bfphi}{\boldsymbol{\varphi}}

\newcommand{\bfxi}{\boldsymbol{\xi}}
\newcommand{\bftau}{\boldsymbol{\tau}}
\newcommand{\bfsigma}{\boldsymbol{\sigma}}
\newcommand{\bn}{\boldsymbol{n}}

\newcommand{\bE}{\mathbf{E}}
\newcommand{\bY}{\mathbf{Y}}

\newcommand{\Om}{\Omega}

\newcommand{\eps}{\epsilon}
\def\ltt{\lesssim}

\newcommand{\dt}{{\Delta t}}
\newcommand{\Th}{{\mathcal{T}_h}}

\DeclareMathOperator{\dist}{dist}

\DeclareMathOperator{\sdist}{Sdist}
\DeclareMathOperator{\supp}{supp}
\newcommand{\dx}{\mathrm{d}x}
\newcommand{\Dt}{\mathrm{d}t} 
\newcommand{\dtee}{\mathrm{d}_t}
\newcommand{\dtt}{\mathrm{d}_{tt}}

\usepackage{todonotes}

\usepackage{comment}
\excludecomment{confidential} 

\allowdisplaybreaks

\begin{document}



\title
{Analysis and finite element approximation of a diffuse interface approach to the Stokes--Biot coupling}

\author{Francis R.~A.~Aznaran\thanks{
        Department of Applied and Computational Mathematics and Statistics,
        University of Notre Dame,
            46556, Indiana, United States, \texttt{$\{$faznaran, mbukac$\}$@nd.edu}
    } ~~~
    Martina Buka{\v{c}}\footnotemark[1] ~~~
Boris Muha\thanks{
        Department of Mathematics, Faculty of Science, University of Zagreb,
            Croatia, \texttt{borism@math.hr}
    } ~~~
Abner J.~Salgado\thanks{
        Department of Mathematics,
        University of Tennessee,
            Knoxville, 
            Tennessee,
            37996,
            United States, \texttt{asalgad1@utk.edu}
    }
}





\maketitle 

\begin{abstract}
We consider
the interaction between a poroelastic structure, described using the Biot model in primal form, and a free-flowing fluid, modelled with the time-dependent incompressible Stokes equations.
We propose a diffuse interface model in which a phase field function is used to write each integral in the weak formulation of the coupled problem on the entire domain containing both the Stokes and Biot regions. The phase field function continuously transitions from one to zero over a diffuse region of width $\mathcal{O}(\varepsilon)$ around the interface; this allows the equations to be posed uniformly across the domain, and obviates tracking the subdomains or the interface between them.
We prove convergence in weighted norms of a finite element discretisation of the diffuse interface model to the continuous diffuse model; here the weight is a power of the distance to the diffuse interface. We in turn prove convergence of the continuous diffuse model to the standard, sharp interface, model.
Numerical examples verify the proven error estimates, and illustrate application of the method to 
fluid flow through a complex network, describing blood circulation in the circle of Willis.
\end{abstract}


\section{Introduction}

The coupling of the Stokes and Biot equations describes the interaction between a free-flowing fluid, and a poroelastic material, 
where the latter is itself governed by a mechanical law describing the elastic phase, and Darcy's law describing the fluid phase.
It is used to describe problems arising in many applications, including the environmental sciences, hydrology, geomechanics, and biomedical engineering.
Both theoretically and computationally, it 
inherits the challenges associated with both the Stokes--Darcy fluid-fluid coupling, and classical fluid-structure interaction problems.
There has been a recent surge of interest in the numerical analysis and development of numerical schemes for fluid-poroelastic structure interaction problems~\cite{Ambartsumyan2018, Ruiz2022, Seboldt2021}, and their well-posedness analyses~\cite{Bociu2021, Kuan2024, Avalos2024}.

In classical approaches, the computational mesh is often aligned with the interface between the two regions; 
we refer to this as the ``sharp'' interface approach.
The diffuse interface method (also known as the phase field, or diffuse domain, approach)~\cite{Anderson1998, Burger2017, Du2020}
instead 
uses phase field functions, parameterised by a variable $\epsilon$, as approximate indicator functions of each subdomain, 
taking the value of unity on most of one domain, zero in most of the other, and then continuously but rapidly transitioning between those two values in a diffuse region of width $\mathcal{O}(\eps)$,
which approximates the interface between the subdomains. 
The resulting solution variables 
can in some instances
be proven to
converge in the limit $\eps\to 0$ to those obtained via the sharp interface formulation.
In the multiphase context, the coupled regions undergo different physics and are described by different equations, 
so the phase field does more than simply interpolate between different coefficients in the same equation.

The diffuse interface approach 
allows the integrals in the weak formulation to be written over a single domain (in this context, the union of the fluid and poroelastic subdomains), 
so that the integration is formally done without reference to the subdomains. As such, there is 
no need to track the locations of each subdomain 
or 
the interface between them,
which are instead implicitly stored in the phase field function.
This in particular allows for increasingly complex domain geometry 
but 
without 
a comparable increase in mesh complexity. 
Even though in this work we consider infinitesimal displacements, and a linear fluid-poroelastic structure interaction problem (the solid and fluid domains remain fixed), this approach is still attractive in more complex applications.
As one example, we mention medical application settings, 
where
an appropriate 
phase field function may be 
inferred 
directly from the colour pixels in an image. In such cases, the domains and interfaces in question  are highly approximative and noisy 
anyway,
in particular incorporating interfaces of nonzero ``width'',
so that 
the exact interface may have to be reconstructed from this smoothed region, or
even the notion of sharp interface solution may no longer be clear.
Furthermore, we intend for this study 
to enable future investigation of
regimes where the deformations are no longer 
assumed 
infinitesimal, and so the solid subdomain can undergo large deformation and/or evolve in time, and may possibly incorporate contact phenomena; 
approximation of the sharp interface solution would thus require the expense of re-meshing the domain, 
so that it becomes increasingly difficult to justify resolving the domains exactly via the sharp interface approach.

The diffuse interface approach 
is popular in engineering contexts (e.g.~\cite{Saylor2016, Stoter2017, Mokbel2018}
),
and now encompasses a broad range of related methods such as 
the level set method~\cite{Osher2003}
and the diffuse Nitsche method~\cite{Nguyen2018},
but 
the underlying theory has yet to catch up with its prevalence in applications. In particular, the 
critical 
question of whether, and 
in what sense, 
the diffuse solutions converge to the sharp solutions in the limit as the interface width tends to zero, is not yet known for many phase field methods;
we term this the convergence of \textit{modelling error}. Theoretically, the diffuse interface approach has been widely studied for elliptic problems and two-phase flow problems~\cite{Burger2017, Burger2015, Schlottbom2016, Franz2012, Li2009, Nguyen2018, Abels2024}, but not as much for multiphysics coupled problems. 
The diffuse interface approach is appealing from the point of view of analysis, as it 
sidesteps issues concerning the singularity or regularity of the interface, which are typically the main bottlenecks for theoretical analysis of fluid-structure interaction, at the cost of working 
in the nonstandard functional setting of
weighted Sobolev spaces.
Namely, the phase field function may be interpreted as weighting the integrals arising in the variational formulation of the 
governing 
equations, so that 
(as observed in e.g.~\cite{Abels2015})
it is natural to pose the resulting variational formulations in Sobolev spaces weighted by the phase field. 
This approach was used in~\cite{Bukac2023} to analyse
the convergence of the diffuse interface model to the sharp interface model for the Stokes--Darcy problem. 
The strong convergence of the modelling error
in Sobolev norms
has been rigorously proved in~\cite{Abels2018} for the Stokes/Allen--Cahn system. We also mention recent work by Abels~\cite{Abels2021, Abels2021a} on convergence and approximate solutions of the Stokes/Cahn--Hilliard system. 

In this paper, we formulate sharp and diffuse interface formulations for the interaction between a fluid and a poroelastic structure. We show that the diffuse interface problem is well-posed. 
In the well-posedness theory for the sharp interface problem, one of the challenges is 
that the structure velocity does not have enough regularity to admit a well-defined trace.
While this obstacle is not present in the same form in the diffuse interface problem, it manifests in a different way, which we resolve by using a specific type of trace inequality which allows us to obtain an estimate uniform in $\epsilon$. To show that 
our discrete 
solution converges to the continuous sharp interface solution, we split the error into two parts: the error between the discrete and continuous diffuse interface solutions (approximation error), and the error between the continuous diffuse and sharp interface solutions (modelling error). We use energy estimates in weighted spaces to obtain 
convergence rates with respect to time and space discretisation parameters for the diffuse interface problem
assuming $\epsilon > 0$. Then we derive rates of convergence of the continuous diffuse interface solutions to the continuous sharp interface solution with respect to $\epsilon$. 

This work is structured as follows. 
In Section~\ref{sec:Notation}, we collect notation and results concerning weighted Sobolev spaces.
Section~\ref{sec:Models} states the standard (sharp interface) model of the Stokes--Biot coupling, before introducing
the 
diffuse interface formulation 
and proving well-posedness of the diffuse weak formulation;
we then prove well-posedness and convergence of a finite element discretisation thereof in Section~\ref{sec:Discretisation}.
We derive rates of convergence of the continuous diffuse interface formulation to the continuous diffuse interface formulation in Section~\ref{sec:ModellingError}, for both an exact and regularised phase field function.
Finally, Section~\ref{sec:numerics} demonstrates the proven rates of convergence of the numerical scheme,
compares the scheme with the sharp interface method in a 3D example,
and finally 
applies the method to 
fluid flow through a complex 
3D 
network, describing blood circulation in a patient-specific model of the circle of Willis
in the brain.

\section{Notation and preliminaries}\label{sec:Notation} 

During the course of our discussion $\Omega \subset \mathbb{R}^d$ with $d\in\{2, 3\}$ is an open bounded domain with Lipschitz boundary. When dealing with discretisation, we shall also assume that it is polytopal. This will be used to denote the domain where the fluid and porous medium interaction takes place. 
Let $\lesssim$ denote domination up to a constant which may depend on mesh regularity, but not on any
other
discretisation or phase field parameters.
We mostly adhere to standard notation and terminology with regards to function spaces and their properties. 
For a Banach space $X$ we denote its dual by $X^*$.

We say that an a.e.~positive function $\omega \in L^1_{\text{loc}}(\mathbb{R}^d)$ is a weight. Every weight $\omega$ induces a measure with density $\omega \dx$ over the Borel subsets of $\mathbb{R}^d$, which for simplicity will also be denoted by $\omega$. In other words, for $E \subset \mathbb{R}^d$ a Borel set, we let $\omega(E) = \int_E \omega \dx$.
For $r\in (1, \infty)$, $\omega$ a weight, and $D \subset \mathbb{R}^d$ a bounded domain, we define weighted Lebesgue spaces and their norms:
\begin{equation*}
    L^r (D, \omega) \coloneqq \left\{ \psi : D \rightarrow \mathbb{R} \; : \; |\psi|^r \omega \in L^1(D) \right\},
    \qquad
    \| \psi \|_{L^r(D, \omega)}^r \coloneqq \int_D |\psi|^r \omega.
\end{equation*}
Associated with weighted $L^r$--spaces, we define the weighted Sobolev spaces
\begin{equation*}
    \begin{aligned}
        W^{k, r}(D, \omega) &\coloneqq \left\{ \psi \in L^r (D, \omega) \; : \; \partial^{\alpha} \psi \in L^r(D, \omega) \; \forall \alpha \in \mathbb{N}_0^d : |\alpha| \leq k \right\},
        \\
        \|\phi\|_{W^{k, r}(D, \omega)}^r &\coloneqq \sum_{|\alpha|\leq k}\|\partial^{\alpha}\psi\|^r_{L^r(D, \omega)}.
    \end{aligned}
\end{equation*}
All 
these 
spaces are complete.
As usual, we set $H^k(D, \omega) = W^{k, 2}(D, \omega)$ for any $k \in \mathbb{N}_0$. 

While it is possible to develop a fairly general theory of weighted Sobolev spaces \cite{Kufner1984}, in what follows we shall only be concerned with a very specific type of weight, namely, a power of the distance to a piece of the boundary. Let us make this concrete and discuss some consequences of this choice.

Let $D \subset \R^d$ be a bounded domain with Lipschitz boundary. Assume that $\Gamma \subset \partial D$ has positive and finite $(d - 1)$--dimensional Hausdorff measure $\mathcal{H}_{d - 1}(\Gamma)$. For $\alpha \in (0,1)$ define the weight
\begin{equation}\label{eq:distisA2}
    \omega(x) \coloneqq \dist_\Gamma(x)^\alpha,
\end{equation}
where $\dist_\Gamma$ is the distance function to $\Gamma$. Owing to this definition, we have the following properties.
\begin{itemize}
    \item The weight $\omega$ belongs to the Muckenhoupt class $A_2$~\cite{Nochetto2016}; see~\cite{Aimar2014} and~\cite[Lemma 2.3(vi)]{Farwig1997}. As a consequence, for every $k \in \mathbb{N}$, the space  $H^k(D, \omega)$ is Hilbert and separable; see~\cite{Kufner1984, Nochetto2016}.
    \item \emph{Traces}: There is a continuous trace operator $\gamma_\Gamma \colon H^1(D, \omega) \to L^2(\Gamma)$; see~\cite{Nekvinda1993}. Therefore, it is legitimate to define $H^1_\Gamma(D, \omega)$ as the subspace of $H^1(D, \omega)$ of functions which vanish on $\Gamma$.
    \item \emph{Poincar\'e--Friedrichs--type inequalities}: In light of the trace results mentioned above, there is a constant $C_P > 0$ such that, for all $v \in H^1(D, \omega)$, we have
    \begin{equation}\label{eq:WeightedPoincare}
        \| v \|_{H^1(D, \omega)}^2 \leq C_P \left( \left| \int_\Gamma v \right| + \| \nabla v \|_{L^2(D, \omega)^d} \right)^2.
    \end{equation}
    As a consequence, $v \mapsto \| \nabla v \|_{L^2(D, \omega)^d}$ is an equivalent norm on $H^1_\Gamma(D, \omega)$.
    \item \emph{Bogovski{\u\i} operator}: There is, see \cite[Theorem 3.9]{Bukac2023}, $\beta > 0$ such that, for all $\psi \in L^2(D, \omega^{-1})$,
    \begin{equation}\label{eq:Bogovskii}
        \beta \| \psi \|_{L^2(D, \omega^{-1})} \leq \sup_{\boldsymbol0 \neq \bv \in H^1_\Gamma(D, \omega)^d} \frac{\int_D \psi \nabla\cdot \bv}{ \|\nabla \bv \|_{L^2(D, \omega)^{d \times d}}}.
    \end{equation}
    \item \emph{Korn's inequality}: For $\bv \in H^1(D, \omega)^d$ we denote by $\bD(\bv) = \tfrac12 \left( \nabla \bv + \nabla \bv^T \right)$ its symmetric gradient. There is a constant $\overline{C}_K > 0$ such that
    \begin{equation}\label{eq:wKorn1}
        \| \nabla \bv \|_{L^2(D, \omega)^{d \times d}} \leq \overline{C}_K \left( \| \bv \|_{L^2(D, \omega)^d} + \| \bD(\bv) \|_{L^2(D, \omega)^{d \times d}} \right), \quad \forall \bv \in H^1(D, \omega)^d;
    \end{equation}
    see~\cite[Theorem 5.17]{Diening2010}. Mimicking for instance~\cite[Theorem 7.3.2]{Dautray1988}, we conclude that there is a constant $C_K > 0$ such that, for all $\bv \in H_\Gamma^1(D, \omega)^d$, we have
    \begin{equation}\label{eq:weightedKorn}
        \| \bv \|_{H^1(D, \omega)^d}^2 \leq C_K \| \bD(\bv) \|_{L^2(D, \omega)^{d \times d}}^2.
    \end{equation}
\end{itemize}
The constants in all the statements above depend on $\omega$ only through the so-called Muckenhoupt characteristic $[\omega]_{A_2}$.

\section{The mathematical models}\label{sec:Models}

\subsection{The sharp interface model}\label{sub:SharpModel}

Let $\Omega_F$ denote the fluid domain and $\Omega_B$ the poroelastic structure reference domain. These are bounded domains in $\mathbb{R}^d$ which have Lipschitz boundary, and form a partition of $\Omega$. In other words, $\Omega_F \cap \Omega_B = \emptyset$ and $\overline{\Omega} = \overline{\Omega_F} \cup \overline{\Omega_B}$. Finally, we assume that $\Omega_B$ is not encapsulated by $\Omega_F$. By this we mean that both $\partial\Omega_F \cap \partial\Omega$ and $\partial\Omega_B \cap \partial\Omega$ are not only nonempty, but have positive $(d - 1)$--Hausdorff measure. Finally, $\Gamma$ is the common boundary between the fluid and poroelastic structure domains, i.e.~$\Gamma = \partial \overline{\Omega_F} \cap \partial \overline{\Omega_B}$. We assume that $\Gamma$ is sufficiently smooth, and that it has positive and finite $(d - 1)$--Hausdorff measure. Finally, we assume that $\partial \Omega = \Gamma_F^1 \cup \Gamma_F^2 \cup \Gamma_B^1$.

\subsubsection{The fluid}

To model the fluid flow, we use the time-dependent Stokes equations, given as follows:
\begin{align}\label{eq:NS}
    \rho_F \partial_t\bu &= \nabla\cdot\bfsigma_F(\bu, \pi) + \rho_F \bF_F &{\rm in}\ \Omega_F \times (0, T), \\
    \label{eq:Incomp}
    \nabla\cdot\bu &= 0&{\rm in}\;\Omega_F\times (0, T),
\end{align}
where $\bf{u}$ is the fluid velocity, the constant $\rho_F > 0$ the fluid density, $\bfsigma_F(\bu, \pi) = 2\mu_F \bD(\bu) - \pi \pmb{\mathbf{I}}$ the fluid Cauchy stress tensor, 
$\pmb{\mathbf{I}}$ the identity tensor, 
$\mu_F$ the fluid viscosity, $\bD(\bu)$ the strain rate tensor, $\pi$ denotes the
fluid pressure, and $\bF_F$ is the density of volumetric forces.

\subsubsection{The poroelastic material}

We describe the poroelastic material using Biot's poroelasticity equations, and incorporate the viscoelastic properties using the Kelvin--Voigt linear model.
The pressure and the deformation are mutually dependent and fully coupled
in Biot's model, which is given as
follows~\cite{Biot1956, Chen1994, Chen1994a}:
\begin{align}
    \rho_B \partial_{tt} \bfeta &= \nabla \cdot \bfsigma_B(\bfeta, p) + \rho_B\bF_B & \textrm{in}\; \Omega_B \times(0, T)
    , \label{eq:B1}\\
    c_0 \partial_t p + \alpha \nabla \cdot \partial_t\bfeta - \nabla \cdot (\boldsymbol \kappa \nabla p)
    &= g &\textrm{in}\; \Omega_B\times(0, T), \label{eq:B2}
\end{align}
where
$\bfeta$ is the poroelastic structure displacement, and
$p$ is the fluid pore pressure.
The constant $\rho_B > 0$ is
the structure density, $\boldsymbol \kappa$ the permeability
tensor, $c_0$ the storage coefficient, and $\alpha$
the
(constant)
Biot--Willis parameter accounting for the coupling strength between
the fluid and solid phase.
The density of external forces on the structure is denoted by $\bF_B$, and $g$ is a source.
We assume $\boldsymbol\kappa$ is symmetric, and uniformly bounded and positive definite, so that $0 < k_* \leq \lambda \leq k^*$
a.e.~for $\lambda\in\lambda(\boldsymbol\kappa)$, where $\lambda(\boldsymbol\kappa)$ is the spectrum of $\boldsymbol\kappa$.
The total Cauchy stress tensor of the poroelastic medium is constituted by
\begin{equation*}
    \bfsigma_B(\bfeta, p) = \bfsigma_E(\bfeta) - \alpha p \pmb{\mathbf{I}},
\end{equation*}
where $\bfsigma_{E}(\bfeta)$ denotes the
elasticity stress tensor.
For
an isotropic, homogeneous, elastic material, using a
linearised
Saint Venant--Kirchhoff model, we have
\begin{equation}\label{eq:elasticity-tensor}
    \bfsigma_{E}(\boldsymbol \eta) = 2\mu_B\bD(\bfeta) + \lambda_B(\nabla\cdot\bfeta) \pmb{\mathbf{I}} \eqqcolon \mathbb{C}\bD(\bfeta),
\end{equation}
where $\lambda_B$ and $\mu_B$ are the material-dependent Lam\'e's first and second parameters, respectively, and $\mathbb{C}$ is the elasticity tensor.

\subsubsection{Coupling conditions}

To couple the fluid and poroelastic material, we prescribe the following interface conditions~\cite{Showalter2005, Murad2001, Badia2009}:
\begin{itemize}
    \item The conservation of mass:
    \begin{equation*}
        \bu\cdot\bn = (\partial_t \bfeta - \boldsymbol\kappa \nabla p ) \cdot\bn \qquad {\rm{on}}\; \Gamma \times(0, T),
    \end{equation*}
    where $\bn$ is the unit normal to $\Gamma$ which points towards $\Omega_B$.
    \item The Beavers--Joseph--Saffman condition:
    \begin{equation*}
        \alpha_{BJ}(\bu - \partial_t \bfeta)\cdot\bftau_{i} + \bfsigma_F(\bu, \pi) \bn \cdot\bftau_{i} = 0, \qquad i = 1, \dots, d - 1, \qquad{\rm on}\; \Gamma \times(0, T),
    \end{equation*}
    where, for 
    a.e.~$x\in \Gamma$, $\{\boldsymbol{\tau}_{i}\}_{i = 1}^{d - 1}$ is an orthonormal basis for the tangent space $T_x(\Gamma)$, and $\alpha_{BJ} > 0$ is the Beavers--Joseph--Saffman--Jones coefficient.
    \item The balance of pressure:
    \begin{equation*}
        -\bfsigma_F(\bu, \pi)\bn \cdot\bn = p \qquad{\rm on}\; \Gamma \times (0, T).
    \end{equation*}
    \item The balance of contact forces:
    \begin{equation*}
        \bfsigma_F(\bu, \pi)\bn = \bfsigma_B(\bfeta, p)\bn \qquad{\rm on}\; \Gamma \times (0, T).
    \end{equation*}
\end{itemize}

\subsubsection{Boundary and initial conditions}

We split the boundaries as: $\partial\Omega_F = \Gamma\cup\Gamma_F^1\cup\Gamma_F^2$, $\partial\Omega_B = \Gamma\cup\Gamma_B^1$, and prescribe the following boundary conditions:
\begin{equation*}
    \begin{aligned}
        \bu &= \boldsymbol0\quad {\rm on}\quad \Gamma_F^1\times (0, T),
        &
        \bfsigma_F(\bu, \pi) \bn&= \boldsymbol0\quad{\rm on}\quad \Gamma_F^2\times (0, T),
        \\
        \bfeta &= \boldsymbol0\quad {\rm on}\quad \Gamma_B^1\times (0, T),
        &
        \boldsymbol{\kappa} \nabla p \cdot \bn &= 0 \quad {\rm on}\quad \Gamma_B^1\times (0, T).
    \end{aligned}
\end{equation*}
Finally, we supplement the problem with the following initial conditions:
\begin{equation*}
    \bu(\cdot, 0) = \bu_0 \ {\rm in } \; \Omega_F, \qquad
    \bfeta(\cdot, 0) = \bfeta_0 \ {\rm in } \; \Omega_B, \qquad
    \partial_t \bfeta(\cdot, 0) = \bfxi_0 \ {\rm in } \; \Omega_B, \qquad
    p(\cdot, 0) = p_{0} \ {\rm in } \; \Omega_B.
\end{equation*}

\subsubsection{Weak formulation}

To accommodate for boundary conditions, we introduce the following function spaces:
\begin{equation*}
    \mathcal{V}_i \coloneqq \{ \bv \in H^1(\Omega_i)^d \; : \; \bv = \boldsymbol0 \; \textrm{on} \; \Gamma_i^1 \}, \quad i\in\{F, B\},
    \quad
    \mathcal{X} \coloneqq H^1(\Omega_B).
\end{equation*}

\begin{definition}[weak solution of the sharp interface problem]\label{def:WS}
We say that the tuple $(\bu, \pi, \bfeta, p )$ is a weak solution to our problem if
\begin{equation*}
    \begin{aligned}
        \bu &\in L^\infty(0, T; L^2(\Omega_F)^d) \cap L^2(0, T; \mathcal{V}_F),
        & \bfeta & \in L^{\infty}(0, T; \mathcal{V}_B)\cap W^{1, \infty}(0, T; L^2(\Omega_B)^d),
        \\
        \pi &\in H^{-1}(0, T; L^2(\Omega_F)),
        & p &\in L^\infty(0, T; L^2(\Omega_B)) \cap L^2(0, T; \mathcal{X}),
    \end{aligned}
\end{equation*}
$\bfeta(0) = \bfeta_0$, and, in addition, for every tuple $(\bv, \zeta, \bfphi, q)$ such that
\begin{equation*}
    \begin{aligned}
        \bv &\in C^1_c([0, T); \mathcal{V}_F),
        & \bfphi &\in C^1_c([0, T); \mathcal{V}_B), 
        \\
        \zeta &\in C^1_c([0, T); L^2(\Omega_F)),
        & q &\in C^1_c([0, T); \mathcal{X}),
    \end{aligned}
\end{equation*}
the following equality is satisfied:
\begin{equation}\label{eq:weakSI}
    \begin{aligned}
        &- \rho_F \int_0^T\int_{\Omega_F} \bu \cdot \partial_t \bv
        - \rho_B \int_0^T\int_{\Omega_B} \partial_{t} \bfeta\cdot\partial_t\bfphi
        - \int_0^T\int_{\Omega_B}(c_0 p + \alpha\nabla\cdot \bfeta)\partial_t q
        + 2\mu_F\int_0^T\int_{\Omega_F}\bD(\bu):\bD(\bv)
        \\ &
        - \langle \pi, \nabla\cdot\bv \rangle_{(0, T) \times \Omega_F}
        + \int_0^T\int_{\Omega_F} \zeta \nabla\cdot\bu
        + \int_0^T\int_{\Omega_B}\bfsigma_B(\bfeta, p):\bD(\bfphi)
        + \int_0^T\int_{\Omega_B} \boldsymbol{\kappa}\nabla p \cdot \nabla q
        \\ &
        - \int_0^T\int_\Gamma (q\bu + \partial_t q\bfeta)\cdot\bn
        + \int_0^T\int_{\Gamma} p (\bv - \bfphi) \cdot \bn
        + \alpha_{BJ} \sum_{i = 1}^{d - 1} \int_0^T\int_{\Gamma}(\bu \cdot\bftau_i) (\bv - \bfphi)\cdot\bftau_i
        + (\bfeta \cdot\bftau_i )\partial_t(\bv - \bfphi)\cdot\bftau_i
        \\
        &= \rho_F\int_0^T\int_{\Omega_F}\bF_F\cdot\bv
        + \rho_B\int_0^T\int_{\Omega_B}\bF_B\cdot\bfphi
        + \int_0^T\int_{\Omega_B}gq
        + \rho_F\int_{\Omega_F}\bu_0\cdot\bv(0)
        + \rho_B \int_{\Omega_B}\bfxi_0\cdot\bfphi(0)
        \\ &
        + \int_{\Omega_B}(c_0p_{0} + \nabla\cdot\bfeta_0) q(0)
        + \int_{\Gamma}q(0)\bfeta_0\cdot\bn
        - \alpha_{BJ}\sum_{i = 1}^{d - 1} \int_{\Gamma}(\bfeta_0 \cdot\bftau_i)(\bv - \bfphi)(0)\cdot\bftau_i.
    \end{aligned}
\end{equation}
Here, by $\langle \cdot, \cdot \rangle_{(0, T) \times \Omega_F}$ we denote the duality pairing between $H^{-1}(0, T; L^2(\Omega_F))$ and $H^{1}_0(0, T; L^2(\Omega_F))$.
\end{definition}

\begin{remark}[trace of $\partial_t \bfeta$]
	Notice that, to define the weak solution, we integrated by parts the interface terms containing $\partial_t\bfeta$. This allows us to avoid having to define the trace of $\partial_t\bfeta$ on $\Gamma$; see~\cite[Section 2.2]{Avalos2024} for more details.
\end{remark}
To simplify notation, we introduce
\begin{equation*}
    \mathcal{H} \coloneqq \left\{ \bv \in L^2(\Omega_F)^d \; : \; \nabla \cdot \bv = 0 \; {\rm in } \; \Omega_F, \; \bv \cdot \bn = 0 \; {\rm on} \; \Gamma_F^1 \right\},
\end{equation*}
and the energy seminorm induced by the elasticity tensor~\eqref{eq:elasticity-tensor} for the Biot displacement, given by
\begin{equation}\label{eq:disp-energy-norm}
    \| \bfphi \|_E^2 \coloneqq 
    \|\mathbb{C}^{\frac12}\bD(\bfphi)\|^2_{L^2(\Omega_B)^{d\times d}} =
    2\mu_B \| \bD(\bfphi) \|^2_{L^2(\Omega_B)^{d\times d}} + \lambda_B \| \nabla \cdot \bfphi \|^2_{L^2(\Omega_B)}. 
\end{equation}
Notice that this is an equivalent norm on $\mathcal{V}_B$.

The theory of weak solutions for the Stokes--Biot system has been developed only recently. 

\begin{theorem}[{well-posedness of the sharp formulation~\cite[Theorem 2.3]{Avalos2024}}]\label{thm:WellPosedSharpInterface}
For every set of initial conditions
\begin{equation*}
    (\bu_0, \bfeta_0, \bfxi_0, p_{0}) \in \mathcal{H} \times \mathcal{V}_B \times L^2(\Omega_B)^d \times L^2(\Omega_B),
\end{equation*}
and right hand sides
\begin{equation*}
    (\bF_F, \bF_B, g) \in  L^2( (0, T) \times \Omega_F )^d \times L^2((0, T) \times \Omega_B)^d \times L^2((0, T) \times \Omega_B),
\end{equation*}
there is a weak solution to our problem in the sense of Definition~\ref{def:WS}. Moreover, this solution satisfies the following so-called energy inequality
    \begin{align*}
        \frac{1}{2}\left(
        \rho_F\|\bu\|^2_{L^2(\Omega_F)^d}
        + \rho_B \|\partial_t \bfeta\|^2_{L^2(\Omega_B)^d}
        + c_0\|p\|^2_{L^2(\Omega_B)}
        + \| \bfeta \|_E^2
        \right)(t)
        \\
        + \int_0^t \left ( 2\mu_F\|\bD(\bu)\|^2_{L^2(\Omega_F)^{{d \times d}}}
        + \alpha_{BJ} \sum_{i = 1}^{d - 1} \| (\bu - \partial_t \bfeta) \cdot\bftau_i \|^2_{L^2(\Gamma)}
        + \|\boldsymbol \kappa^{\frac12} \nabla p\|^2_{L^2(\Omega_B)^{{d}}}\right )
        \\
        \lesssim  \int_0^t \left( \rho_F\|\bF_F\|_{L^2(\Omega_F)^d}^2
        + \rho_B\|\bF_B\|_{L^2(\Omega_B)^d}^2
        + \|g\|^2_{L^2(\Omega_B)} \right) 
        \\
        + \rho_F \|\bu_0\|^2_{L^2(\Omega_F)^d}
        + \rho_B \|\bfeta_0\|_{H^1(\Omega_B)^d}^2
        + \rho_B \|\bfxi_0\|_{L^2(\Omega_B)^d}^2
        + c_0 \|p_{0}\|_{L^2(\Omega_B)}^2,
    \end{align*}
for a.e.~$t\in (0, T)$.
\end{theorem}

\begin{remark}[geometric configuration]
	Strictly speaking, in~\cite{Avalos2024} a specific geometric configuration with two stacked boxes is considered. However, as the authors remark, the existence result is valid for more general geometries and this will be presented in their forthcoming work. 
\end{remark}

\subsection{The diffuse interface formulation}\label{sub:PF}

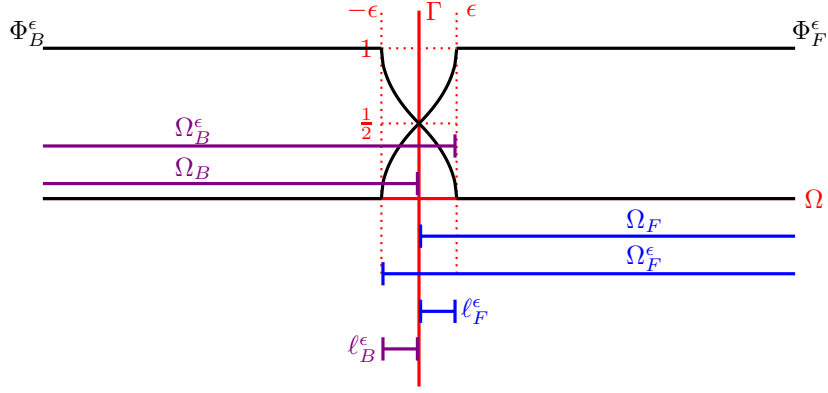
\begin{figure}
    \begin{center}
        \begin{tikzpicture}
            \draw[red, very thick] (-5, 0) -- (5, 0)
            node at (5.25, 0) {$\Omega$};
            \draw[red, very thick] (0, -2.5) -- (0, 2.5)
            node at (0.2, 2.5) {$\Gamma$};
            \draw[red, thick, dotted] (0.5, -1) -- (0.5, 2.5)
            node at (0.7, 2.5) {$\eps$};
            \draw[red, thick, dotted] (-0.5, -1) -- (-0.5, 2.5)
            node at (-0.75, 2.5) {$-\eps$};
            \draw[red, thick, dotted] (-0.5, 1) --(0.5, 1)
            node at (-0.7, 1) {$\frac1{2}$};
            \draw[red, thick, dotted] (-0.5, 2) --(0.5, 2)
            node at (-0.7, 2) {$1$};

            \draw[black, very thick] (-5, 0) -- (-0.5, 0);
            \draw[black, very thick] (0.5, 2) -- (5, 2)
            node at (5.2, 2.2) {$\Phi_F^\eps$};
            \draw[black, very thick, smooth, domain = -0.5:0] plot ({\x, sqrt(2.*\x + 1 ) });
            \draw[black, very thick, smooth, domain = 0:0.5] plot ({\x, 2 - sqrt(1 - 2.*\x ) });
      
            \draw[black, very thick] (-5, 2) -- (-0.5, 2)
            node at (-5.2, 2.2) {$\Phi_B^\eps$};;
            \draw[black, very thick] (0.5, 0) -- (5, 0);
            \draw[black, very thick, smooth, domain = 0:0.5] plot ({\x, sqrt(-2.*\x + 1 ) });
            \draw[black, very thick, smooth, domain = -0.5:0] plot ({\x, 2-sqrt(1 + 2.*\x ) });
      
            \draw[blue, very thick,|-] (0, -0.5) -- (5, -0.5) 
            node at (3, -0.3) {$\Omega_F$};
            \draw[blue, very thick,|-] (-0.5, -1) -- (5, -1) 
            node at (3, -0.8) {$\Omega_F^\eps$};
        
            \draw[blue, very thick, |-|] (0, -1.5) -- (0.5, -1.5)
            node at (0.75, -1.5) {$\ell_F^\eps$};
            \draw[violet, very thick, |-|] (-0.5, -2) -- (0, -2)
            node at (-0.75, -2) {$\ell_B^\eps$};
      
            \draw[violet, very thick, -|] (-5, 0.2) -- (0, 0.2)
            node at (-3, 0.4) {$\Omega_B$};
            \draw[violet, very thick, -|] (-5, 0.7) -- (0.5, 0.7)
            node at (-3, 0.9) {$\Omega_B^\eps$};
        \end{tikzpicture}    
    \end{center}
    \caption{Graphical representation of the diffuse interface approach (inspired by~\cite{Bukac2023}), in which $\Omega_F$ and $\Omega_F^\eps$ are, respectively, the sharp and diffuse fluid domains. Similarly $\Omega_B$ and $\Omega_B^\eps$ are the sharp and diffuse poroelastic domains, and $\Phi_F^\eps$ and $\Phi_B^\eps$ are the distance functions that define the diffuse fluid and poroelastic domains, respectively. Finally, $\ell_F^\eps$ and $\ell_B^\eps$ are transitional layers.}\label{fig:DomainsAndLayers}
\end{figure}

We now briefly describe our diffuse interface approach; this is explained in detail in~\cite[Sec.~3.3]{Bukac2023}, but we summarise it for self-containment. 

For $\eps > 0$ and $\beta\in (0, 1)$, we define
\begin{equation}\label{eq:defofPhiPsi}
    \cS(t) \coloneqq \begin{dcases}
         -1, & t \leq -1, \\
         (t + 1)^\beta -1, & t \in (-1, 0], \\
         1 - (1 - t)^\beta, & t \in (0, 1], \\
         1, & t > 1,
    \end{dcases}
    \qquad\text{ and }\qquad
    \renewcommand{\arraystretch}{2.0}
    \begin{array}{l}
        {\displaystyle 
            \Phi_F^\eps(x) \coloneqq \frac12\left( 1 + \mathcal{S}\left( \frac{\sdist_\Gamma(x)}\eps\right) \right),
        }
        \\
        \Phi_B^\eps(x) \coloneqq 1 - \Phi_F^\eps(x),
    \end{array}
\end{equation}
where $\sdist_\Gamma$ denotes the signed distance function to $\Gamma$ which is positive on $\Omega_F$. 
We then have, for $\eps$ sufficiently small, that $\Phi_F^{\eps}\approx 1$ in $\Omega_F$, $\Phi_F^{\eps}\approx 0$ in $\Omega_B$, and $\Phi_F^{\epsilon}$ transitions between these two values on a ``diffuse'' layer of width $\mathcal{O}(\epsilon)$. A similar reasoning applies to $\Phi_B^\eps$. We introduce, for $i\in\{F, B\}$, the following domains, see Figure~\ref{fig:DomainsAndLayers},
\begin{equation}\label{eq:DefOfDiffuseDomains}
    \begin{aligned}
        \Omega^{\eps}_i \coloneqq \{ x \in \Omega: \Phi_i^{\epsilon}(x) > 0\}, \quad
        \ell_i^\eps \coloneqq \left\{x \in \Omega_i^\eps : \Phi_i^\eps(x) \in (\tfrac1{2}, 1) \right\}, \quad  
        \ell^\eps \coloneqq \ell_F^\eps \cup \Gamma \cup \ell_B^\eps.
    \end{aligned}
\end{equation}
Since 
$\Omega_F \subset \Omega_F^\eps$ and $\Omega_B \subset \Omega_B^\eps$, these are diffuse versions of our fluid and poroelastic domains, respectively. The domains $\ell_F^\eps$ and $\ell_B^\eps$ are transitional layers, and for $\eps$ sufficiently small, we have
\begin{equation}\label{eq:MeasureOfLayer}
    |\ell^\eps| = |\ell_F^\eps| + |\ell_B^\eps| \lesssim \eps \mathcal{H}_{d - 1}(\Gamma),
\end{equation}
where the implied constant is independent of $\eps$ and $\mathcal{H}_{d - 1}(\Gamma)$.
The reason for this particular construction of a phase field function is the following result.

\begin{proposition}[$\Phi_F^\eps \in A_2$]\label{prop:PhiisA2}
    The function $\Phi_F^\eps$, when restricted to $\Omega_F^\eps$, is such that $\Phi_F^\eps \in A_2$. Similarly, the restriction of $\Phi_B^\eps$ to $\Omega_B^\eps$ defines an $A_2$ weight.
    More importantly, for $\eps$ sufficiently small, and $x \in \ell_B^\eps$, we have
    \begin{equation}\label{eq:PhiisDist}
        \Phi_F^\eps(x) = \frac12 \left( \frac{ \dist_\Gamma(x) + \eps}\eps \right)^\beta = \frac1{2\eps^\beta} \dist_{\Gamma^\eps}(x)^\beta,
    \end{equation}
    with $\Gamma^\eps = \Omega \cap \partial\Omega_F^\eps$.
    Therefore, $[\Phi_F^\eps]_{A_2}$ and $[\Phi_B^\eps]_{A_2}$ are independent of $\eps$ .
\end{proposition}
\begin{proof}
We refer the reader to~\cite[Prop.~3.3]{Bukac2023} for a proof that $\Phi_F^\eps, \Phi_B^\eps \in A_2$, and the claimed equality on $\ell_B^\eps$.
Then clearly $[\Phi_F^\eps]_{A_2} = [\dist_{\Gamma^\eps}^\beta]_{A_2}$. To conclude we should observe that although $\Gamma^\eps$ changes with $\eps$, it is just by translation in the normal direction. Thus, for $\eps$ sufficiently small, we have $[\dist_{\Gamma^\eps}^\beta]_{A_2} = [\dist_{\Gamma}^\beta]_{A_2}$.
\end{proof}

Let us now, for $k\in\N$, define the following spaces:
\begin{equation*}
    \begin{aligned}
        \mathcal{V}_i^{k, \eps} &\coloneqq \left\{ \bv \in H^k(\Omega^{\eps}_i, \Phi_i^{\eps})^d \; : \; \bv = \boldsymbol0 \; \textrm{on} \; \Gamma_i^1 \right\},  
        &\mathcal{V}_i^{\eps} &\coloneqq \mathcal{V}_i^{1, \eps}, 
        \qquad i\in\{F, B\},
        \\
        \mathcal{M}^{k, \eps} &\coloneqq H^k(\Omega_F^\eps, \Phi_F^\eps), &\mathcal{M}^\eps &\coloneqq \mathcal{M}^{0, \eps}, 
        \\
        \mathcal{Q}^{k, \eps} &\coloneqq H^k(\Omega_F^\eps, 1/\Phi_F^\eps), &\mathcal{Q}^\eps &\coloneqq \mathcal{Q}^{0, \eps},
        \\
        \mathcal{X}^{k, \eps} &\coloneqq  H^k(\Omega_B^{\eps}, \Phi_B^{\eps}),
        & \mathcal{X}^{\eps} &\coloneqq \mathcal{X}^{1, \eps}, 
    \end{aligned}
\end{equation*}
Owing to~\eqref{eq:PhiisDist}, all the results regarding weighted Sobolev spaces we mentioned in Section~\ref{sec:Notation} apply. In particular, the constants in the Poincar\'e, Bogovski{\u\i}, and Korn inequalities~\eqref{eq:WeightedPoincare}--\eqref{eq:weightedKorn} are 
now taken as 
independent of $\eps$, 
with the understanding that $\eps$ is small enough that Proposition~\ref{prop:PhiisA2} applies. 
The spaces $\mathcal{V}^{\eps}_F$ and $\mathcal{M}^\eps$ will be associated with the fluid velocity and pressure, respectively; the space $\mathcal{Q}^\eps$ is auxiliary. The spaces $\mathcal{V}^{\eps}_B$ and $\mathcal{X}^{\eps}$ will be associated with the poroelastic displacement and Biot pressure, respectively.

To handle integrals at the diffuse interface we shall need the following result,
a diffuse analogue of the trace inequality $H^1(\Omega)^d\to L^2(\Gamma)$ for $\Gamma$ sufficiently regular.

\begin{lemma}[{diffuse trace inequality~\cite[Lemma 3.6]{Bukac2023}}]\label{lem:nablaPhi}
    Let $\epsilon_0$ be sufficiently small. Then 
    for $0 < \epsilon < \epsilon_0$, 
    the restriction 
    \begin{equation*}
        \bv\mapsto\bv|_{\ell^\eps}
        \text{ is continuous as a map }\mathcal{V}^\eps_i\to L^2\left(\ell^\eps, \tfrac{1}{2\eps}|\nabla\dist_\Gamma|\right)^d,
    \end{equation*}
    i.e.~there exists a constant $C_{\rm tr} > 0$ 
    independent of $\eps$
    such that for 
    $\bv \in \mathcal{V}_i^{\eps}$, we have
    \begin{equation*}
        \frac1{2\eps} \int_{\ell^\eps} |\bv|^2 |\nabla \dist_\Gamma | \leq C_{\rm tr} \| \bv \|^2_{\mathcal{V}_i^\eps}, \qquad i\in\{F, B\}.
    \end{equation*}
    An analogous estimate holds for $\psi \in \mathcal{X}^\eps$.
\end{lemma}

We are now ready to introduce the notion of weak solution for the diffuse interface problem. We shall assume that there is $\eps_0 > 0$ such that, for every $\eps \in (0, \eps_0]$, we have $\bu_0^\eps$, $\bfeta_0^{\eps}$, $\bfxi_0^{\eps}$, $p_{0}^\eps$, $\bF_F^\eps$, $\bF_B^\eps$, $g^\eps$, and $\boldsymbol{\kappa}_\eps$, which are suitable extensions of $\bu_0$, $\bfeta_0$, $\bfxi_0$, $p_{0}$, $\bF_F$, $\bF_B$, $g$, and $\boldsymbol\kappa$, to $\Omega_F^\eps$ and $\Omega_B^\eps$, respectively.
In addition, we assume that $\boldsymbol{\kappa}_\eps$ obeys the same spectral bounds as $\boldsymbol{\kappa}$.

Define $\mathcal{A}_{\eps}: (\mathcal{V}^\eps_ F \times \mathcal{V}^\eps_B \times \mathcal{X}^\eps)^2\to\mathbb{R}$ via
\begin{equation}\label{eq:A_eps}
    \begin{split}
        \cA_\eps((\bu, \bfeta, p), (\bv, \bfphi, q))
        \coloneqq 2\mu_F\int_{\Omega_F^\eps}\bD(\bu):\bD(\bv)\Phi_F^\eps
        &
        + \int_{\Omega_B^\eps}\bfsigma_E(\bfeta):\bD(\bfphi) \Phi_B^\eps
        + \int_{\Omega_B^\eps} \boldsymbol\kappa_\eps\nabla p \cdot \nabla q \Phi_B^\eps
        \\ &
        + \frac1{2\eps}\int_{\ell^\eps} q \bu \cdot \nabla \dist_\Gamma
        - \frac1{2\eps} \int_{\ell^\eps} p \bv \cdot \nabla \dist_\Gamma.
    \end{split}
\end{equation}
For almost every $t \in (0, T)$, the linear form $\mathcal{F}_\eps(t): \mathcal{V}_F^{\epsilon}\times\mathcal{V}_B^\eps\times\mathcal{X}^{\epsilon} \to \R$ is defined as
\begin{equation}\label{eq:F_eps}
    \langle\mathcal{F}_\eps(t), (\bv, \bfphi, q)\rangle \coloneqq \rho_F \int_{\Omega_F^\eps}\bF_F^\eps(t)\cdot\bv\Phi_F^\eps + \rho_B\int_{\Omega^\eps_B}\bF_B^\eps(t)\cdot\bfphi\Phi_B^\eps + \int_{\Omega^\eps_B}g^\eps(t)q\Phi_B^\eps.
\end{equation}
Our notion of solution reads as follows.

\begin{definition}[weak solution of the diffuse interface problem]\label{def:SBFP}
We say that the tuple $(\bu^{\eps}, \pi^{\eps}, \bfeta^\eps,  p^\eps)$ is a diffuse interface weak solution to our problem if
\begin{equation*}
    \begin{aligned}
        \bu^{\eps} &\in L^\infty(0, T; L^2(\Omega_F^\eps, \Phi_F^\eps)^d) \cap L^2(0, T; \mathcal{V}_F^{\eps}), 
        & \bfeta^{\eps}&\in L^{\infty}(0, T; \mathcal{V}_B^{\eps})\cap W^{1, \infty}(0, T; L^2(\Omega_B^\eps, \Phi_B^\eps)^d), 
        \\
        \pi^{\eps} &\in H^{-1}(0, T; \mathcal{M}^\eps), 
        & p^{\eps} &\in L^\infty(0, T; L^2(\Omega_B^\eps, \Phi_B^\eps )) \cap L^2(0, T; \mathcal{X}^{\eps}),
    \end{aligned}
\end{equation*}
$\bfeta^{\eps}(0) = \bfeta_0^{\eps}$, and in addition, for every tuple $(\bv, \zeta, \bfphi, q)$ such that
\begin{equation*}
    \begin{aligned}
        \bv &\in C^1_c([0, T)\times \overline{\Omega_F^\eps} \setminus \Gamma_F^1)^d,
        & \bfphi &\in C^1_c([0, T)\times \overline{\Omega_B^\eps} \setminus \Gamma_B^1)^d, \\
        \zeta &\in C^1_c([0, T)\times \Omega_F^\eps ),
        & q &\in C^1_c([0, T)\times \overline{\Omega_B^\eps} ),
    \end{aligned}
\end{equation*}
the following equality is satisfied:
    \begin{align}\label{eq:SBFP_weak}
        &- \rho_F\int_0^T \int_{\Omega_F^\eps}\bu^{\eps}\cdot \partial_t\bv \Phi_F^\eps
        - \rho_B \int_0^T\int_{\Omega_B^\eps}\partial_{t} \bfeta^{\eps}\cdot \partial_t\bfphi \Phi_B^\eps
        - \int_0^T\int_{\Omega_B^\eps} (c_0p^{\eps} + \alpha\nabla\cdot\bfeta^{\eps}) \partial_t q \Phi_B^\eps 
        \\ &
        + \int_0^T \cA_\eps((\bu^\eps, \bfeta^\eps, p^\eps), (\bv, \bfphi, q))
        - \langle \pi^{\eps}, \nabla \cdot \bv \rangle_{(0, T)\times\Omega_F^\eps, \Phi_F^\eps} 
        + \int_0^T\int_{\Omega_F^\eps} (\nabla \cdot \bu^{\eps}) \zeta \Phi_F^\eps
        \\ &
        - \alpha\int_0^T\int_{\Omega_B^\eps} p^{\eps} \nabla \cdot \bfphi \Phi_B^\eps
        - \frac1{2\eps}\int_0^T\int_{\ell^\eps}q \partial_t \bfeta^{\eps} \cdot \nabla \dist_\Gamma
        + \frac1{2\eps} \int_0^T\int_{\ell^\eps} p^{\eps} \bfphi \cdot \nabla \dist_\Gamma 
        \\ &
        + \frac{\alpha_{BJ}}{2\eps} \sum_{i = 1}^{d - 1} \int_0^T\int_{\ell^\eps}((\bu^{\eps} - \partial_t \bfeta^{\eps}) \cdot \tilde{\bftau}_i) ((\bv - \bfphi) \cdot \tilde{\bftau}_i) |\nabla\dist_\Gamma|
        = 
        \int_0^T \langle\mathcal{F}_\eps, (\bv, \bfphi, q)\rangle
        \\ & 
        + \rho_F\int_{\Omega_F^\eps}\bu_0^{\eps}\cdot \bv(0)\Phi_F^\eps
        + \rho_B\int_{\Omega_B^\eps}\bfxi_0^{\eps}\cdot \bfphi(0)\Phi_B^\eps
        + \int_{\Omega_B^\eps}(c_0p_{0}^{\eps} + \nabla\cdot\bfeta_0^{\eps})q(0)\Phi_B^\eps.
    \end{align}
Above, by $\langle \cdot , \cdot \rangle_{(0, T)\times\Omega_F^\eps, \Phi_F^\eps}$, we denote the duality pairing between $H^{-1}(0, T; L^2(\Omega_F^\eps, \Phi_F^\eps))$ and $H^1_0(0, T;L^2(\Omega_F^\eps, 1/\Phi_F^\eps))$.
Diffuse, normalised approximations $\{\tilde{\bftau}_i\}_i$ to the tangent vector fields $\{\bftau_i\}_i$ may be obtained directly from the phase field function using 
formulae from~\cite{Stoter2017}. 
\end{definition}

Notice that the integrals 
over $\ell^\eps$ in the diffuse weak formulation
correspond to the interface integrals in the sharp interface one. 
We can invoke Lemma~\ref{lem:nablaPhi} to show that the integrals involving the unknowns $\bu^\eps$ and $p^\eps$ are well-defined. However, since we only have that $\partial_t\bfeta^\eps(t) \in L^2(\Omega_B^{\eps},\Phi_B^{\eps})^d$, Lemma~\ref{lem:nablaPhi} does not apply. Notice that a similar issue was faced in the sharp interface formulation. This was circumvented by formally integrating by parts in time and thus applying the time derivative to the test function. In the diffuse interface formulation we use the fact that, since $\Phi_B^\eps \in A_2$, we have $L^2(\Omega_B^{\eps}, \Phi_B^{\eps})\hookrightarrow L^s(\Omega_B^{\eps})$ for some $s > 1$; see for instance~\cite[Corollary 3.3]{Adimurthi2021}. 
Therefore we have, for instance,
\begin{equation*}
    \left |\frac1{2\eps} \int_0^T \int_{\ell^\eps} q \partial_t \bfeta^\eps \cdot \nabla \dist_\Gamma\right |
    \leq \frac{T}{2\eps}\| \nabla \dist_\Gamma\|_{L^{\infty}(\Omega_B^{\eps})^d}\|\partial_t \bfeta^\eps \|_{L^\infty(0, T;~L^{1}(\Omega_B^{\eps})^d)} \| q \|_{L^\infty((0, T)\times\Omega_B^\eps)}.
\end{equation*}
Hence, for every $\eps \in (0, \eps_0]$ this integral is well-defined. The same reasoning can be applied to the other integral arising from the Beavers--Joseph--Saffman condition. Of course, the above estimate is not uniform in the diffuse interface parameter $\eps$, but this fact will not cause any issues because in the subsequent analysis we will use a suitable cancellation property to obtain estimates which are uniform in $\epsilon$.
We note also that in all spatial integrals in the diffuse weak form can equivalently be written over the entire domain $\Omega$. In that case, the formulation is more suitable for numerical implementation, as described in Section~\ref{sec:numerics}. 

Our main goal in this paper is to study the convergence of suitable discretisations of~\eqref{eq:SBFP_weak} to solutions of our problem in the sharp interface sense, i.e.~according to Definition~\ref{def:WS}. This will be done in two steps. In the first step, see Section~\ref{sec:numerical-error}, we fix $\epsilon$ and the phase field functions $\Phi_F^\eps, \Phi_B^\eps$, and prove error estimates for a finite element approximation of the diffuse interface formulation~\eqref{eq:SBFP_weak}. In the second step, see Section~\ref{sec:ModellingError}, we analyse the convergence of the continuous diffuse interface  formulation to the continuous sharp interface formulation.

\subsection{Well-posedness}\label{sub:DiffuseDomainWellPosed}

Let us prove the well-posedness for the diffuse interface formulation. Owing to all the existing theory regarding weighted Sobolev spaces, the proof is similar to the well-posedness proof for the diffuse interface Stokes--Darcy system~\cite[Theorem 3.9]{Bukac2023}. Therefore, here we just outline the main steps of the proof. We define the function spaces
\begin{equation*}
    \mathcal{H}^\eps_{\rm div} \coloneqq \left\{ \bv \in L^2(\Omega_F^\eps, \Phi_F^\eps)^d\; : \; \nabla\cdot\bv = 0 \; \textrm{in} \; \Omega_F^\eps \;, \; \bv\cdot \bn = 0 \; \textrm{on} \; \Gamma_F^1 \right\},
    \qquad
    \mathcal{V}^{\epsilon}_{\rm div} \coloneqq \left\{\bv \in \mathcal{V}_F^{\eps} : \nabla\cdot\bv = 0 \right\}.
\end{equation*}
Notice that, by definition, $\bu^\eps \in L^2(0, T; \mathcal{V}^{\epsilon}_{\rm div})$. 
In analogy to~\eqref{eq:disp-energy-norm}, we also define the diffuse energy seminorm
\begin{equation}\label{eq:diffuse-energy-norm}
    \|\bfphi\|^2_{E, \eps} \coloneqq 2\mu_B\|\bD(\bfphi)\|^2_{L^2(\Omega^\eps_B, \Phi^\eps_B)^{d\times d}} + \lambda_B\|\nabla\cdot\bfphi\|^2_{L^2(\Omega^\eps_B, \Phi^\eps_B)}.
\end{equation}
This is an equivalent norm on $\mathcal{V}^\eps_B$.
We can now prove well-posedness.

\begin{theorem}[well-posedness of the diffuse formulation]\label{thm:WellPosedness}
There is $\eps_0 > 0$ such that, for every $\eps \in (0, \eps_0]$, problem~\eqref{eq:SBFP_weak} is well-posed in the following sense. For every set of initial conditions
\begin{equation*}
    (\bu_0^\eps, \bfeta^{\eps}_0, \bfxi_0^\eps, p_0^\eps )
    \in \mathcal{H}^\eps_{\rm div} \times \mathcal{V}_{B}^{\eps} \times L^2(\Omega_B^\eps, \Phi_B^\eps )^d \times L^2(\Omega_B^\eps, \Phi_B^\eps ), 
\end{equation*}
and right hand sides
\begin{equation*}
    (\bF_F^\eps, \bF_B^\eps, g^\eps) \in
    L^2(0,T; L^2(\Omega_F^\eps, \Phi_F^\eps)^d) \times 
    L^2(0,T; L^2(\Omega_B^\eps,\Phi_B^\eps)^d) \times 
    L^2(0,T; L^2(\Omega_B^\eps,\Phi_B^\eps)),
\end{equation*}
the Stokes--Biot diffuse interface problem has a unique weak solution in sense of Definition~\ref{def:SBFP}. Moreover, this solution satisfies the following energy inequality
    \begin{align}\label{eq:EE}
        &\frac{1}{2}\left(
    	\rho_F\|\bu^\eps\|^2_{L^2(\Omega_F^\eps, \Phi_F^{\eps})^d}
        + \rho_B\|\partial_t \bfeta^\eps\|^2_{L^2(\Omega_B^\eps, \Phi_B^{\eps})^d}
        + c_0\|p^\eps\|^2_{L^2(\Omega_B^\eps, \Phi_B^\eps)}
        + \|\bfeta^\eps\|^2_{E, \eps}
        \right)(t)
        \\
	    + \int_0^t &\left (
        2\mu_F\|\bD(\bu^\eps)\|^2_{L^2(\Omega_F^\eps, \Phi_F^{\eps})^{{d \times d}}}
        + \alpha_{BJ} \sum_{i = 1}^{d - 1}\|(\bu^{\eps} - \partial_t \bfeta^{\eps}) \cdot\tilde\bftau_i\|^2_{L^2\left(\ell^\eps, \frac{1}{2\eps}|\nabla\dist_\Gamma|\right)}
        + \|\boldsymbol{\kappa}_\eps^{\frac12} \nabla p^\eps\|^2_{L^2(\Omega_B^\eps, \Phi_B^{\eps})^d}\right )
        \\
        \lesssim 
        \int_0^t&\left( \rho_F \|\bF_F^{\eps}\|_{L^2(\Omega_F^\eps, \Phi_F^{\eps})^d}^2
        + \rho_B\|\bF_B^\eps\|_{L^2(\Omega_B^\eps, \Phi_B^\eps)^d}^2
        + \|g^\eps\|_{L^2(\Omega_B^\eps, \Phi_B^\eps)}^2 \right) 
        \\
        &
        + \rho_F\|\bu_0^\eps\|^2_{L^2(\Omega_F^\eps, \Phi_F^{\eps})^d}
        + \rho_B\|\bfeta_0^\eps\|_{H^1(\Omega_B^\eps, \Phi_B^{\eps})^d}^2 
        + \rho_B \|\bfxi_0^\eps\|_{L^2(\Omega_B^\eps, \Phi_B^{\eps})^d}^2
        + c_0 \|p_{0}^\eps\|_{L^2(\Omega_B^\eps, \Phi_B^{\eps})}^2,
    \end{align}
for a.e.~$t\in (0, T)$.
\end{theorem}
\begin{proof}
The proof is rather standard and so here we just outline the main steps and emphasise several points characteristic to the diffuse interface Stokes--Biot problem. We use a Galerkin method to construct the solution, where the finite dimensional fluid velocity spaces are taken to be solenoidal, i.e.~they are subspaces of $\mathcal{H}^{\eps}_{\rm div}$. This allows us to ignore the fluid pressure at this stage. Let us denote by $\{\bu^{\eps}_n, \bfeta^{\eps}_n, p^{\eps}_n \}_{n > 0}$ the family of Galerkin approximations. By taking $(\bv, \bfphi, q) = (\bu^{\eps}_n, \partial_t\bfeta^{\eps}_n, p^{\eps}_n)$ in the weak formulation~\eqref{eq:SBFP_weak} we obtain, after an application of Gr\"onwall's lemma, estimate~\eqref{eq:EE} for the Galerkin approximations. By the weighted Poincar\'e inequality~\eqref{eq:WeightedPoincare} we conclude that 
\begin{itemize}
    \item $\{ \bu^{\eps}_n \}_{n > 0}$ is uniformly bounded in $L^\infty(0, T; \mathcal{H}^{\eps}_{\rm div}) \cap L^2(0, T; \mathcal{V}_F^{\eps})$, 
    \item $\{ \bfeta^{\eps}_n \}_{n > 0}$ is uniformly bounded in $L^{\infty}(0, T; \mathcal{V}_F^{\eps})\cap W^{1, \infty}(0, T; L^2(\Omega_B^\eps, \Phi_B^\eps)^d)$.
\end{itemize}
However, since we do not have a Poincar\'e-type inequality in $\mathcal{X}^{\eps}$, we have to estimate the mean value of each $p^{\eps}_n$. Let $\gamma = \gamma(t)\in C^1_c([0, T))$ be such that $\gamma'(t) = -1$. Let us take $(\bv, \bfphi, q) = (0, 0, \gamma)$ as test functions in~\eqref{eq:SBFP_weak} to obtain
\begin{equation*}
    \begin{aligned}
        c_0\int_0^T\int_{\Omega_B^{\eps}}p^{\eps}_n\Phi_B^{\eps}
        =& - \alpha\int_0^T\int_{\Omega_B^{\eps}}\nabla\cdot\bfeta_n^{\eps}\Phi_B^{\eps}
        - \frac{1}{2\eps}\int_0^T\gamma \int_{\ell^{\eps}}(\bu^{\eps}_n - \partial_t\bfeta^{\eps}_n)\cdot\nabla\dist_{\Gamma} \\
        &+ \gamma(0)\int_{\Omega_B^\eps} \left( c_0 p_{0, n}^{\eps} + \alpha \nabla\cdot\bfeta_{0, n}^{\eps} \right)\Phi_B^\eps.
    \end{aligned}
\end{equation*}
Therefore we have:
\begin{equation*}
    \begin{aligned}
        \left|\int_0^T\int_{\Omega_B^{\eps}}p^{\eps}_n\Phi_B^{\eps}\right| &\leq 
        C\int_0^T\left(
        \|\nabla \cdot \bfeta^{\eps}_n\|_{L^1(\Omega_B^\eps,\Phi_B^{\eps})}
        + \frac1{2\eps} \|\bu^{\eps}_n\|_{L^1(\ell^\eps)^d} + \frac{1}{2\eps}\|\partial_t\bfeta^{\eps}_n\|_{L^1(\ell^{\eps})^d} 
        \right) \\
        &+ C\|p_0^\eps \|_{L^1(\Omega_B^\eps, \Phi_B^{\eps})} + C\|\nabla \cdot \bfeta_0^\eps \|_{L^1(\Omega_B^\eps, \Phi_B^\eps)} \\
        &\leq C\int_0^T\left (\|\bfeta^\eps_n \|_{\mathcal{V}_{B}^{\eps}} + \frac1\eps \|\bu^\eps_n\|_{L^2(\Omega_F^{\eps}, \Phi_F^\eps)^d}
        + \frac{1}{\eps} \|\partial_t\bfeta^\eps_n\|_{L^2(\Omega_B^\eps, \Phi_B^{\eps})^d}  \right ) + C_0,
    \end{aligned}
\end{equation*}
where $C_0$ depends only on the initial data. Here, again, we used the embeddings $L^2(\Omega_i^{\eps}, \Phi_i^{\eps})\hookrightarrow L^1(\Omega_i^{\eps})$, $i \in \{F, B\}$, and the fact that $\Phi_B^\eps$ is bounded and positive so that, for instance,
\begin{equation*}
    \| \nabla \cdot \bfeta_n^\eps \|_{L^1(\Omega_B^\eps, \Phi_B^\eps)} \leq \| \nabla \cdot \bfeta_n^\eps \|_{L^2(\Omega_B^\eps, \Phi_B^\eps)} \left( \int_{\Omega_B^\eps} \Phi_B^\eps \right)^{1/2} \leq C \| \nabla \cdot \bfeta_n^\eps \|_{L^2(\Omega_B^\eps, \Phi_B^\eps)}.
\end{equation*}
We combine the obtained inequality with~\eqref{eq:EE} to conclude (notice that at this point $\eps$ is fixed)
\begin{itemize}
    \item $\{ p^{\eps}_n \}_{n > 0}$ is uniformly bounded in $L^2(0, T; \mathcal{X}^{\eps})$.
\end{itemize}

Now, since the problem is linear, we can pass to the limit using weak and weak$^\star$ convergence in a standard manner. Let us illustrate how to pass to the limit in a boundary term characteristic to the diffuse interface Stokes--Biot coupling. Again, using embeddings of weighted Lebesgue spaces we have $\partial_t\bfeta^\eps_n \rightharpoonup \partial_t\bfeta^{\eps}$ weakly in $L^s(\Omega_B^{\eps})$ for some $s > 1$. Therefore
\begin{equation*}
    \frac{1}{2\eps}\int_{\ell^{\eps}} q \partial_t\bfeta^{\eps}_n\cdot\nabla\dist_{\Gamma}
    \to
    \frac{1}{2\eps}\int_{\ell^{\eps}} q \partial_t\bfeta^{\eps}\cdot\nabla\dist_{\Gamma}.
\end{equation*}

Finally, we must show the existence of the corresponding fluid pressure $\pi^{\eps}$. However, this construction is completely analogous to the corresponding construction in the diffuse interface Stokes--Darcy system, so we refer the reader to the proof of~\cite[Theorem 3.9]{Bukac2023}.
\end{proof}

\section{Discretisation}\label{sec:Discretisation}

Having obtained the well-posedness of our diffuse interface problem, in this section we proceed with its discretisation. We provide a numerical scheme and its error analysis.

\subsection{Time discretisation}\label{sub:TimeDiscr}

For time discretisation we will employ the backward Euler method. 
Let $N \in \mathbb{N}$ be the number of timesteps, and $\dt = T/N$ the timestep. We define, for $n = 0, \ldots, N$, the discrete times $t^n \coloneqq n \dt$. Let $X$ be a normed space. For a given function $W:[0, T] \to X$, we will compute sequences $W^{\dt} = \{W^n\}_{n = 0}^N$ so that $W^n \approx W(t^n)$. The discrete time derivative is defined as $ \dtee W^{n + 1} \coloneqq \frac{W^{n + 1} - W^n}{\Delta t}$,
and 
the second and third as $\dtt \coloneqq \dtee\circ \dtee$, 
$\mathrm{d}_{ttt} \coloneqq \dtee\circ\dtt$.
Over such sequences, we define the following norms
\begin{equation*}
    \| W^\dt \|_{L^2_\dt(0, T; X)}^2 \coloneqq \dt\sum_{n = 0}^{N-1} \| W^{n + 1} \|_X^2,
    \qquad
    \| W^\dt \|_{L^\infty_\dt(0, T; X)} \coloneqq \max_{0\leq n \leq N} \| W^n \|_X.
\end{equation*}

\subsection{Space discretisation}\label{sub:FEM}

To discretise in space, we use the finite element method. Since we have assumed that $\Om$ is a polytope, it can be meshed exactly. Let $\{\Th\}_{h > 0}$ be a family of conforming, simplicial triangulations of $\Om$ that is quasiuniform in the usual finite element sense~\cite{Ciarlet2002, Girault1986}. The parameter $h > 0$ denotes the 
characteristic
mesh size. Notice that we do not assume that the mesh is in any way aligned with $\Gamma$
or $\ell^\eps$.
We set
\begin{equation}\label{eq:submeshes}
    \mathcal{T}_{h, i}^\eps \coloneqq \bigcup\left\{ T \in \Th : T \cap \Omega_i^\eps \neq \emptyset \right\}, \qquad i\in\{F, B\}.
\end{equation}
We assume that, for each $h > 0$, we have at hand spaces
\begin{equation*}
    \widetilde{\mathcal{V}}_{h, i}^\eps \subset W^{1, \infty}(\Omega)^d, \ i \in \{F, B\}, \qquad
    \widetilde{\mathcal{Q}}_h^\eps \subset L^\infty(\Omega), \qquad
    \widetilde{\mathcal{X}}_h^\eps \subset W^{1, \infty}(\Omega),
\end{equation*}
which consist of piecewise polynomials subordinate to $\Th$. On the basis of these we define
    \begin{align*}
        \mathcal{V}_{h, i}^\eps &\coloneqq \left\{ \bv_h \in \widetilde{\mathcal{V}}_{h, i}^\eps : \supp \bv_h \subset \mathcal{T}_{h, i}^\eps, \ \bv_h = \boldsymbol0 \text{ on } \Gamma_i^1 \right\}
        \subset \mathcal{V}^\eps_i, \qquad i\in\{F, B\}, 
        \\
        \mathcal{Q}_h^\eps &\coloneqq \left\{ \zeta_h \in \widetilde{\mathcal{Q}}_{h}^\eps : \supp \zeta_h \subset \mathcal{T}_{h, F}^\eps \right\} 
        \subset \mathcal{Q}^\eps, 
        \\
        \mathcal{X}_h^\eps &\coloneqq \left\{ q_h \in \widetilde{\mathcal{X}}_{h}^\eps : \supp q_h \subset \mathcal{T}_{h, B}^\eps \right\} 
        \subset \mathcal{X}^\eps.
    \end{align*}
These spaces will be used to approximate the fluid velocity, poroelastic displacement, 
Biot 
pressure, and a quantity related to the fluid pressure, respectively. 
Note that our finite element spaces are defined with respect to submeshes $\mathcal{T}_{h,i}^\eps$~\eqref{eq:submeshes} which themselves depend on $\eps$, but we shall assume $\eps$ to be fixed for the purpose of this section.

We must require that 
a discrete analogue of~\eqref{eq:Bogovskii}, holds uniformly in $h$: 
the fluid velocity and pressure spaces are compatible in the following sense:
there is $b > 0$ such that for all $h > 0$,
\begin{equation}\label{eq:DiscrInfSupWeirdBCs}
    b\| \zeta_h \|_{\mathcal{Q}^\eps} \leq \sup_{\boldsymbol0 \neq \bv_h \in \mathcal{V}^\eps_{h, F}} \frac{ \int_{\Omega_F^\eps} \zeta_h \nabla \cdot \bv_h }{ \| \bv_h \|_{\mathcal{V}^\eps_F}}, \qquad \forall \zeta_h \in \mathcal{Q}_h^\eps.   
\end{equation}
We refer the reader to~\cite[Lemma 4.1]{Bukac2023} for a proof of this inequality for a wide class of classical velocity-pressure finite element pairs. We 
also 
recall that~\eqref{eq:DiscrInfSupWeirdBCs} is equivalent to the existence of a Fortin projection $P_h:\mathcal{V}^\eps_F\to\mathcal{V}^\eps_{h, F}$.

For error analysis, we take the standard approach of using discrete projections of the exact solutions as intermediate approximants, and then concluding the final error bounds by the triangle inequality.
These projections will be with respect to the bilinear forms induced by 
the elliptic and coupled Stokes-elastic subproblems of the overall Stokes--Biot problem, defined in the next lemma.
Since these subproblems are already weighted by the phase field functions $\Phi^\eps_F, \Phi^\eps_B$, the approximation properties of the projections with respect to weighted norms readily follow, so that there is no need for results on the weighted-norm approximation of the unweighted problems as in for example~\cite{Nochetto2016}.

\begin{lemma}[weighted Stokes-like and Ritz projections]\label{lem:weighted-ritz}
    Assuming~\eqref{eq:DiscrInfSupWeirdBCs}, there exist unique bounded projections
    \begin{equation*}
        \mathcal{P}_{FB}: \mathcal{V}^\eps_F\times\mathcal{Q}^\eps \times \mathcal{X}^\eps \to \mathcal{V}^\eps_{F, h}\times\mathcal{Q}^\eps_h \times \mathcal{X}^\eps_h, \quad
        \mathcal{P}_B:\mathcal{V}^\eps_B\to\mathcal{V}^\eps_{B, h}, 
    \end{equation*}
    such that 
    for each $ (\mathfrak{u}, \vartheta, \mathfrak{p}) \in \mathcal{V}^\eps_F\times\mathcal{Q}^\eps \times \mathcal{X}^\eps $ and all $(\bv_h, \zeta_h, q_h)\in\mathcal{V}^\eps_{h, F}\times\mathcal{Q}^\eps_h\times\mathcal{X}^\eps_h$,
    denoting $\mathcal{P}_{FB}(\mathfrak{u}, \vartheta, \mathfrak{p}) = (\underline{\mathfrak{u}}, \underline{\vartheta}, \underline{\mathfrak{p}})$,
    \begin{equation}\label{eq:weighted-stokes-proj}
        \begin{aligned}
            \int_{\Omega^\eps_B}\boldsymbol\kappa_\eps\nabla(\underline{\mathfrak{p}} - \mathfrak{p})\cdot\nabla q_h\Phi_B^\eps &+ \frac{1}{2\eps}\int_{\ell^\eps}q_h(\underline{\mathfrak{u}} - \mathfrak{u})\cdot\nabla\dist_\Gamma & &= 0,
            \\
            - \frac{1}{2\eps}\int_{\ell^\eps}(\underline{\mathfrak{p}} - \mathfrak{p})\bv_h\cdot\nabla\dist_\Gamma &+ \int_{\Omega^\eps_F}2\mu_F\bD(\underline{\mathfrak{u}} - \mathfrak{u}):\bD(\bv_h)\Phi_F^\eps & - \int_{\Omega^\eps_F} (\nabla\cdot\bv_h)(\underline{\vartheta} - \vartheta) &= 0,
            \\
            & - \int_{\Omega^\eps_F}(\nabla\cdot(\underline{\mathfrak{u}} - \mathfrak{u}))\zeta_h & &= 0,
        \end{aligned}
    \end{equation}
    and for each $\Upsilon\in\mathcal{V}^\eps_B$ and all $\bfphi_h\in\mathcal{V}^\eps_{h, B}$,
    \begin{equation}\label{eq:weighted-ritz-proj}
        \int_{\Omega^\eps_B}\bfsigma_E(\mathcal{P}_B(\Upsilon) - \Upsilon):\bD(\bfphi_h)\Phi^\eps_B = 0.
    \end{equation}
    Moreover,
    these projections optimally approximate their arguments in weighted norms, i.e.~there is $r \in \N$ such that whenever $k \leq r$, if $(\mathfrak{u}, \vartheta, \mathfrak{p})\in \mathcal{V}^{k + 1, \eps}_F\times\mathcal{Q}^{k, \eps} \times \mathcal{X}^{k + 1, \eps}$ and $\Upsilon\in\mathcal{V}^{k + 1, \eps}_B$, 
    then
    \begin{equation}\label{eq:projection-error-estimate}
        \begin{aligned}
        \|\underline{\mathfrak{u}} - \mathfrak{u}\|_{\mathcal{V}^\eps_F} + \|\underline{\vartheta} - \vartheta\|_{\mathcal{Q}^\eps} + \|\underline{\mathfrak{p}} - \mathfrak{p}\|_{\mathcal{X}^\eps} &\lesssim h^k\left(\|\mathfrak{u}\|_{\mathcal{V}^{k + 1, \eps}_F} + \|\vartheta\|_{\mathcal{Q}^{k, \eps}} + \| \mathfrak{p}\|_{\mathcal{X}^{k + 1, \eps}}\right),
        \\
        \|\mathcal{P}_B(\Upsilon) - \Upsilon\|_{\mathcal{V}^\eps_B} &\lesssim h^k\|\Upsilon\|_{\mathcal{V}^{k + 1, \eps}_B},
        \end{aligned}
    \end{equation}
    with implied constants independent of $h$.
\end{lemma}
\begin{proof}
    That each projection exists, is unique, and is bounded follows 
    from the well-posedness of each variational 
    problem~\eqref{eq:weighted-stokes-proj}--\eqref{eq:weighted-ritz-proj};
    this is clear for~\eqref{eq:weighted-ritz-proj},
    and for the problem~\eqref{eq:weighted-stokes-proj} follows from coercivity of the 
    bilinear 
    form
    \begin{equation*}
        (\mathcal{V}^\eps_F \times \mathcal{X}^\eps)^2\ni
        ((\mathfrak{u}, \mathfrak{p}), (\bv,q)) \mapsto \int_{\Omega^\eps_B}\boldsymbol\kappa_\eps\nabla\mathfrak{p}\cdot\nabla q\Phi_B^\eps + \frac{1}{2\eps}\int_{\ell^\eps}(q\mathfrak{u} - \mathfrak{p}\bv)\cdot\nabla\dist_\Gamma + \int_{\Omega^\eps_F}2\mu_F\bD(\mathfrak{u}):\bD(\bv)\Phi_F^\eps
    \end{equation*}
    and inf-sup compatibility of the pair $\mathcal{V}^\eps_{F, h}\times\mathcal{Q}^\eps_h$.
    Approximation properties of $\mathcal{P}_{FB}$ and $\mathcal{P}_B$ follow immediately from their stability, the fact that they are projections,
    and those of the corresponding finite element interpolants
    $I_{\mathcal{V}, i}:\mathcal{V}^\eps_i \to \mathcal{V}^\eps_{h, i}$, for $i\in\{F, B\},
    I_{\mathcal{X}} : \mathcal{X}^\eps \to \mathcal{X}^\eps_h$,
    and
    $I_{\mathcal{Q}} : \mathcal{Q}^\eps \to \mathcal{Q}_h^\eps$.
    Owing to the fact that $\Phi_F^\eps, \Phi_B^\eps \in A_2$, the work~\cite{Nochetto2016} has shown that standard piecewise polynomial finite element spaces 
    admit such operators, and has provided an explicit construction thereof.
\end{proof}

\subsection{The scheme}\label{sub:TheNumScheme}

We are now ready to present the scheme and discuss its basic properties. We begin by setting
\begin{equation}\label{eq:discrete-initial-data}
    \begin{aligned}
        (\bu_h^0, p_h^0) &= (\mathcal{P}_{FB}(\bu_0^\eps, 0, p^\eps_0)_1, \mathcal{P}_{FB}( \bu_0^\eps, 0, p^\eps_0)_2) \in \mathcal{V}^\eps_{h, F} \times \mathcal{X}^\eps_h, 
        \\
        (\bfeta_h^0, \bfxi_h^0) &= (\mathcal{P}_B(\bfeta_0^\eps), \mathcal{P}_B(\bfxi_0^\eps))\in\mathcal{V}^\eps_{h, B}\times\mathcal{V}^\eps_{h, B},
    \end{aligned}
\end{equation}
where $\mathcal{P}_{FB}$ and $\mathcal{P}_B$ were defined in Lemma~\ref{lem:weighted-ritz}.
In order to define the data term evaluated at timestep $n$, 
for simplicity,
we assume of the right hand sides that
\begin{equation*}
    (\bF_F^\eps, \bF_B^\eps, g^\eps) \in
    C^0([0,T]; L^2(\Omega_F^\eps, \Phi_F^\eps)^d) \times 
    C^0([0,T]; L^2(\Omega_B^\eps,\Phi_B^\eps)^d) \times 
    C^0([0,T]; L^2(\Omega_B^\eps,\Phi_B^\eps)).
\end{equation*}
It is then legitimate to define, for $n = 1, \ldots, N$, the linear form $\mathcal{F}^n_\eps \coloneqq \mathcal{F}_\eps(t^n)$; see~\eqref{eq:F_eps}.

The scheme is then given  as follows. For $n = 0, \ldots, N - 1$, we seek $(\bu_h^{n + 1}, \theta_h^{n + 1}, \bfeta_h^{n + 1}, p_h^{n + 1} )\in \mathcal{V}^\eps_{h, F} \times \mathcal{Q}_h^\eps \times \mathcal{V}^\eps_{h, B} \times \mathcal{X}^\eps_h$ such that, for every $(\bv_{h}, \zeta_h, \bfphi_h, q_h) \in \mathcal{V}_{h, F}^{\eps} \times \mathcal{Q}_h^\eps \times \mathcal{V}^\eps_{h, B} \times \mathcal{X}^{\eps}_{h} $, we have:
\begin{equation}\label{eq:approxWF}
    \begin{aligned}
        &\rho_F \int_{\Omega_F^\eps} \dtee \bu_{h}^{n + 1} \cdot\bv_{h} \Phi_F^\eps + \rho_B\int_{\Omega_B^\eps} \dtt \bfeta^{n + 1}_h \cdot \bfphi_h \Phi_B^\eps + c_0\int_{\Omega_B^\eps} \dtee p_{h}^{n + 1} q_{h} \Phi_B^\eps
        - \frac1{2\eps}\int_{\ell^\eps} q_h \dtee \bfeta_h^{n + 1} \cdot \nabla \dist_\Gamma
        \\ &
        + \frac1{2\eps} \int_{\ell^\eps} p_h^{n + 1} \bfphi_h \cdot \nabla \dist_\Gamma
        + \frac{ \alpha_{BJ} }{2\eps} \sum_{i = 1}^{d - 1} \int_{\ell^\eps}((
        \bu_h^{n + 1}
        - \dtee\bfeta_h^{n + 1}
        )\cdot \tilde{\bftau}_i) ((\bv_h - \bfphi_h)\cdot \tilde{\bftau}_i)|\nabla\dist_\Gamma|
        \\ &
        - \alpha\int_{\Omega_B^\eps} \nabla \cdot \bfphi_h p_h^{n + 1} \Phi_B^\eps
        + \alpha\int_{\Omega_B^\eps} \nabla \cdot \dtee\bfeta_h^{n + 1} q_h \Phi_B^\eps
        + \cA_{\eps}( (\bu_h^{n + 1}, \bfeta_h^{n + 1}, p_h^{n + 1}), (\bv_h, \bfphi_h, q_h) ) 
        \\ &
        - \int_{\Omega^\eps_F}(\nabla\cdot\bv_h)\theta_h^{n + 1} + \int_{\Omega^\eps_F}(\nabla\cdot\bu_h^{n + 1})\zeta_h
        = \langle \mathcal{F}_\eps^{n + 1}, (\bv_h, \bfphi_h, q_h) \rangle.
    \end{aligned}
\end{equation}
Here for the centred finite difference approximation $\dtt\bfeta_h^{n + 1}$ to $\partial_{tt}\bfeta^\eps$,
it is understood that $\bfeta_h^{-1} = \bfeta_h^0 - \dt\bfxi_h^0$.
For 
further
details on implementation, the reader is referred to Section~\ref{sec:numerics}.

Essentially, the bilinear form $\cA_\eps$, defined in~\eqref{eq:A_eps}, collects together terms to which purely spatial derivatives are applied, whose approximation therefore follows by standard Galerkin projection; most of the remaining terms approximate time derivatives.
Note that terms of both type arise on the interface.

\begin{remark}[fluid pressure recovery]\label{rem:PiVsTheta}
    Notice that, in scheme~\eqref{eq:approxWF}, the terms involving incompressibility and the fluid pressure are missing the weight $\Phi_F^\eps$. This is 
    because 
    we have the compatibility condition~\eqref{eq:DiscrInfSupWeirdBCs}, but not 
    an analogous
    one involving the weight. Nevertheless, if one wishes to compute an approximate fluid pressure, this can 
    trivially 
    be obtained \emph{a posteriori} via
    \begin{equation*}
        \pi_h^\dt \coloneqq \frac1{\Phi_F^\eps} \theta_h^\dt.
    \end{equation*}
\end{remark}

\begin{theorem}[discrete stability and well-posedness]\label{thm:DiscrStability}
In the setting of Theorem~\ref{thm:WellPosedness} we have that, for every $\dt$ and $h$, scheme~\eqref{eq:approxWF} has a unique solution. Moreover, this solution satisfies the following energy identity 
\begin{multline*}
    \frac{\dt}{2}\left(
    \rho_F \dtee \|\bu_h^{n + 1}\|^2_{L^2(\Omega_F^\eps, \Phi_F^\eps)^d} + \rho_F \dt \| \dtee \bu_h^{n + 1} \|^2_{L^2(\Omega_F^\eps, \Phi_F^\eps)^d}
    + \rho_B \dtee \|\dtee\bfeta_h^{n + 1}\|^2_{L^2(\Omega_B^\eps, \Phi_B^\eps)^d} 
    \right. 
    \\ 
    \left.
    + \rho_B \dt \| \dtt \bfeta_h^{n + 1} \|^2_{L^2(\Omega_B^\eps, \Phi_B^\eps)^d} 
    + c_0 \dtee \|p_h^{n + 1}\|^2_{L^2(\Omega_B^\eps, \Phi_B^\eps )} + c_0 \dt \|\dtee p_h^{n + 1}\|^2_{L^2(\Omega_B^\eps, \Phi_B^\eps )}
    + \dtee\|\bfeta_h^{n + 1}\|^2_{E, \eps}
    \right. 
    \\ 
    \left.
    + \dt\|\dtee\bfeta_h^{n + 1}\|^2_{E, \eps}
    \right)
    + \dt\left(
    2\mu_F\|\bD(\bu_h^{n + 1})\|^2_{L^2(\Omega_F^\eps, \Phi_F^\eps)^{d \times d}}
    + \alpha_{BJ} \sum_{i = 1}^{d - 1} \|(\bu_h^{n + 1} - \dtee\bfeta_h^{n + 1})\cdot\tilde \bftau_i\|^2_{L^2\left(\ell^\eps, \frac{1}{2\eps}|\nabla\dist_\Gamma|\right)}
    \right.
    \\
    \left.
    + \|\boldsymbol \kappa_\eps^{\frac12} \nabla p_h^{n + 1}\|^2_{L^2(\Omega_B^\eps, \Phi_B^\eps )^d}
    \right)
    = 
    \dt \langle \mathcal{F}_\eps^{n + 1}, (\bu_h^{n + 1}, \dtee\bfeta^{n + 1}_h, p_h^{n + 1} )\rangle.
\end{multline*}
\end{theorem}
\begin{proof}
As in the proof of Theorem~\ref{thm:WellPosedness}, we first show the existence of $(\bu_h^{n + 1}, \bfeta_h^{n + 1}, p_h^{n + 1})$
by restricting~\eqref{eq:approxWF} at first to the subspace of discrete fluid velocities $\bv_h\in\mathcal{V}^\eps_{h, F}$ for which
\begin{equation*}
    \int_{\Omega^\eps_F}(\nabla\cdot\bv_h)\zeta_h = 0, \qquad \forall \zeta_h\in\mathcal{Q}^\eps_h.
\end{equation*}
This allows us to ignore $\theta^{n + 1}_h$ for now.
Isolating 
the bilinear terms, we 
can thus 
write~\eqref{eq:approxWF} as
\begin{equation}\label{eq:discrete_pre-divergence_unscaled}
    \mathcal{G}_\eps((\bu_h^{n + 1}, \bfeta_h^{n + 1}, p_h^{n + 1}), (\bv_h, \bfphi_h, q_h)) = \langle\mathfrak{F}_{n + 1, \eps}, (\bv_h, \bfphi_h, q_h)\rangle \quad~\forall~(\bv_h, \bfphi_h, q_h)\in\mathcal{V}_{h, F}^\eps\times\mathcal{V}_{h, B}^\eps\times\mathcal{X}_h^\eps,
\end{equation}
where
    \begin{align*}
        &\mathcal{G}_\eps((\bu_h, \bfeta_h, p_h), (\bv_h, \bfphi_h, q_h)) \coloneqq \frac{\rho_F}{\dt} \int_{\Omega_F^\eps} \bu_{h} \cdot\bv_{h} \Phi_F^\eps + \frac{\rho_B}{(\dt)^2}\int_{\Omega_B^\eps} \bfeta_h \cdot \bfphi_h \Phi_B^\eps + \frac{c_0}{\dt}\int_{\Omega_B^\eps} p_{h} q_{h} \Phi_B^\eps
        \\
        &+ \cA_\eps((\bu_h, \bfeta_h, p_h), (\bv_h, \bfphi_h, q_h))
        - \alpha\int_{\Omega^\eps_B}\left(\nabla\cdot\bfphi_h p_h - \nabla\cdot\frac{\bfeta_h}{\dt}q_h\right)\Phi^\eps_B
        - \frac1{2\eps}\int_{\ell^\eps}\left(q_h\frac{\bfeta_h}{\dt} - p_h\bfphi_h\right)\cdot\nabla\dist_\Gamma
        \\ &
        + \frac{\alpha_{BJ}}{2\eps} \sum_{i = 1}^{d - 1} \int_{\ell^\eps}\left(\left(\bu_h - \frac{\bfeta_h}{\dt}\right)\cdot \tilde{\bftau}_i\right) ((\bv_h - \bfphi_h)\cdot \tilde{\bftau}_i) |\nabla\dist_\Gamma|, 
        \\
        &\langle\mathfrak{F}_{n + 1, \eps}, (\bv_h, \bfphi_h, q_h) \rangle
        \coloneqq \langle \mathcal{F}_\eps^{n + 1}, (\bv_h, \bfphi_h, q_h) \rangle 
        + \rho_F \int_{\Omega_F^\eps} \frac{\bu_{h}^n}{\dt} \cdot\bv_{h} \Phi_F^\eps 
        + \rho_B\int_{\Omega_B^\eps} \frac{2\bfeta_h^n - \bfeta_h^{n - 1}}{(\dt)^2} \cdot \bfphi_h \Phi_B^\eps 
        \\ &
        + c_0\int_{\Omega_B^\eps} \frac{p_h^n}{\dt} q_{h} \Phi_B^\eps
        - \frac1{2\eps}\int_{\ell^\eps} q_h \frac{\bfeta_h^n}{\dt} \cdot \nabla \dist_\Gamma 
        - \frac{\alpha_{BJ}}{2\eps} \sum_{i = 1}^{d - 1} \int_{\ell^\eps}\left(\frac{\bfeta_h^n}{\dt}\cdot \tilde{\bftau}_i\right) ((\bv_h - \bfphi_h)\cdot \tilde{\bftau}_i) |\nabla\dist_\Gamma| 
        \\ &
        + \alpha\int_{\Omega^\eps_B}\nabla\cdot\frac{\bfeta^n_h}{\dt}q_h\Phi_B^\eps.
    \end{align*}
Mimicking 
a trick 
from~\cite{Bociu2021}, we show the equivalence of this discretised problem to a coercive one by rescaling a subset of the test functions to induce cancellations, or positive sign, in the rescaled form when testing its coercivity.
Consider 
\begin{equation}\label{eq:rescaled-discrete-form}
    \mathbb{G}_\eps((\bu_h, \bfeta_h, p_h), (\bv_h, \bfphi_h, q_h)) \coloneqq \mathcal{G}_\eps\left((\bu_h, \bfeta_h, p_h), \left(\bv_h, \frac{\bfphi_h}{\dt}, q_h\right)\right).
\end{equation}
The problem~\eqref{eq:discrete_pre-divergence_unscaled} is, by the universal quantifier, equivalent to
\begin{equation*}
    \mathbb{G}_\eps((\bu_h^{n + 1}, \bfeta_h^{n + 1}, p_h^{n + 1}), (\bv_h, \bfphi_h, q_h)) = 
    \left\langle\mathfrak{F}_{n + 1, \eps}, \left(\bv_h, \frac{\bfphi_h}{\dt}, q_h\right)\right\rangle
    \quad \forall~(\bv_h, \bfphi_h, q_h)\in\mathcal{V}_{h, F}^\eps\times\mathcal{V}_{h, B}^\eps\times\mathcal{X}_h^\eps.
\end{equation*}
We have
\begin{equation*}
    \begin{aligned}
        &\mathbb{G}_\eps((\bu_h, \bfeta_h, p_h), (\bu_h, \bfeta_h, p_h)) 
        = \frac{\rho_F}{\dt}\|\bu_h\|^2_{L^2(\Omega^\eps_F, \Phi_F^\eps)^d} 
        + \frac{\rho_B}{(\dt)^3}\|\bfeta_h\|^2_{L^2(\Omega^\eps_B, \Phi_B^\eps)^d} 
        + \frac{c_0}{\dt}\|p_h\|^2_{L^2(\Omega^\eps_B, \Phi_B^\eps)}
        \\ &
        + 2\mu_F\|\bD(\bu_h)\|^2_{L^2(\Omega_F^\eps, \Phi_F^\eps)^{d\times d}}
        + \frac{1}{\dt}\|\bfeta_h\|^2_{E, \eps}
        + \|\boldsymbol \kappa_\eps^{\frac12}\nabla p_h\|^2_{L^2(\Omega^\eps_B, \Phi_B^\eps)^d} 
        \\ &
        + \alpha_{BJ}\sum_{i = 1}^{d - 1}\left\|\left(\bu_h - \frac{\bfeta_h}{\dt}\right)\cdot\tilde{\bftau}_i\right\|^2_{L^2\left(\ell^\eps, \frac{1}{2\eps}|\nabla\dist_\Gamma|\right)}
        \geq \frac{c_0}{\dt}\|p_h\|_{L^2(\Omega^\eps_B, \Phi_F^\eps)}^2 
        + 2\mu_F\|\bD(\bu_h)\|^2_{L^2(\Omega_F^\eps, \Phi_F^\eps)^{d\times d}}
        \\ &
        + \frac{2\mu_B}{\dt}\|\bD(\bfeta_h)\|^2_{L^2(\Omega^\eps_B, \Phi_B^\eps)^{d\times d}}
        + \|\boldsymbol \kappa_\eps^{\frac12}\nabla p_h\|^2_{L^2(\Omega^\eps_B, \Phi_B^\eps)^d}
        \geq D\|(\bu_h, \bfeta_h, p_h)\|^2_{\mathcal{V}^\eps_F\times\mathcal{V}^\eps_B\times\mathcal{X}^\eps}, 
    \end{aligned}
\end{equation*}
where
\begin{equation}\label{eq:rescaled-coercivity-constant}
    D = \min\left\{\frac{2\mu_F}{C^{F}_K}, \frac{2\mu_B}{C^{B}_K\dt}, \frac{c_0}{\dt}, k_*\right\},
\end{equation}
and $C^F_K, C^B_K$ are the Korn constants from~\eqref{eq:weightedKorn} associated with $\mathcal{V}^\eps_F, \mathcal{V}^\eps_B$, respectively, which we may assume are independent of $\eps$.
In addition, the compatibility condition~\eqref{eq:DiscrInfSupWeirdBCs} yields the existence and uniqueness of $\theta_h^{n + 1}$.

Finally, to obtain the energy identity, it suffices to set $(\bv_h, \zeta_h, \bfphi_h, q_h) = \dt (\bu_h^{n + 1}, \theta_h^{n + 1}, \dtee\bfeta_h^{n + 1}, p_h^{n + 1})$ and apply the identity
$\int_D 2(\dtee W^{n + 1})W^{n + 1}\omega = \dtee\|W^{n + 1}\|^2_{L^2(D, \omega)} + \dt\|\dtee W^{n + 1}\|^2_{L^2(D, \omega)}$.
\begin{confidential}
\begin{equation*}
    \begin{aligned}
        &\rho_F\dt\int_{\Omega_F^\eps}\bF^\eps\cdot\bu_h^{n + 1}\Phi_F^\eps
        + \rho_B\dt\int_{\Omega^\eps_B}\bF_B^\eps\cdot\dtee\bfeta^{n + 1}_h\Phi_B^\eps
        + \rho_B\dt\int_{\Omega^\eps_B}g^\eps p_h^{n + 1}\Phi_B^\eps
        \\ &
        = \dt\rho_F\int_{\Omega^\eps_F}2\dtee\bu_h^{n + 1}\cdot\bu_h^{n + 1}\Phi_F^\eps 
        + \dt\rho_B\int_{\Omega^\eps_B}2\dtt\bfeta_h^{n + 1}\cdot \dtee\bfeta_h^{n + 1}\Phi_B^\eps 
        + c_0\dt \int_{\Omega^\eps_B}2(\dtee p^{n + 1}_h)p^{n + 1}_h\Phi_B^\eps
        \\ &
        + 2\dt\mu_F\int_{\Omega_F^\eps}\vert\bD(\bu_h^{n + 1})\vert^2\Phi_F^\eps 
        + \dt\int_{\Omega^\eps_B}2\mu_B\bD(\bfeta_h^{n + 1}):\bD(\dtee\bfeta_h^{n + 1})\Phi_B^\eps 
        + \lambda_B(\nabla\cdot\bfeta_h^{n + 1})(\nabla\cdot \dtee\bfeta_h^{n + 1})\Phi_B^\eps 
        \\ &
        + \dt\int_{\Omega^\eps_B}{\boldsymbol\kappa_\eps}\vert\nabla p^{n + 1}_h\vert^2\Phi_B^\eps
        \\ &
        + \underbrace{\frac{\dt}{2\eps}\int_{\ell^\eps}p_h^{n + 1}\bu_h^{n + 1}\cdot\nabla\dist_\Gamma - \frac{\dt}{2\eps}\int_{\ell^\eps}p_h^{n + 1}\dtee\bfeta_h^{n + 1}\cdot\nabla\dist_\Gamma - \frac{\dt}{2\eps}\int_{\ell^\eps}p_h^{n + 1}\bu_h^{n + 1}\cdot\nabla\dist_\Gamma + \frac{\dt}{2\eps}\int_{\ell^\eps}p_h^{n + 1}\dtee\bfeta_h^{n + 1}\cdot\nabla\dist_\Gamma}_{= 0}\\
        &+ \frac{\alpha_{BJ}}{2\eps} \dt \sum_{i = 1}^{d - 1} \int_{\ell^\eps} | (\bu_h^{n + 1} - \dtee\bfeta_h^{n + 1})\cdot\tilde \bftau_i |^2 |\nabla \dist_\Gamma|\\
        &\underbrace{- \alpha\dt\int_{\Omega^\eps_B}(\nabla\cdot \dtee\bfeta_h^{n + 1})p_h^{n + 1}\Phi_B^\eps 
        + \alpha\dt\int_{\Omega^\eps_B}(\nabla\cdot \dtee\bfeta_h^{n + 1})p_h^{n + 1}\Phi_B^\eps
    }_{= 0} 
        \underbrace{- \dt\int_{\Omega^\eps_F}(\nabla\cdot\bu^{n + 1}_h)\theta_h^{n + 1} + \dt\int_{\Omega^\eps_F}(\nabla\cdot\bu^{n + 1}_h)\theta_h^{n + 1}}_{= 0}\\
        &= \frac{\dt}{2}\left(\B{\rho_F\int_{\Omega^\eps_F}2\dtee\bu_h^{n + 1}\cdot\bu_h^{n + 1}\Phi_F^\eps + \rho_B\int_{\Omega^\eps_B}2\dtt\bfeta_h^{n + 1}\cdot \dtee\bfeta_h^{n + 1}\Phi_B^\eps + c_0 \int_{\Omega^\eps_B}2\dtee p^{n + 1}_hp^{n + 1}_h\Phi_B^\eps}\right)
        \\ &
        + 2\mu_F\dt\|\bD(\bu_h^{n + 1})\|^2_{L^2(\Omega^\eps_B, \Phi_F^\eps)^{d\times d}} 
        + \dt\B{\int_{\Omega^\eps_B}2\mu_B\bD(\bfeta_h^{n + 1}):\dtee\bD(\bfeta_h^{n + 1})\Phi_B^\eps 
        + \lambda_B(\nabla\cdot\bfeta_h^{n + 1})\dtee(\nabla\cdot\bfeta_h^{n + 1})\Phi_B^\eps}\\
        &+ \frac{\alpha_{BJ}}{2\eps} \dt \sum_{i = 1}^{d - 1} \int_{\ell^\eps} |(\bu_h^{n + 1} - \dtee\bfeta_h^{n + 1})\cdot\tilde \bftau_i|^2 |\nabla \dist_\Gamma|
        + \dt \|\boldsymbol \kappa_\eps^{\frac12} \nabla p_h^{n + 1}\|^2_{L^2(\Omega_B^\eps, \Phi_B^\eps)^d}.
    \end{aligned}
\end{equation*}
The energy identity is given by applying the following identity to the \B{blue terms}:
\begin{equation*}
    \begin{aligned}
        &\int_D 2(\dtee W^{n + 1})W^{n + 1}\omega = \frac{1}{\dt}\int_D\left(2\vert W^{n + 1}\vert^2 - 2W^{n + 1}W^n\right) \\
        &= \frac{1}{\dt}\left(\int_D\vert W^{n + 1}\vert^2\omega - \int_D\vert W^n\vert^2\omega + \int_D\left(\vert W^{n + 1}\vert^2 - 2W^{n + 1}W^n + \vert W^n\vert^2\right)\omega\right) \\
        &= \frac{1}{\dt}\left(\int_D\vert W^{n + 1}\vert^2\omega - \int_D\vert W^n\vert^2\omega\right) + \dt\int_D\left\vert\frac{W^{n + 1} - W^n}{\dt}\right\vert^2\omega 
        = \dtee\|W^{n + 1}\|^2_{L^2(D, \omega)} + \dt\|\dtee W^{n + 1}\|^2_{L^2(D, \omega)}.
    \end{aligned}
\end{equation*}
This may be seen as an application of the polarisation identity $(a - b)a = \frac12 a^2 - \frac12 b^2 + \frac12 (a - b)^2$.
\end{confidential}
\end{proof}

Observe that (in contrast to~\cite[Remark 5.2]{Bociu2021}) the coercivity constant~\eqref{eq:rescaled-coercivity-constant} is nonsingular (and in fact, improves) with small $\dt$, but this is arbitrary because (for example) the alternative test function rescaling $(\bv_h, \bfphi_h, q_h)\mapsto ((\dt)\bv_h, \bfphi_h, (\dt)q_h)$
gives the singular constant $(\dt)D$;
in any case,
this proves an inf-sup condition for 
the original form $\mathcal{G}_\eps$, 
which 
is not (or is not clearly) coercive,
and the dependence of later estimates on $\dt$ will be unaffected.

\subsection{Error analysis}\label{sec:numerical-error}

We now proceed to carry out an error analysis for scheme~\eqref{eq:approxWF}. To begin, let us introduce $\theta^\eps \coloneqq \pi^\eps \Phi_F^\eps \in H^{-1}(0, T; \mathcal{Q}^\eps)$; compare with Remark~\ref{rem:PiVsTheta}. Next, we must make suitable smoothness assumptions. Namely, we assume that there is $k \leq r$, where $r$ was defined in Lemma~\ref{lem:weighted-ritz}, for which
\begin{equation}\label{eq:SolIsSmooth}
    \bu^\eps \in C^2([0, T]; \mathcal{V}^{k + 1, \eps}_F),
    \quad \theta^\eps \in C^1([0, T]; \mathcal{Q}^{k, \eps}), 
    \quad \bfeta^\eps \in C^4([0, T]; \mathcal{V}^{k + 1, \eps}_B), 
    \quad p^\eps \in C^2([0, T]; \mathcal{X}^{k + 1, \eps}). 
\end{equation}

Under assumption \eqref{eq:SolIsSmooth} we first observe that, in \eqref{eq:SBFP_weak}, the terms involving the fluid pressure can be equivalently rewritten. Namely,
\begin{equation}\label{eq:NewWeakDiffuse}
    \begin{aligned}
        -\langle \pi^\eps, \nabla \cdot \bv \rangle_{(0,T) \times \Omega_F^\eps, \Phi_F^\eps} + \int_0^T \int_{\Omega_F^\eps} (\nabla \cdot \bu^{\eps}) \zeta \Phi_F^\eps &=
        - \int_0^T\int_{\Omega_F^\eps} (\nabla \cdot \bv) \pi^{\eps} \Phi_F^\eps
        + \int_0^T \int_{\Omega_F^\eps} (\nabla \cdot \bu^{\eps}) \zeta \Phi_F^\eps \\
         &=
        - \int_0^T\int_{\Omega_F^\eps} (\nabla \cdot \bv) \theta^\eps
        + \int_0^T\int_{\Omega_F^\eps} (\nabla \cdot \bu^{\eps}) \zeta,
    \end{aligned}
\end{equation}
where now $\zeta \in \mathcal{Q}^\eps$. Next, the regularity assumptions~\eqref{eq:SolIsSmooth} allow us to write the diffuse spacetime weak formulation~\eqref{eq:SBFP_weak} as a weak formulation over only the spatial domains, holding for a.e.~time. By density, we extend the spatial test functions to the Sobolev spaces in which the solution tuples take values: for all $(\bv, \zeta, \bfphi, q)\in \mathcal{V}_F^{\eps} \times \mathcal{Q}^\eps \times \mathcal{V}_B^{\eps} \times \mathcal{X}^{\eps}$,
\begin{multline}\label{eq:diffuse-without-time-integrals}
    \rho_F \int_{\Omega_F^\eps}\partial_t\bu^{\eps}\cdot\bv \Phi_F^\eps
    + \rho_B\int_{\Omega_B^\eps}\partial_{tt} \bfeta^{\eps}\cdot \bfphi \Phi_B^\eps
    + \int_{\Omega_B^\eps} (c_0\partial_t p^{\eps} + \alpha\nabla\cdot\partial_t\bfeta^{\eps})q \Phi_B^\eps
    \\
    + \cA_\eps((\bu^\eps(t), \bfeta^\eps(t), p^\eps(t)), (\bv, \bfphi, q))
    - \alpha\int_{\Omega_B^\eps} \nabla \cdot \bfphi p^{\eps} \Phi_B^\eps
    - \frac1{2\eps}\int_{\ell^\eps} q \partial_t \bfeta^{\eps} \cdot \nabla \dist_\Gamma
    + \frac1{2\eps} \int_{\ell^\eps} p^{\eps} \bfphi \cdot \nabla \dist_\Gamma
    \\
    + \frac{ \alpha_{BJ} }{2\eps} \sum_{i = 1}^{d - 1} \int_{\ell^\eps}((\bu^{\eps} -\partial_t \bfeta^{\eps}) \cdot \tilde{\bftau}_i) ((\bv - \bfphi) \cdot \tilde{\bftau}_i) |\nabla\dist_\Gamma|
    - \int_{\Omega_F^\eps} (\nabla \cdot \bv) \theta^{\eps}
    + \int_{\Omega_F^\eps} (\nabla \cdot \bu^{\eps}) \zeta 
    \\
    = \langle\mathcal{F}_\eps, (\bv, \bfphi, q)\rangle.
\end{multline}
It is with these modifications that our problem can be compared to the scheme~\eqref{eq:approxWF}. 

We can now begin the error analysis \emph{per se}. We construct sequences $\bu^\dt = \{ \bu^{n} = \bu^\eps(t^n) \}_{n = 0}^N$, $\theta^\dt = \{ \theta^{n} = \theta^\eps(t^n) \}_{n = 0}^N$, $\bfeta^\dt = \{\bfeta^n = \bfeta^\eps(t^n) \}_{n = 0}^N$, and $p^\dt = \{ p^{n} = p^\eps(t^n) \}_{n = 0}^N$. We define the errors $\be_\bu^\dt = (\bu^\eps)^\dt - \bu_h^\dt$,
$e_\theta^\dt = (\theta^\eps)^\dt - \theta_h^\dt$,
$\be_{\bfeta}^\dt = (\bfeta^\eps)^\dt - \bfeta_h^\dt$, and 
$e_p^\dt  = (p^\eps)^\dt - p_h^\dt$.
Next, as is standard, we decompose the error into the \textit{projection} and \textit{truncation} errors
as
\begin{equation}
    \begin{aligned}\label{eq:truncation-errors}
        \be_{\bu}^\dt &= ((\bu^\eps)^\dt - \bU_h^\dt ) + (\bU_h^\dt - \bu_h^\dt) =: \bY_{\bu}^\dt + \bE_{\bu, h}^\dt, \\
        e_\theta^\dt &= (\theta^\dt - \Theta_h^\dt) + (\Theta_h^\dt - \theta_h^\dt) =: Y_{\theta}^\dt + E_{\theta, h}^\dt, \\
        \be_{\bfeta}^\dt &= ((\bfeta^\eps)^\dt - \bfETA_h^\dt) + (\bfETA_h^\dt - \bfeta_h^\dt) =: \bY_{\bfeta}^\dt + \bE_{\bfeta, h}^\dt, \\
        e_p^\dt &= ((p^\eps)^\dt - P_h^\dt) + (P_h^\dt - p_h^\dt) =: Y_{p}^\dt + E_{p, h}^\dt.
    \end{aligned}
\end{equation}
Here, $(\bU_h^\dt, \Theta_h^\dt, P_h^\dt) = \mathcal{P}_{FB}((\bu^\eps)^\dt, (\theta^\eps)^\dt, (p^\eps)^\dt)$
and $\bfETA_h^\dt = \mathcal{P}_{B}((\bfeta^\eps)^\dt)$. 

We can now proceed with our main error estimate.

\begin{theorem}[error estimate]
Let $(\bu^\eps, \theta^\eps, \bfeta^\eps, p^\eps)$ solve~\eqref{eq:SBFP_weak} and satisfy~\eqref{eq:SolIsSmooth}. Let $(\bu_h^\dt, \theta_h^\dt, \bfeta_h^\dt, p_h^\dt)$ solve~\eqref{eq:approxWF} with initial data~\eqref{eq:discrete-initial-data}.
The following error estimate holds
\begin{multline}\label{eq:overall-error-estimate}
    \| (\bu^\eps)^\dt - \bu_h^\dt \|_{L^\infty_\dt(0, T; L^2(\Omega_F^\eps, \Phi_F^\eps)^d)} 
    + \| (\bfeta^\eps)^\dt - \bfeta_h^\dt \|_{L^\infty_\dt(0, T; \mathcal{V}^\eps_B )} 
    \\
    + \|(\partial_t\bfeta^\eps)^\dt - \dtee\bfeta^\dt_h\|_{L^\infty_\dt(0, T; L^2(\Omega^\eps_B, \Phi^\eps_B)^d)} 
    + \| (p^\eps)^\dt - p_h^\dt \|_{L^\infty_\dt(0, T; L^2(\Omega_B^\eps, \Phi_B^\eps ))} 
    \\
    + \| (\bu^\eps)^\dt - \bu_h^\dt \|_{L^2_\dt(0, T; \mathcal{V}_F^\eps)} 
    + \| \nabla((p^\eps)^\dt - p_h^\dt) \|_{L^2_\dt(0, T; L^2(\Omega^\eps_B, \Phi^\eps_B)^d)} 
    \lesssim \dt + h^{k},
\end{multline}
where the implied constant depends on the material parameters, the final time $T$, and the domain $\Omega$. They are independent of $\dt$ and $h$, and depend on $\eps$ only through higher order norms of the exact solution.
\end{theorem}
\begin{proof}
Upon restricting, in the new formulation~\eqref{eq:diffuse-without-time-integrals}, 
the test functions to lie on the corresponding finite element spaces, we obtain
\begin{multline*}
    \rho_F \int_{\Omega_F^\eps} \dtee \be_{\bu}^{n + 1} \cdot\bv_{h} \Phi_F^\eps 
    + \rho_B\int_{\Omega_B^\eps}\dtt\be_{\bfeta}^{n + 1}\cdot\bfphi_h \Phi_B^\eps 
    + c_0\int_{\Omega_B^\eps} \dtee e_{p}^{n + 1} q_{h} \Phi_B^\eps
    + \alpha\int_{\Omega_B^\eps} \nabla \cdot \dtee\be_{\bfeta}^{n + 1} q_h \Phi_B^\eps
    \\
    - \alpha\int_{\Omega^\eps_B}(\nabla\cdot\bfphi_h)e^{n + 1}_p\Phi_B^\eps
    - \frac1{2\eps}\int_{\ell^\eps} q_h \dtee \be_{\bfeta}^{n + 1} \cdot \nabla \dist_\Gamma
    + \frac{1}{2\eps}\int_{\ell^\eps}e^{n + 1}_p\bfphi_h\cdot\nabla\dist_\Gamma
    \\
    + \frac{ \alpha_{BJ} }{2\eps} \sum_{i = 1}^{d - 1} \int_{\ell^\eps}((\be^{n + 1}_\bu - \dtee\be_{\bfeta}^{n + 1})\cdot \tilde{\bftau}_i) ((\bv_h - \bfphi_h)\cdot \tilde{\bftau}_i)|\nabla\dist_\Gamma|
    + \cA_{\eps}( (\be_{\bu}^{n + 1}, \be_{\bfeta}^{n + 1}, e_p^{n + 1}), (\bv_h, \bfphi_h, q_h) )
    \\
        - \int_{\Omega^\eps_F}(\nabla\cdot\bv_h)e_\theta^{n + 1} + \int_{\Omega^\eps_F}(\nabla\cdot\be_{\bu}^{n + 1})\zeta_h
    = \langle \mathcal{R}_\eps^{n + 1}, (\bv_h, \bfphi_h, q_h) \rangle,
\end{multline*}
where the time consistency error $\cR^\dt_\eps$ is defined and estimated in Lemma~\ref{lem:ConsistencyDerivatives}.
Note 
for the truncation errors~\eqref{eq:truncation-errors}
that $\bE^0_{\bu, h} = \bE^0_{\bfeta, h} = \boldsymbol0$, $E^0_{p, h} = E^0_{\theta, h} = 0$. 
Moreover we have, for all $n = 0, \ldots, N - 1$,
\begin{equation*}
    \cA_\eps((\bY_{\bu}^{n + 1}, \bY_{\bfeta}^{n + 1}, Y_{p}^{n + 1}), (\bv_h, \bfphi_h, q_h))
    - \int_{\Omega^\eps_F}(\nabla\cdot\bv_h)Y_{\theta}^{n + 1} + \int_{\Omega^\eps_F}(\nabla\cdot\bY_{\bu}^{n + 1})\zeta_h = 0,
\end{equation*}
where $(\bv_h, \zeta_h, \bfphi_h, q_h) \in \mathcal{V}_{h, F}^\eps \times \mathcal{Q}^\eps_h \times \mathcal{V}^\eps_{h, B} \times \mathcal{X}_h^\eps$ are arbitrary. 
The decomposition~\eqref{eq:truncation-errors} thus yields
\begin{multline}\label{eq:decomposed-error-eqns}
    \rho_F \int_{\Omega_F^\eps} \dtee \bE_{\bu, h}^{n + 1} \cdot\bv_{h} \Phi_F^\eps 
    + \rho_B\int_{\Omega_B^\eps}\dtt\bE_{\bfeta, h}^{n + 1}\cdot\bfphi_h \Phi_B^\eps 
    + c_0\int_{\Omega_B^\eps} \dtee E_{p, h}^{n + 1} q_{h} \Phi_B^\eps 
    + \alpha\int_{\Omega^\eps_B}\nabla\cdot\bE_{\bfeta, h}^{n + 1}q_h\Phi^\eps_B
    \\
    - \alpha\int_{\Omega^\eps_B}(\nabla\cdot\bfphi_h)E_{p, h}^{n + 1}\Phi^\eps_B
    - \frac1{2\eps}\int_{\ell^\eps} q_h \dtee \bE_{\bfeta, h}^{n + 1} \cdot \nabla \dist_\Gamma
    + \frac{1}{2\eps}\int_{\ell^\eps}E^{n + 1}_{p, h}\bfphi_h\cdot\nabla\dist_\Gamma
    \\
    + \frac{\alpha_{BJ}}{2\eps} \sum_{i = 1}^{d - 1} \int_{\ell^\eps}((\bE^{n + 1}_{\bu, h} - \dtee\bE_{\bfeta, h}^{n + 1})\cdot \tilde{\bftau}_i) ((\bv_h - \bfphi_h)\cdot \tilde{\bftau}_i)|\nabla\dist_\Gamma|
    + \cA_\eps((\bE_{\bu, h}^{n + 1}, \bE^{n + 1}_{\bfeta, h}, E_{p, h}^{n + 1}), (\bv_h, \bfphi_h, q_h))
    \\
    - \int_{\Omega^\eps_F}(\nabla\cdot\bv_h)E_{\theta, h}^{n + 1} 
    + \int_{\Omega^\eps_F}(\nabla\cdot\bE_{\bu, h}^{n + 1})\zeta_h
    = \langle \mathcal{R}_\eps^{n + 1}, (\bv_h, \bfphi_h, q_h) \rangle
    \\
    - \rho_F \int_{\Omega_F^\eps} \dtee \bY_{\bu}^{n + 1} \cdot\bv_{h} \Phi_F^\eps 
    - \rho_B\int_{\Omega_B^\eps}\dtt\bY_{\bfeta}^{n + 1}\cdot\bfphi_h \Phi_B^\eps
    - c_0\int_{\Omega_B^\eps} \dtee Y_{p}^{n + 1} q_{h} \Phi_B^\eps
    - \alpha\int_{\Omega^\eps_B}\nabla\cdot\dtee\bY_{\bfeta}^{n + 1}q_h\Phi^\eps_B
    \\
    + \alpha\int_{\Omega^\eps_B}(\nabla\cdot\bfphi_h)Y_{p}^{n + 1}\Phi^\eps_B
    + \frac1{2\eps}\int_{\ell^\eps} q_h \dtee\bY_{\bfeta}^{n + 1} \cdot \nabla \dist_\Gamma
    - \frac{1}{2\eps}\int_{\ell^\eps}Y_{p}^{n + 1}\bfphi_h\cdot\nabla\dist_\Gamma
    \\
    - \frac{\alpha_{BJ}}{2\eps} \sum_{i = 1}^{d - 1} \int_{\ell^\eps}((\bY_{\bu}^{n + 1} -\dtee\bY_{\bfeta}^{n + 1})\cdot \tilde{\bftau}_i) ((\bv_h - \bfphi_h)\cdot \tilde{\bftau}_i)|\nabla\dist_\Gamma|.
\end{multline}
Now choose $(\bv_h, \zeta_h, \bfphi_h, q_h) = (\bE^{n + 1}_{\bu, h}, E^{n + 1}_{\theta, h}, \dtee\bE^{n + 1}_{\bfeta, h}, E^{n + 1}_{p, h})$.
The divergence terms become
\begin{equation*}
    (\pm)\int_{\Omega^\eps_F}(\nabla\cdot\bE^{n + 1}_{\bu, h})E^{n + 1}_{\theta_h} 
    = \int_{\Omega^\eps_F}(\nabla\cdot\bU^{n + 1}_h)E^{n + 1}_{\theta, h} 
    \int_{\Omega^\eps_F}(\nabla\cdot\bu^{n + 1}_h)E^{n + 1}_\theta
    =
    - \int_{\Omega^\eps_F}(\nabla\cdot (\bu^\eps)^{n + 1})E^{n + 1}_{\theta, h} = 0.
\end{equation*}
Using (the proof of) the discrete energy balance, Theorem~\ref{thm:DiscrStability}, we multiply the decomposed error equations~\eqref{eq:decomposed-error-eqns} by $\dt$. 
Now choosing $1 \leq K \leq N$ and taking the sum over $n = 0, \ldots, K - 1$ (denoted $\sum_0^{K - 1}$),
we obtain
\begin{multline*}
    \frac{1}2 \left(
        \rho_F\|\bE_{\bu, h}^{K}\|^2_{L^2(\Omega_F^\eps, \Phi_F^\eps)^d} +
        \rho_B \|\dtee\bE_{\bfeta, h}^{K}\|^2_{L^2(\Omega_B^\eps, \Phi_B^\eps)^d} +
        2 \mu_B \|\bD(\bE_{\bfeta, h}^{K})\|^2_{L^2(\Omega_B^\eps, \Phi_B^\eps)^{d \times d}} 
        + \lambda_B \|\nabla\cdot\bE_{\bfeta, h}^{K}\|^2_{L^2(\Omega^\eps_B, \Phi_B^\eps)} 
    \right. \\ \left.
    + c_0 \|E_{p, h}^{K}\|^2_{L^2(\Omega_B^\eps, \Phi_B^\eps )} 
    \right) 
    + \frac{(\dt)^2}2 \sum_0^{K - 1} \left(
        \rho_F \| \dtee \bE_{\bu, h}^{n + 1} \|^2_{L^2(\Omega_F^\eps, \Phi_F^\eps)^d}
        + \rho_B \| \dtt \bE_{\bfeta, h}^{n + 1} \|^2_{L^2(\Omega_B^\eps, \Phi_B^\eps)^d}
    \right. \\ \left.
        + 2\mu_B \|\bD(\dtee\bE_{\bfeta, h}^{n + 1})\|^2_{L^2(\Omega^\eps_B, \Phi_B^\eps)^{d\times d}} 
    + \lambda_B \|\nabla\cdot \dtee\bE_{\bfeta, h}^{n + 1}\|^2_{L^2(\Omega^\eps_B, \Phi_B^\eps)}+
        c_0 \|\dtee E_{p, h}^{n + 1}\|^2_{L^2(\Omega_B^\eps, \Phi_B^\eps )}
    \right)
    \\
    + \dt \sum_0^{K-1} \left(
        2\mu_F \|\bD(\bE_{\bu, h}^{n + 1})\|^2_{L^2(\Omega_F^\eps, \Phi_F^\eps)^{d \times d}} +
        \|\boldsymbol\kappa_\eps^{\frac12} \nabla E_{p, h}^{n + 1}\|^2_{L^2(\Omega_B^\eps, \Phi_B^\eps )^d} 
    \right.\\
    \left.
        + \alpha_{BJ} \sum_{i = 1}^{d - 1} \|(\bE_{\bu, h}^{n + 1} - \dtee\bE_{\bfeta, h}^{n + 1})\cdot\tilde \bftau_i\|^2_{L^2\left(\ell^\eps, \frac{1}{2\eps}|\nabla\dist_\Gamma|\right)}
    \right) 
    = \dt\sum_0^{K - 1}\langle \mathcal{R}_\eps^{n + 1}, (\bE^{n + 1}_{\bu, h}, \dtee\bE^{n + 1}_{\bfeta, h}, E^{n + 1}_{p, h}) \rangle
    \\
    - \rho_F \dt\int_{\Omega_F^\eps} \sum_0^{K - 1} \dtee \bY_{\bu}^{n + 1} \cdot\bE_{\bu, h}^{n + 1} \Phi_F^\eps
    - \rho_B\dt \int_{\Omega_B^\eps}\sum_0^{K - 1} \dtt\bY_{\bfeta}^{n + 1}\cdot\dtee\bE_{\bfeta, h}^{n + 1} \Phi_B^\eps
    - c_0\dt \int_{\Omega_B^\eps}\sum_0^{K - 1} \dtee Y_{p}^{n + 1} E_{p, h}^{n + 1} \Phi_B^\eps
    \\
    - \alpha\dt\int_{\Omega^\eps_B}\sum_0^{K - 1}\nabla\cdot\dtee\bY_{\bfeta}^{n + 1}E_{p, h}^{n + 1}\Phi^\eps_B
    + \alpha\dt\int_{\Omega^\eps_B}\sum_0^{K - 1}(\nabla\cdot\dtee\bE_{\bfeta, h}^{n + 1})Y_{p}^{n + 1}\Phi^\eps_B
    + \frac{\dt}{2\eps}\int_{\ell^\eps} \sum_0^{K - 1} E_{p, h}^{n + 1} \dtee\bY_{\bfeta}^{n + 1} \cdot \nabla \dist_\Gamma
    \\
    - \frac{\dt}{2\eps}\int_{\ell^\eps}\sum_0^{K - 1} Y_{p}^{n + 1}\dtee\bE_{\bfeta, h}^{n + 1}\cdot\nabla\dist_\Gamma
    - \frac{\alpha_{BJ}}{2\eps}\sum_{i = 1}^{d - 1} \dt\sum_0^{K - 1}\int_{\ell^\eps}((\bY_{\bu}^{n + 1} -\dtee\bY_{\bfeta}^{n + 1})\cdot \tilde{\bftau}_i) ((\bE_{\bu, h}^{n + 1} - \dtee\bE_{\bfeta, h}^{n + 1})\cdot \tilde{\bftau}_i)|\nabla\dist_\Gamma|,
\end{multline*}
where we applied a telescoping sum and used that the truncation errors vanish at $n = 0$.
Applying Lemma~\ref{lem:ConsistencyDerivatives} 
and denoting the term with constant $C$ therein by $C^*(\dt)^2$,
applying Young's inequality, 
discrete integration by parts (see e.g.~\cite[eq.~A2]{Bukac2015}),
and letting $\|\nabla\cdot\bv\|_{L^2(D, \omega)}^2\leq C_d\|\nabla\bv\|^2_{L^2(D, \omega)^{d\times d}}$,
we have
\begin{multline*}
    \frac{\rho_F}{2}\|\bE^K_{\bu, h}\|^2_{L^2(\Omega^\eps; \Phi^\eps_F)^d}
    + \frac{\rho_B}{2}\|\dtee\bE^K_{\bfeta, h}\|^2_{L^2(\Omega^\eps_B, \Phi^\eps_B)^d} 
    + \frac{3\mu_B}{5C_K}\|\bE^K_{\bfeta, h}\|^2_{\mathcal{V}^\eps_B}
    + \frac{c_0}{2}\|E^K_{p, h}\|^2_{L^2(\Omega^\eps_B, \Phi^\eps_B)}
    \\
    + \frac{\mu_F}{C_K}\dt\sum_0^{K - 1}\|\bE^{n + 1}_{\bu, h}\|^2_{\mathcal{V}^\eps_F}
    + \frac{k_*}{2}\dt\sum_0^{K - 1}\|\nabla E^{n + 1}_{p, h}\|^2_{L^2(\Omega^\eps_B, \Phi^\eps_B)^d}
    \\
    \leq C^*(\dt)^2
    + \left(
        \frac{\alpha^2}{\lambda_B}
        + \frac{5C_KC_{\rm tr}^2}{2\mu_B}
    \right)\| Y^K_{p}\|^2_{\mathcal{X}^\eps}
    + \frac{5C_K\alpha_{BJ}^2C_{\rm tr}^2(d - 1)^2}{2\mu_B}\|\bY^K_{\bu}\|^2_{\mathcal{V}^\eps_F}
    \\
    + \frac{5C_K\alpha_{BJ}^2C_{\rm tr}^2(d - 1)^2}{2\mu_B}\|\dtee\bY^K_{\bfeta}\|^2_{\mathcal{V}^\eps_B}
    + \left(\frac{\rho_F^2C_K}{\mu_F} + \frac{\alpha_{BJ}C_{\rm tr}(d - 1)}{2}\right)\dt\sum_0^{K - 1}\|\dtee\bY^{n + 1}_{\bu}\|^2_{\mathcal{V}^\eps_F}
    \\
    + \left(\frac{\rho_B}{2} + \frac{\alpha_{BJ}C_{\rm tr}(d - 1)}{2}\right)\dt\sum_0^{K - 1}\|\dtt\bY^{n + 1}_{\bfeta}\|^2_{\mathcal{V}^\eps_B}
    + \frac{1}{2}(c_0 + \alpha + C_{\rm tr})\dt\sum_0^{K - 1}\|\dtee Y^{n + 1}_{p}\|^2_{\mathcal{X}^\eps}
    \\
    + \left(\frac{\alpha C_d}{2} + \frac{2C_{\rm tr}^2}{k_*} + \frac{\mu_F}{4C_K}\right)\dt\sum_0^{K - 1}\|\dtee\bY^{n + 1}_{\bfeta}\|^2_{\mathcal{V}^\eps_B}
    + \frac{C_K\alpha_{BJ}^2(d - 1)^2C_{\rm tr}^2}{\mu_F}\dt\sum_0^{K - 1}\|\bY^{n + 1}_{\bu}\|^2_{\mathcal{V}^\eps_F} 
    \\
    + \left(
        \frac{\mu_B}{10C_K}
        + 
        \frac{C_{\rm tr}}{2}\left(1 + 2\alpha_{BJ}(d - 1)\right) + \frac{\alpha C_d}{2}\right)\dt\sum_0^{K - 1}\|\bE^{n + 1}_{\bfeta, h}\|^2_{\mathcal{V}^\eps_B}
    \\
    + \frac{\rho_B}{2}\dt\sum_0^{K - 1}\|\dtee\bE^{n + 1}_{\bfeta, h}\|^2_{L^2(\Omega^\eps_B, \Phi^\eps_B)^d}
    + \frac{1}{2}(c_0 + \alpha + k_*)\dt\sum_0^{K - 1}\|E^{n + 1}_{p, h}\|^2_{L^2(\Omega^\eps_B, \Phi^\eps_B)}
    + \dt\sum_0^{K - 1}\|\bE^{n + 1}_{\bu, h}\|^2_{L^2(\Omega_F^\eps,\Phi_F^\eps)^d}.
\end{multline*}
By the discrete Gr\"onwall lemma,
\begin{multline*}
        \|\bE^K_{\bu, h}\|^2_{L^2(\Omega^\eps; \Phi^\eps_F)^d}
    + 
        \|\dtee\bE^K_{\bfeta, h}\|^2_{L^2(\Omega^\eps_B, \Phi^\eps_B)^d}
    + 
        \|\bE^K_{\bfeta, h}\|^2_{\mathcal{V}^\eps_B}
    + \|E^K_{p, h}\|^2_{L^2(\Omega^\eps_B, \Phi^\eps_B)}
    + \dt\sum_0^{K - 1}\|\bE^{n + 1}_{\bu, h}\|^2_{\mathcal{V}^\eps_F}
    \\
    + 
        \dt\sum_0^{K - 1}\|\nabla E^{n + 1}_{p, h}\|^2_{L^2(\Omega^\eps_B, \Phi^\eps_B)^d}
    \\
    \lesssim 
    {\rm e}^T \left(
        (\dt)^2
    + 
        \| Y^K_{p}\|^2_{\mathcal{X}^\eps}
    + 
        \|\bY^K_{\bu}\|^2_{\mathcal{V}^\eps_F}
    + 
        \|\dtee\bY^K_{\bfeta}\|^2_{\mathcal{V}^\eps_B}
    + 
        \dt\sum_0^{K - 1}\|\dtee\bY^{n + 1}_{\bu}\|^2_{\mathcal{V}^\eps_F}
    + 
        \dt\sum_0^{K - 1}\|\dtt\bY^{n + 1}_{\bfeta}\|^2_{\mathcal{V}^\eps_B}
    \right. \\ \left.
    + 
        \dt\sum_0^{K - 1}\|\dtee Y^{n + 1}_{p}\|^2_{\mathcal{X}^\eps}
    + 
        \dt\sum_0^{K - 1}\|\dtee\bY^{n + 1}_{\bfeta}\|^2_{\mathcal{V}^\eps_B}
    + 
        \dt\sum_0^{K - 1}\|\bY^{n + 1}_{\bu}\|^2_{\mathcal{V}^\eps_F}
    \right).
\end{multline*}

Finally, on the left-hand-side, we either remove terms summed over $n$ before taking $\max_{K = 0}^{N}$, 
or remove terms evaluated at time $K$ before taking $K = N$. The result corresponds then to each term on the left hand side of~\eqref{eq:overall-error-estimate}, which is finally obtained by the triangle inequality and Lemmas~\ref{lem:weighted-ritz} and~\ref{lem:more-projection-errors}.
\end{proof}

\begin{remark}[(sub)optimality]
    The first four terms on the left hand side of~\eqref{eq:overall-error-estimate} appear to converge at a suboptimal rate in space. This is because, in the error analysis, we need to handle the terms involving integrals on the diffuse interface $\ell^\eps$ via the weighted trace inequality of Lemma~\ref{lem:nablaPhi}. Notice that not even the projection errors are superconvergent, as the definition of $\mathcal{P}_{FB}$, see Lemma~\ref{lem:weighted-ritz}, includes diffuse interface integrals. Nevertheless, this is consistent with existing error estimates for the sharp case (e.g.~\cite{Bukac2015}), and the modelling errors we present below.
\end{remark}

\section{Modelling error}\label{sec:ModellingError}

\subsection{Positive Lipschitz weights}
In this section, we consider regularised phase field functions, which  yield positive Lipschitz weights. In particular, we define
\begin{equation}\label{eq:regularizacija}
    \Phi_F^{\eps, \delta} \coloneqq (1 - 2\delta) \Phi_F^{\eps} + \delta = \Phi_F^{\eps} + \delta(1 - 2\Phi_F^{\eps}),
    \qquad
    \Phi_B^{\eps, \delta} \coloneqq 1 - \Phi_F^{\eps,\delta},
\end{equation}
where $\delta > 0$ is a small regularisation parameter, and $\Phi_F^\eps$ and $\Phi_B^\eps$ are defined in the same way as in~\eqref{eq:defofPhiPsi}, but using a different function from the one used in~\eqref{eq:defofPhiPsi}. Namely, 
we replace  $\cS$ by $\cS_r$, where $\cS_r$ satisfies conditions (S1)--(S3) from~\cite[Section 3]{Burger2017}. Some examples of such functions are
\begin{equation*}
    \cS_r(t) \coloneqq \begin{dcases}
               t, & |t| \leq 1, 
               \\
               \frac{t}{|t|}, & |t| > 1,
        \end{dcases}
    \qquad\qquad \textrm{or} \qquad\qquad
    \cS_r(t) \coloneqq \tanh(t).
\end{equation*}
Notice that since both $\Phi_F^{\eps, \delta}$ and $\Phi_B^{\eps, \delta}$ are strictly positive they satisfy all conditions from Sections~\ref{sub:DiffuseDomainWellPosed} and~\ref{sec:Discretisation} with $\Omega_F^{\eps} = \Omega_B^{\eps} = \Omega$. In addition, $\nabla\Phi_F^{\eps}$ and $\nabla\Phi_F^{\eps, \delta}$ are colinear and therefore we can define tangents $\{\tilde{\bftau}_i\}_{i = 1}^{d - 1}$ using the regularised phase field functions in the same way as before.
In analogy to~\eqref{eq:diffuse-energy-norm}, we define
\begin{equation*}
    \|\bfphi\|^2_{E, \eps, \delta} \coloneqq 2\mu_B\|\bD(\bfphi)\|^2_{L^2(\Omega, \Phi^{\eps, \delta}_B)^{d\times d}} + \lambda_B\|\nabla\cdot\bfphi\|^2_{L^2(\Omega, \Phi^{\eps, \delta}_B)}.
\end{equation*}

In this section, our goal is to derive rates of convergence of the continuous diffuse interface method to the continuous sharp interface method in terms of the parameters $\epsilon$ and $\delta$. This will be done by subtracting the weak form of the diffuse interface formulation from the weak form of the sharp interface formulation, taking suitable test functions, and using the convergence properties of diffuse integrals. For the sake of clarity, in this section we introduce subscripts ``\emph{s}'' and ``\emph{d}'' to denote the solutions to the sharp interface and the diffuse interface formulations, respectively. The modelling errors are then denoted by
\begin{equation*}
    \be_\bu^{\eps, \delta} \coloneqq \bu_s - \bu_d^{\eps, \delta}, \qquad \be_\bfeta^{\eps, \delta} \coloneqq \bfeta_s - \bfeta_d^{\eps, \delta}, \qquad e_p^{\eps, \delta} \coloneqq p_s - p_d^{\eps, \delta}.
\end{equation*}
Our main result of this section is given as follows. 

\begin{theorem}[modelling error I]
Let $(\bu_s, \pi_s, \bfeta_s, p_s)$ be the solution of the sharp interface Stokes--Biot problem in the sense of Definition~\ref{def:WS}, and $(\bu^{\epsilon, \delta}_d, \pi^{\eps, \delta}_d, \bfeta^{\epsilon, \delta}_d, p^{\epsilon, \delta}_d)$ be the solution of the diffuse interface problem in sense of Definition~\ref{def:SBFP}. Assume that the sharp interface solution and the forcing terms are sufficiently smooth, along with their extensions to the corresponding diffuse domains, that $\epsilon>0$, $\delta\geq 0$, 
and that the error of the initial data is $\mathcal{O}(\eps^{3/2} + \delta^{1/2})$.
Then, there exists a constant $C$ independent of $\eps$ and $\delta$ such that the following estimate holds:
\begin{multline*}
    \frac{1}{4}
    \left(     
    \rho_F \| \be_\bu^{\eps,\delta}(t)\|^2_{L^2(\Omega, \Phi_F^{\eps, \delta})^d}
    + \rho_B \| \partial_t {\be_\bfeta^{\eps, \delta}(t)}\|^2_{L^2(\Omega, \Phi_B^{\eps, \delta})^d}
    + c_0 \|{e_p^{\eps, \delta}(t)} \|^2_{L^2(\Omega, \Phi_B^{\eps, \delta})}
    + \|{\be_\bfeta^{\eps, \delta}(t)}\|^2_{E, \eps, \delta}
    \right) 
    \\
    + 2\mu_F \|\bD({\be_\bu^{\epsilon, \delta}})\|^2_{L^2(0, t; L^2(\Omega, \Phi_F^{\eps, \delta})^{d \times d})}
    + \| \boldsymbol\kappa^{1/2}\nabla {e_p^{\epsilon, \delta}}\|^2_{L^2(0, t; L^2(\Omega, \Phi_B^{\eps, \delta})^d)}
    \\
    + \alpha_{BJ} \sum_{i = 1}^{d - 1} \int_0^t \int_{\Omega}\left |( {\be_\bu^{\epsilon, \delta} - \partial_t \be_\bfeta^{\epsilon, \delta}})\cdot \tilde{\bftau}_i \right |^2 |\nabla \Phi_F^{\epsilon, \delta}|
    \ltt 
    (\eps^{3} + {\delta}) e^{t} b(t),
\end{multline*}
where $b$ depends on higher order norms of the problem data and sharp interface solution. 
\begin{confidential}
\begin{align*}
    &b(t) = \|\bu_s\|_{L^2(0, t; W^{3, \infty}(\Omega)^d)}^2 
    + \|\bu_s\|_{L^{\infty}(0, t; H^{3}(\Omega)^d)}^2
    + \|\bu_s\|_{H^1(0, t; W^{1, \infty}(\Omega)^d)}^2
    \\
    &+\|\pi_s\|_{L^2(0,t;W^{2, \infty}(\Omega)^d)}^2 
    + \|\pi_s\|_{L^{\infty}(0, T; H^2(\Omega))}^2
    + \|p_s\|_{L^2(0, t;W^{3, \infty}(\Omega)^d)}^2 
    + \|p_s\|_{L^{\infty}(0, T; W^{2, \infty}(\Omega))}^2
    \\
    &+ \| \bfeta_s \|_{L^{\infty}(0, T; W^{3, \infty}(\Omega)^d)}^2
    + \| \bfeta_s \|_{H^{3}(0, T; W^{1, \infty}(\Omega)^d)}^2
    + \| \bfeta_s \|_{H^2(0, T; H^2(\Omega)^d)}^2
    + \| \partial_t \bfeta_s \|_{L^{\infty}(0, T; H^2(\Omega)^d)}^2
    \\
    & + \|\bF_F\|_{L^2(0, t; H^1(\Omega)^d)}^2
    + \|g\|_{L^2(0, t; H^1(\Omega))}^2
    + \| \bF_B \|_{L^{\infty}(0, t; H^1(\Omega)^d)}.
    
\end{align*}
\end{confidential}
\end{theorem}
\begin{proof}
To simplify notation, we set $\chi = \chi_{\Omega_F}$. 
We start by subtracting~\eqref{eq:SBFP_weak} from~\eqref{eq:weakSI}. Using a divergence-free extension operator to extend the sharp interface solution to the whole domain $\Omega$, we obtain, upon assuming that $\nabla \cdot \bv = 0$, the following equation:
\begin{equation}
    \begin{aligned}\label{eq:modelerror}
        & \rho_F \int_0^t \int_{\Omega}\partial_t {\be_\bu^{\epsilon, \delta}}\cdot\bv\Phi_F^{\epsilon, \delta}
        + \rho_F \int_0^t \int_{\Omega}\partial_t\bu_s\cdot \bv (\chi - \Phi_F^{\epsilon, \delta})
        + \rho_B \int_0^t \int_{\Omega} \partial_{tt} {\be_\bfeta^{\epsilon, \delta} }\cdot \bfphi \Phi_B^{\epsilon, \delta}
        \\
        & \; + \rho_B \int_0^t \int_{\Omega}\partial_{tt} \bfeta_s\cdot \bfphi (\Phi_F^{\epsilon, \delta} - \chi)
        + c_0 \int_0^t \int_{\Omega}\partial_t {e_p^{\epsilon, \delta} } q\Phi_B^{\epsilon, \delta}
        + c_0 \int_0^t \int_{\Omega}\partial_t p_s q (\Phi_F^{\epsilon, \delta} - \chi)
        \\
        & \; + 2 \mu_F \int_0^t \int_{\Omega} \bD({\be_\bu^{\epsilon, \delta} }) : \bD(\bv)\Phi_F^{\epsilon, \delta}
        + 2 \mu_F \int_0^t \int_{\Omega} \bD(\bu_s) : \bD( \bv)(\chi - \Phi_F^{\epsilon, \delta})
        \\
        & \;   
        + 2 \mu_B \int_0^t \int_{\Omega} \bD({\be_\bfeta^{\epsilon, \delta}}) : \bD( \bfphi)\Phi_B^{\epsilon, \delta}
        + 2 \mu_B \int_0^t \int_{\Omega} \bD(\bfeta_s ) : \bD( \bfphi)(\Phi_F^{\epsilon, \delta} - \chi)
        \\
        & \;
        + \lambda_B \int_0^t \int_{\Omega} \nabla \cdot ({\be_\bfeta^{\epsilon, \delta}}) \nabla \cdot \bfphi \Phi_B^{\epsilon, \delta}
        + \lambda_B \int_0^t \int_{\Omega} \nabla \cdot \bfeta_s \nabla \cdot \bfphi (\Phi_F^{\epsilon, \delta} - \chi)
        + \int_0^t \int_{\Omega} \boldsymbol\kappa \nabla {e_p^{\epsilon, \delta} }\cdot \nabla q \Phi_B^{\epsilon, \delta}
        \\
        &\;
        + \int_0^t \int_{\Omega}\boldsymbol\kappa \nabla p_s\cdot\nabla q (\Phi_F^{\epsilon, \delta} - \chi)
        - \alpha \int_0^t \int_{\Omega} {e_p^{\epsilon,\delta} } \nabla\cdot\bfphi \Phi_B^{\epsilon, \delta}
        - \alpha \int_0^t \int_{\Omega} p_s \nabla\cdot\bfphi (\Phi_F^{\epsilon,\delta}-\chi)
        \\
        & \;
        + \alpha \int_0^t \int_{\Omega} q \nabla\cdot \partial_t {\be_\bfeta^{\epsilon, \delta} } \Phi_B^{\epsilon, \delta}
        + \alpha \int_0^t \int_{\Omega} q \nabla\cdot \partial_t \bfeta_s (\Phi_F^{\epsilon, \delta} - \chi) 
        \\
        & \;
        + \int_0^t \int_{\Gamma} q \partial_t \bfeta_s \cdot\bn 
        - \int_0^t \int_{\Gamma} q \bu_s \cdot\bn 
        - \int_0^t \int_{\Gamma} p_s \bfphi \cdot \bn 
        + \int_0^t \int_{\Gamma} p_s \bv \cdot \bn 
        \\
        & \;                    
        + \alpha_{BJ} \sum_{i = 1}^{d - 1} \int_0^t \int_{\Gamma} ( \bu_s - \partial_t \bfeta_s) \cdot\bftau_i (\bv - \bfphi)\cdot\bftau_i
        + \int_0^t \int_{\Omega} q (\partial_t \bfeta^{\epsilon, \delta}_d - \bu^{\epsilon, \delta}_d) \cdot \nabla\Phi_F^{\epsilon, \delta}
        \\
        & \;   
        - \int_0^t \int_{\Omega} p^{\epsilon, \delta}_{d} (\bfphi - \bv) \cdot \nabla\Phi_F^{\epsilon, \delta}
        - \alpha_{BJ} \sum_{i = 0}^{d - 1} \int_0^t \int_{\Omega}(\bu^{\epsilon, \delta}_d - \partial_t \bfeta^{\epsilon, \delta}_d )\cdot \tilde{\bftau}_i(\bv - \bfphi )\cdot \tilde{\bftau}_i |\nabla\Phi_F^{\epsilon, \delta}|  
        \\
        &
        = \int_0^t \int_{\Omega}\left (\rho_F\bF_F \cdot\bv - \rho_B\bF_B \cdot\bfphi - g q\right )(\chi - \Phi_F^{\eps, \delta}).
    \end{aligned}
\end{equation}
Taking $(\bv, \bfphi, q) = ({\be_\bu^{\epsilon, \delta}, \partial_t \be_\bfeta^{\epsilon, \delta}, e_p^{\epsilon, \delta}})$, which is admissible since $\nabla \cdot {\be_\bu^{\epsilon,\delta} } = 0$, we obtain, after integrating in time,
    \begin{align*}
        & \frac{1}{2}
        \left (     
        \rho_F \|\be_\bu^{\epsilon, \delta}(t)\|^2_{L^2(\Omega, \Phi_F^{\eps, \delta})^d}
        + \rho_B \|\partial_t \be_\bfeta^{\epsilon, \delta}(t)\|^2_{L^2(\Omega, \Phi_B^{\eps, \delta})^d}
        + c_0 \|e_p^{\epsilon, \delta}(t) \|^2_{L^2(\Omega, \Phi_B^{\eps, \delta})}
        + \| \be_\bfeta^{\epsilon, \delta}(t)\|^2_{E, \eps, \delta}
        \right)
        \\
        & + 2\mu_F \|\bD(\be_\bu^{\epsilon, \delta})\|^2_{L^2(0, t; L^2(\Omega, \Phi_F^{\eps, \delta})^{d \times d})}
        + \|\boldsymbol\kappa \nabla e_p^{\epsilon, \delta}\|^2_{L^2(0, t; L^2(\Omega, \Phi_B^{\eps, \delta})^d)}
        = \sum_{j = 1}^6 \cI_j,
    \end{align*}
where
    \begin{align*}
        \cI_1 &\coloneqq - \rho_F\int_0^t\int_{\Omega}\partial_t\bu_s\cdot \bv (\chi - \Phi_F^{\epsilon, \delta})
        - \rho_B \int_0^t \int_{\Omega}\partial_{tt} \bfeta_s\cdot \bfphi (\Phi_F^{\epsilon, \delta} - \chi)
        - c_0 \int_0^t \int_{\Omega}\partial_t p_s q (\Phi_F^{\epsilon, \delta} - \chi),
        \\
        \cI_2 &\coloneqq - 2 \mu_F \int_0^t \int_{\Omega}\bD(\bu_s):\bD( \bv)(\chi - \Phi_F^{\epsilon, \delta})
        - \int_0^t\int_{\Omega}\bfsigma_B(\bfeta_s, p_s):\nabla \bfphi (\Phi_F^{\epsilon, \delta} - \chi) 
        ,\\
        \cI_3 &\coloneqq - \alpha \int_0^t \int_{\Omega} q \nabla\cdot \partial_t \bfeta_s (\Phi_F^{\epsilon, \delta} - \chi)  
        - \int_0^t\int_{\Omega} \boldsymbol\kappa \nabla p_s\cdot\nabla q (\Phi_F^{\epsilon, \delta} - \chi)
        \\
        \cI_4 &\coloneqq - \int_0^t \int_{\Gamma} \left(q \left(\partial_t \bfeta_s - \bu_s \right) \cdot\bn - p_s \left(\bfphi - \bv \right)\cdot \bn \right),
        \\
        &- \int_0^t \int_{\Omega} \left(q \left(\partial_t \bfeta^{\epsilon, \delta}_d - \bu^{\epsilon, \delta}_d \right) \cdot \nabla\Phi_F^{\epsilon, \delta}
        - p^{\epsilon, \delta}_{d} \left(\bfphi - \bv \right) \cdot \nabla\Phi_F^{\epsilon, \delta} \right)
        \\
        \cI_5 &\coloneqq - \alpha_{BJ} \sum_{i = 0}^{d - 1} \int_0^t\left (\int_{\Gamma}(\bu_s - \partial_t \bfeta_s) \cdot\bftau_i (\bv - \bfphi) \cdot\bftau_i
        - \int_{\Omega}(\bu^{\epsilon, \delta}_d - \partial_t \bfeta^{\epsilon, \delta}_d )\cdot \tilde{\bftau}_i (\bv - \bfphi)\cdot \tilde{\bftau}_i |\nabla\Phi_F^{\epsilon, \delta}| \right )
    	,\\
        \cI_6 &\coloneqq \int_0^t\int_{\Omega}\left (\rho_F\bF_F \cdot\bv- \bF_B \cdot\bfphi - g q \right )(\chi - \Phi_F^{\eps, \delta}).
    \end{align*}

To estimate integrals $\cI_1$ to $\cI_6$, we use the convergence results for the diffuse integrals from~\cite[Section 5]{Burger2017}, and the following estimates 
which are consequences of the Poincar\'e and Korn inequalities
\begin{equation*}
    \|{\be_\bu^{\eps, \delta}}\|^2_{H^1(\Omega, \Phi_F^{\epsilon, \delta})^d}\ltt \int_{\Omega}|{\bD}({\be_\bu^{\eps, \delta}})|^2\Phi_F^{\epsilon, \delta},
    \qquad
    \| {e_p^{\eps, \delta}}\|^2_{H^1 (\Omega, \Phi_B^{\epsilon, \delta})}\ltt \int_{\Omega}|\nabla {e_p^{\eps, \delta}} |^2\Phi_B^{\epsilon, \delta}.
\end{equation*}
Note that since $e_p^{\eps, \delta}$ does not vanish on the boundary, we cannot directly use the Poincar\'e inequality as indicated above. However, we can bound its mean value as in the proof of Theorem~\ref{thm:WellPosedness}. For brevity we skip the details since the estimates are completely analogous. We also note that, because of~\eqref{eq:regularizacija}, we have, for $i \in\{F, B\}$, $\|\cdot \|_{L^2(\Omega, \Phi_i^{\eps})}\ltt \|\cdot \|_{L^2(\Omega,\Phi_i^{\eps, \delta})}$.

\framebox{Bound on $\cI_1$:} We integrate by parts in time the term containing $\bfphi = \partial_t {\be_\bfeta^{\epsilon, \delta}}$.
\begin{confidential}
as follows:
    \begin{align*}
        & \rho_B \int_{\Omega} \int_0^t \partial_{tt} \bfeta_s\cdot \bfphi (\Phi_F^{\epsilon} - \chi)
        \\
        &\qquad
        = - \rho_B \int_{\Omega}\int_0^t \partial_{ttt} \bfeta_s\cdot ( \bfeta_s - \bfeta_d^{\epsilon, \delta}) (\Phi_F^{\epsilon} - \chi)
        + \rho_B \left[\int_{\Omega} \partial_{tt} \bfeta_s \cdot ( \bfeta_s - \bfeta_d^{\epsilon, \delta}) (\Phi_F^{\epsilon} - \chi) \right]_0^t.
    \end{align*}
Integral $I_1$ is now estimated 
\end{confidential}
Using~\eqref{eq:regularizacija} and~\cite[Theorem 5.2]{Burger2017} for $p = 2$, we have:
    \begin{align*}
        |\cI_1| &\leq \left |\rho_F \int_0^t \int_{\Omega}\partial_t\bu_s\cdot \bv (\chi - \Phi_F^{\epsilon})
        - \rho_B\int_{\Omega}\int_0^t \partial_{ttt} \bfeta_s\cdot {\be_\bfeta^{\epsilon, \delta}} (\Phi_F^{\epsilon} - \chi)
        \right.
        \\
        & \qquad \left. + \rho_B \int_{\Omega} \partial_{tt} \bfeta_s(t) \cdot {\be_\bfeta^{\epsilon, \delta}}  (\Phi_F^{\epsilon} - \chi) 
        + c_0 \int_0^t \int_{\Omega}\partial_t p_s  q (\Phi_F^{\epsilon} - \chi)\right |  
        \\
        &
        + \delta \left |\rho_F \int_0^t \int_{\Omega}\partial_t\bu_s\cdot \bv (-1 + 2\Phi_F^{\epsilon})
        - \rho_B  \int_{\Omega}\int_0^t \partial_{ttt} \bfeta_s\cdot {\be_\bfeta^{\epsilon, \delta}} (1 - 2\Phi_F^{\epsilon})
        \right.
        \\
        & 
        \qquad \left.
        + \rho_B \int_{\Omega} \partial_{tt} \bfeta_s(t) \cdot {\be_\bfeta^{\epsilon, \delta}}(t) (1 - 2\Phi_F^{\epsilon})
        + c_0 \int_0^t \int_{\Omega}\partial_t p_s q (1 - 2\Phi_F^{\epsilon})
        \right |
        \\
        &\ltt 
        \epsilon^{3/2} \int_0^t \left (                      
        \|\bD({\be_\bu^{\epsilon, \delta}})\|_{L^2(\Omega, \Phi_F^{\epsilon})^{d \times d}}
        + 
        \|\bD({\be_\bfeta^{\epsilon, \delta}})\|_{L^2(\Omega, \Phi_B^{\epsilon})^{d \times d}}
        + 
        \|\nabla {e_p^{\epsilon, \delta} } \|_{L^2(\Omega,\Phi_B^{\epsilon})^d} \right) 
        \\ 
        & + \epsilon^{3/2}  
        \|\bD( {\be_\bfeta^{\epsilon, \delta}(t)}) \|_{L^2(\Omega, \Phi_B^{\epsilon})^{d \times d}} 
        + \cI_1^{\delta},
    \end{align*}
where $\cI_1^{\delta}$ denotes terms multiplied by $\delta$ and the implicit constant depends on higher order norms of the sharp interface solution.

To estimate $\cI_1^{\delta}$, we proceed as follows:
    \begin{align*}
        |\cI_1^{\delta}| &\ltt \int_0^t \left( \|\partial_t \bu_s \|_{L^{\infty}(\Omega)^d}\int_{\Omega}\delta |\bv| \right)
        + \int_0^t \left(\|\partial_{ttt} \bfeta_s \|_{L^{\infty}(\Omega)^d}\int_{\Omega}\delta | {\be_\bfeta^{\epsilon, \delta} }| \right)
        \\
        & \quad + \int_0^t \left(\|\partial_t p_s\|_{L^{\infty}(\Omega)} \int_{\Omega}\delta | q| \right)
        + \|\partial_{tt} \bfeta_s (t) \|_{L^{\infty}(\Omega)^d}\int_{\Omega}\delta |{\be_\bfeta^{\epsilon, \delta} } (t)| 
        \\
        &\ltt \sqrt{\delta} \int_0^t \left( 
        \|{\be_\bu^{\eps, \delta}}\|_{L^2(\Omega, \Phi_F^{\eps, \delta})^d} 
        + 
        \|{\be_\bfeta^{\epsilon, \delta} }\|_{L^2(\Omega, \Phi_B^{\eps, \delta})^d} 
        + 
        \|{ e_p^{\eps, \delta} }\|_{L^2(\Omega, \Phi_B^{\eps, \delta})} \right)
        + \sqrt{\delta}
        \| {\be_\bfeta^{\epsilon, \delta} }(t)\|_{L^2(\Omega, \Phi_B^{\eps, \delta})^d},
    \end{align*}
where the implicit constant depends on higher order norms of the sharp interface solution.
We note that all other terms in volume integrals that have a factor of $\delta$ can be estimated analogously. Finally, by using $\|\cdot\|_{L^2(\Omega, \Phi_F^{\eps})}\ltt \|\cdot\|_{L^2(\Omega, \Phi_F^{\eps, \delta})}$, the Cauchy--Schwarz inequality, and combining previous estimates we get:
    \begin{align*}
        |\cI_1| &\ltt (\eps^{3/2} + \sqrt{\delta})
        \left (\|\bD( {\be_\bu^{\epsilon, \delta} })\|_{L^2(0, t; L^2(\Omega, \Phi_F^{\epsilon, \delta})^{d \times d})}  
        + \|\nabla {e_p^{\epsilon, \delta} } \|_{L^2(0, t; L^2(\Omega, \Phi_B^{\epsilon, \delta})^d)}
        \right.
        \\
        & + \left. \|\bD({\be_\bfeta^{\epsilon, \delta}})\|_{L^2(0, t; L^2(\Omega, \Phi_B^{\epsilon})^{d \times d})}
        +
        \|\bD( {\be_\bfeta^{\epsilon, \delta}}(t))\|_{L^2(\Omega, \Phi_B^{\epsilon})^{d \times d}} \right).
    \end{align*}

\framebox{Bound on $\cI_6$:} In a similar way, we obtain the following bound
\begin{equation*}
    \frac{ |\cI_6| }{\epsilon^{3/2} + \sqrt{\delta} }
    \ltt 
    \|\bD( {\be_\bu^{\epsilon, \delta} })\|_{L^2(0, t; L^2(\Omega, \Phi_F^{\epsilon, \delta})^{d \times d})}
    + \|\nabla {e_p^{\epsilon, \delta} }\|_{L^2(0, t; L^2(\Omega, \Phi_B^{\epsilon, \delta})^d)}
    + \|\bD( {\be_\bfeta^{\epsilon, \delta} }(t))\|_{L^2(\Omega, \Phi_B^{\epsilon})^{d \times d}}.
\end{equation*}

\framebox{Bound on $\cI_2$:} 
We integrate by parts in space and use again the fact that $\bu_s$ and $\bu_d^{\epsilon, \delta}$ are divergence free.
\begin{confidential}
, obtaining the following:
\begin{align*}
        I_2& = \int_0^t \int_{\Omega}(\nabla\cdot\bfsigma_F(\bu_s, \pi_s))\cdot\bv(\chi - \Phi_F^{\epsilon, \delta})
        - \int_0^t\int_{\Omega}\bfsigma_F(\bu_s, \pi_s)\nabla\Phi_F^{\epsilon, \delta}\cdot\bv
        - \int_0^t\int_{\Gamma }\bfsigma_F(\bu_s, \pi_s)\bn\cdot\bv
        \\
        & + \int_0^t\int_{\Omega}\nabla\cdot\bfsigma_B(\bfeta_s, p_s)\cdot\bfphi (\Phi_F^{\epsilon, \delta} - \chi)
        + \int_0^t\int_{\Omega}\bfsigma_B(\bfeta_s, p_s)\nabla\Phi_F^{\epsilon, \delta}\cdot\bfphi
        + \int_0^t\int_{\Gamma}\bfsigma_B(\bfeta_s, p_s)\bn \cdot\bfphi
        \\
        &= \int_0^t\int_{\Omega}\nabla\cdot\bfsigma_F(\bu_s, \pi_s)\cdot\bv(\chi - \Phi_F^{\epsilon, \delta})
        - \int_0^t\int_{\Omega} \left(\bfsigma_F(\bu_s, \pi_s)\frac{\nabla\Phi_F^{\epsilon, \delta}}{|\nabla\Phi_F^{\epsilon, \delta}|}\cdot\frac{\nabla\Phi_F^{\epsilon, \delta}}{|\nabla\Phi_F^{\epsilon, \delta}|} \right)(\bv\cdot\nabla\Phi_F^{\epsilon, \delta}) 
        \\
        &-\sum_{i = 1}^{d - 1} \int_0^t\int_{\Omega} \left(\bfsigma_F(\bu_s, \pi_s) \nabla \Phi_F^{\epsilon, \delta} \cdot\tilde{\bftau}_i \right)(\bv\cdot \tilde{\bftau}_i )
        + \int_0^t\int_{\Gamma}p_s\bv\cdot\bn
        + \alpha_{BJ} \sum_{i = 1}^{d - 1} \int_0^t \int_{\Gamma}(\bu_s -\partial_t \bfeta_s)\cdot\bftau_i (\bv\cdot\bftau_i)
        \\
        &+ \int_0^t\int_{\Omega}\nabla\cdot\bfsigma_B(\bfeta_s, p_s) \cdot\bfphi (\Phi_F^{\epsilon, \delta} - \chi)
        + \int_0^t\int_{\Omega} \left(\bfsigma_B(\bfeta_s, p_s) \frac{\nabla\Phi_F^{\epsilon, \delta}}{|\nabla\Phi_F^{\epsilon, \delta}|}\cdot\frac{\nabla\Phi_F^{\epsilon, \delta}}{|\nabla\Phi_F^{\epsilon, \delta}|} \right)(\bfphi \cdot\nabla\Phi_F^{\epsilon, \delta}) 
        \\
        &
       + \sum_{i = 1}^{d - 1} \int_0^t\int_{\Omega} \left(\bfsigma_B(\bfeta_s, p_s) \nabla \Phi_F^{\epsilon, \delta} \cdot\tilde{\bftau}_i \right)(\bfphi \cdot \tilde{\bftau}_i )
        - \int_0^t \int_{\Gamma}p_s\bfphi \cdot\bn
        - \alpha_{BJ} \sum_{i = 1}^{d - 1}\int_0^t \int_{\Gamma}(\bu_s -\partial_t \bfeta_s)\cdot\bftau_i (\bfphi \cdot\bftau_i).
\end{align*}
\end{confidential}
Adding and subtracting $\sum_{i = 1}^{d - 1} \int_0^t \int_{\Omega} \bfsigma_F( \bfeta_s, p_s) \nabla\Phi_F^{\epsilon, \delta} \cdot\tilde{\bftau}_i (\bfphi \cdot \tilde{\bftau}_i )$, we get:
    \begin{align*}
        \cI_2& = \int_0^t \int_{\Omega}(\nabla\cdot\bfsigma_F(\bu_s, \pi_s))\cdot\bv(\chi - \Phi_F^{\epsilon, \delta})
        - \int_0^t \int_{\Omega} \left(\bfsigma_F(\bu_s, \pi_s)\frac{\nabla\Phi_F^{\epsilon, \delta}}{|\nabla\Phi_F^{\epsilon, \delta}|}\cdot\frac{\nabla\Phi_F^{\epsilon, \delta}}{|\nabla\Phi_F^{\epsilon, \delta}|} \right)(\bv\cdot\nabla\Phi_F^{\epsilon, \delta}) 
        \\
        &- \sum_{i = 1}^{d - 1} \int_0^t \int_{\Omega} \left(\bfsigma_F(\bu_s, \pi_s) \nabla \Phi_F^{\epsilon, \delta} \cdot\tilde{\bftau}_i \right)(\bv - \bfphi) \cdot \tilde{\bftau}_i 
        + \int_0^t\int_{\Gamma}p_s\bv\cdot\bn
        + \alpha_{BJ} \sum_{i = 1}^{d - 1} \int_0^t \int_{\Gamma}(\bu_s - \partial_t \bfeta_s)\cdot\bftau_i (\bv\cdot\bftau_i)
        \\
        &+ \int_0^t\int_{\Omega}\nabla\cdot\bfsigma_B(\bfeta_s, p_s) \cdot\bfphi (\Phi_F^{\epsilon, \delta} - \chi)
        + \int_0^t \int_{\Omega} \left(\bfsigma_B(\bfeta_s, p_s) \frac{\nabla\Phi_F^{\epsilon, \delta}}{|\nabla\Phi_F^{\epsilon, \delta}|}\cdot\frac{\nabla\Phi_F^{\epsilon, \delta}}{|\nabla\Phi_F^{\epsilon, \delta}|} \right)(\bfphi \cdot\nabla\Phi_F^{\epsilon, \delta}) 
        \\
        &
        + \sum_{i = 1}^{d - 1} \int_0^t \int_{\Omega} \left(\bfsigma_B(\bfeta_s, p_s) \nabla \Phi_F^{\epsilon, \delta} - \bfsigma_F(\bfeta_s, p_s) \nabla\Phi_F^{\epsilon, \delta} \right) \cdot\tilde{\bftau}_i (\bfphi \cdot \tilde{\bftau}_i )
        \\
        &
        - \int_0^t \int_{\Gamma}p_s\bfphi \cdot\bn
        - \alpha_{BJ}\sum_{i = 1}^{d - 1} \int_0^t \int_{\Gamma}(\bu_s -\partial_t \bfeta_s)\cdot\bftau_i (\bfphi \cdot\bftau_i).
    \end{align*}

\framebox{Bound on $\cI_3$:} We integrate by parts in space to obtain:
\begin{equation*}
    \cI_3 = \int_0^t \int_{\Omega} \left( \boldsymbol\kappa \Delta p_s q(\chi - \Phi_F^{\epsilon, \delta})
    + q\boldsymbol\kappa \nabla p_s\cdot \nabla\Phi_F^{\epsilon, \delta} \right)
    + \int_0^t \int_{\Gamma} q \boldsymbol\kappa \nabla p_s \cdot \bn
    - \alpha\int_0^t \int_{\Omega} q \nabla\cdot \partial_t \bfeta_s (\Phi_F^{\epsilon, \delta} - \chi),
\end{equation*}
and we recall that, owing to the coupling conditions, on $\Gamma$ we have
\begin{equation*}
    \boldsymbol\kappa \nabla p_s \cdot \bn = -\bu_s\cdot\bn + \partial_t \bfeta_s\cdot\bn.
\end{equation*}

\framebox{Bound on $\cI_4$:}  Observe that, by the choice of test functions we have
\begin{equation*}
    \int_0^t \int_{\Omega}e_p^{\varepsilon, \delta} (\bv - \bfphi ) \cdot\nabla\Phi_F^{\epsilon, \delta}
    - \int_0^t \int_{\Omega} q(\be_\bu^{\epsilon, \delta} - \partial_t \be_\bfeta^{\epsilon, \delta})\cdot\nabla\Phi_F^{\epsilon, \delta} = 0.
\end{equation*}
Therefore, adding and subtracting $\int_0^t \int_{\Omega}p_s(\bv - \bfphi) \cdot\nabla\Phi_F^{\epsilon, \delta}$ and $\int_0^t \int_{\Omega} q (\bu_s - \partial_t \bfeta_s) \cdot\nabla\Phi_F^{\epsilon, \delta}$, we obtain
\begin{equation*}
    \cI_4 = 
    \int_0^t \int_{\Omega} \left[ q (\bu_s - \partial_t \bfeta_s) \cdot\nabla\Phi_F^{\epsilon, \delta}
    - p_s(\bv - \bfphi) \cdot\nabla\Phi_F^{\epsilon, \delta} \right]
    + \int_0^t \int_{\Gamma} \left[ p_s (\bfphi - \bv)
    - q \left(\partial_t \bfeta_s - \bu_s \right) \right] \cdot\bn.
\end{equation*}

\framebox{Bound on $\cI_5$:} In a similar way, we obtain:
    \begin{align*}
        \cI_5 =& - \alpha_{BJ} \sum_{i = 1}^{d - 1} \int_0^t \int_{\Omega}\left |( {\be_\bu^{\epsilon, \delta} } - \partial_t {\be_\bfeta^{\epsilon, \delta} })\cdot \tilde{\bftau}_i \right |^2 |\nabla \Phi_F^{\epsilon, \delta}|
        \\
        & + \alpha_{BJ}\sum_{i = 1}^{d - 1} \int_0^t \int_{\Omega}(\bu_s - \partial_t \bfeta_s) \cdot \tilde{\bftau}_i (\bv - \bfphi) \cdot \tilde{\bftau}_i |\nabla\Phi_F^{\epsilon, \delta}|
        - \alpha_{BJ} \sum_{i = 1}^{d - 1} \int_0^t \int_{\Gamma}(\bu_s - \partial_t \bfeta_s) \cdot\bftau_i (\bv - \bfphi ) \cdot\bftau_i.
    \end{align*}

\framebox{Bound on $\sum_{j = 2}^5 \cI_j$:}
Adding the estimates for $\cI_2, \ldots, \cI_5$ together, we get the following:
    \begin{align*}
        \sum_{j = 2}^5 \cI_j &=
        - \alpha_{BJ} \sum_{i = 1}^{d - 1} \int_0^t \int_{\Omega}\left |( {\be_\bu^{\epsilon, \delta} - \partial_t \be_\bfeta^{\epsilon, \delta} })\cdot \tilde{\bftau}_i \right |^2 |\nabla \Phi_F^{\epsilon, \delta}|
        \\
        &
        +
        \int_0^t \int_{\Omega}\left ((\nabla\cdot\bfsigma_F(\bu_s, \pi_s))\cdot\bv - \nabla\cdot\bfsigma_B(\bfeta_s, p_s) \cdot\bfphi + \boldsymbol\kappa\Delta p_s q
        + \alpha q \nabla\cdot \partial_t \bfeta_s
        \right )(\chi - \Phi_F^{\epsilon, \delta})
        \\
        &
        + \int_0^t \int_{\Omega} q\left (\bu_s - \partial_t \bfeta_s + \boldsymbol\kappa\nabla p_s \right)\cdot\nabla\Phi_F^{\epsilon, \delta}
        - \int_0^t \int_{\Omega}\left (p_s + \bfsigma_F(\bu_s, \pi_s)\frac{\nabla\Phi_F^{\epsilon, \delta}}{|\nabla\Phi_F^{\epsilon, \delta}|}\cdot\frac{\nabla\Phi_F^{\epsilon, \delta}}{|\nabla\Phi_F^{\epsilon, \delta}|} \right )(\bv\cdot\nabla\Phi_F^{\epsilon, \delta})
        \\
        &
        + \int_0^t \int_{\Omega}\left (p_s + \bfsigma_B(\bfeta_s, p_s)\frac{\nabla\Phi_F^{\epsilon, \delta}}{|\nabla\Phi_F^{\epsilon, \delta}|}\cdot\frac{\nabla\Phi_F^{\epsilon, \delta}}{|\nabla\Phi_F^{\epsilon, \delta}|} \right )(\bfphi\cdot\nabla\Phi_F^{\epsilon, \delta})
        \\
        &+ \sum_{i = 1}^{d - 1} \int_0^t \int_{\Omega} \alpha_{BJ} (\bu_s - \partial_t \bfeta_s) \cdot \tilde{\bftau}_i (\bv - \bfphi)\cdot \tilde{\bftau}_i |\nabla\Phi_F^{\epsilon, \delta}|
        - \sum_{i = 1}^{d - 1} \int_0^t \int_{\Omega} \bfsigma_F(\bu_s, \pi_s) \nabla\Phi_F^{\epsilon, \delta} \cdot\tilde{\bftau}_i (\bv\cdot \tilde{\bftau}_i ) 
        \\
        &+ \sum_{i = 1}^{d - 1} \int_0^t \int_{\Omega} \bfsigma_B(\bfeta_s, p_s) \nabla\Phi_F^{\epsilon, \delta} \cdot\tilde{\bftau}_i (\bfphi \cdot \tilde{\bftau}_i ) = \sum_{j = 1}^6 \cD_j.
    \end{align*}
We see that $\cD_1 \leq 0$, so it will be combined with the terms on  the left hand side. Next, integration by parts in time and an application of~\cite[Theorem 5.2]{Burger2017} reveals that:
    \begin{align*}
        \frac{ |\cD_2| }{ \epsilon^{3/2} + \sqrt{\delta} }
        &\ltt 
        \|\bD( \be_\bu^{\epsilon, \delta})\|_{L^2(0, t; L^2(\Omega, \Phi_F^{\epsilon, \delta})^{d \times d})} 
        + \|\bD(\be_\bfeta^{\epsilon, \delta})\|_{L^2(0, t; L^2(\Omega, \Phi_B^{\epsilon, \delta})^{d \times d})} 
        +
        \| \bD(\be_\bfeta^{\epsilon, \delta}(t))\|_{L^2(\Omega, \Phi_B^{\epsilon, \delta})^{d \times d}} 
        \\
        & +
        \|\nabla e_p^{\epsilon, \delta} \|_{L^2(0, t; L^2(\Omega, \Phi_B^{\epsilon, \delta})^d)},
    \end{align*}
where, again, the implicit constant depends on higher order norms of the sharp interface solution. 
To estimate the remaining terms first we notice that in the tubular neighbourhood of $\Gamma$ we have $\bn = -\displaystyle\frac{\nabla\Phi_F^{\epsilon, \delta}}{|\nabla\Phi_F^{\epsilon, \delta}|} = -\displaystyle\frac{\nabla\Phi_F^{\epsilon}}{|\nabla\Phi_F^{\epsilon}|}$ and 
$\bftau_i = { -\tilde{\bftau}_i}, i = 1, \ldots, d - 1$. Therefore, using the coupling conditions, 
we are able to infer that, on $\Gamma$,
\begin{equation*}
    \left (\bu_s - \partial_t \bfeta_s + \boldsymbol\kappa\nabla p_s \right )\cdot\nabla\Phi_F^{\epsilon, \delta} = 0,
    \quad
    p_s + \bfsigma_F(\bu_s, \pi_s)\frac{\nabla\Phi_F^{\epsilon, \delta}}{|\nabla\Phi_F^{\epsilon, \delta}|}\cdot\frac{\nabla\Phi_F^{\epsilon, \delta}}{|\nabla\Phi_F^{\epsilon, \delta}|} = 0, 
\end{equation*}
and that
\begin{equation*}
    \alpha_{BJ} (\bu_s - \partial_t \bfeta_s) \cdot \tilde{\bftau}_i
    -\bfsigma_F(\bu_s, \pi_s) \frac{\nabla\Phi_F^{\epsilon, \delta}}{|\nabla\Phi_F^{\epsilon, \delta}|} \cdot\tilde{\bftau}_i = 0,
    \quad
    \bfsigma_B(\bfeta_s, p_s) \frac{\nabla\Phi_F^{\epsilon, \delta}}{|\nabla\Phi_F^{\epsilon, \delta}|}  
    -\bfsigma_F(\bu_s, \pi_s) \frac{\nabla\Phi_F^{\epsilon, \delta}}{|\nabla\Phi_F^{\epsilon, \delta}|} = 0.
\end{equation*}
This, combined with~\cite[Theorem 5.6]{Burger2017}, allows us to conclude that
    \begin{align*}
        \sum_{j = 3}^6 |\cD_j| &\lesssim \eps^{3/2} \left( 
        \| e_p^{\eps,\delta} \|_{L^2(0, t; H^1(\Omega, \Phi_B^{\eps, \delta}))} 
        + \| \bD(\be_\bu^{\eps, \delta}) \|_{L^2(0, t; L^2(\Omega, \Phi_F^{\eps, \delta})^{d \times d}) }
        + \| \bD(\be_\bfeta^{\eps, \delta}) \|_{L^2(0, t; L^2(\Omega, \Phi_B^{\eps, \delta})^{d \times d}) }
        \right. \\
        &+ \left. \| \bD(\be_\bfeta^{\eps, \delta}(t)) \|_{L^2(\Omega, \Phi_B^{\eps, \delta})^{d \times d} }
        \right).
    \end{align*}
where we also integrated by parts, in time, the terms that contain $\bfphi = \partial_t \be_\bfeta^{\eps,\delta}$.

\framebox{Final estimate:}
Combining all the obtained estimates, using Young's inequality and Gr\"onwall's inequality, we get the following error estimate:
    \begin{align*}
        & \frac{1}{4}
        \left (     
        \rho_F \|\be_\bu^{\epsilon, \delta}(t)\|^2_{L^2(\Omega, \Phi_F^{\eps, \delta})^d}
        + \rho_B \|\partial_t \be_\bfeta^{\epsilon, \delta}(t)\|^2_{L^2(\Omega, \Phi_B^{\eps, \delta})^d}
        + c_0 \|e_p^{\epsilon, \delta}(t) \|^2_{L^2(\Omega, \Phi_B^{\eps, \delta})}
        + \| \be_\bfeta^{\epsilon, \delta}(t)\|^2_{E, \eps, \delta}
        \right)
        \\
        & \; + \|\bD(\be_\bu^{\epsilon, \delta})\|^2_{L^2(0, t; L^2(\Omega, \Phi_F^{\eps, \delta})^{d \times d})}
        + \| \boldsymbol\kappa^{1/2} \nabla e_p^{\epsilon, \delta} \|^2_{L^2(0, t; L^2(\Omega, \Phi_B^{\eps, \delta})^d)} 
        \\
        & + \alpha_{BJ} \sum_{i = 1}^{d - 1} \int_0^t \int_{\Omega}\left |( {\be_\bu^{\epsilon, \delta} - \partial_t \be_\bfeta^{\epsilon, \delta}})\cdot \tilde{\bftau}_i \right |^2 |\nabla \Phi_F^{\epsilon, \delta}|
        \ltt
        \eps^3 + \delta + \| \bD( \be_\bfeta^{\eps, \delta} ) \|_{L^2(0, t; L^2(\Omega, \Phi_B^\eps)^{d \times d})}^2.
    \end{align*}
\begin{confidential}
Using Gr\"onwall's inequality, we obtain the following error estimate:
    \begin{align*}
        \frac{1}{4}
        \left (     
        \rho_F \|\bu_s(t) - \bu^{\epsilon, \delta}_d(t)\|^2_{L^2(\Omega, \Phi_F^{\eps, \delta})^d}
        + \rho_B \|\partial_t \bfeta_s(t) - \partial_t \bfeta^{\epsilon, \delta}_d(t)\|^2_{L^2(\Omega, \Phi_B^{\eps, \delta})^d}
        + c_0 \|p_s(t) - p^{\epsilon, \delta}_{d}(t) \|^2_{L^2(\Omega, \Phi_B^{\eps, \delta})}
        \right.
        \\
        \left. + 2\mu_B \|\bD(\bfeta_s(t) - \bfeta_d^{\epsilon, \delta}(t))\|^2_{L^2(\Omega, \Phi_B^{\eps, \delta})^{d \times d}}
        + \lambda_B \|\nabla \cdot(\bfeta_s(t) - \bfeta_d^{\epsilon, \delta}(t))\|^2_{L^2(\Omega, \Phi_B^{\eps, \delta})}
        \right)
        \\
        + 2\mu_F \|\bD(\bu_s - \bu_d^{\epsilon, \delta})\|^2_{L^2(0, t; L^2(\Omega, \Phi_F^{\eps, \delta})^{d \times d})}
        + \|\boldsymbol\kappa \nabla(p_s - p^{\epsilon, \delta}_d)\|^2_{L^2(0, t;  L^2(\Omega, \Phi_B^{\eps, \delta})^d)}
        \\
        + \alpha_{BJ} \sum_{i = 1}^{d - 1} \int_0^t \int_{\Omega}\left |(\bu_s - \bu_d^{\epsilon, \delta} - \partial_t \bfeta_s + \bfxi_d^{\epsilon, \delta})\cdot \tilde{\bftau}_i \right |^2 |\nabla \Phi_F^{\epsilon, \delta}|
        \\
        \ltt 
        (\eps^{3} + {\delta}) e^{t} b(t)
    \end{align*}
where
    \begin{align*}
        b(t) = \|\partial_t\bu_s\|_{L^2(0, t; W^{1, \infty}(\Omega)^d)}^2 + \|\partial_t p_s \|_{L^2(0, t; W^{1, \infty}(\Omega))}^2
        + \|\partial_{tt} \bfeta_s (t)\|_{W^{1, \infty}(\Omega)^d}^2
        + \| \bu_s\|_{L^2(0, t; W^{3, \infty}(\Omega)^d)}^2 
        \\
        + \|\pi_s\|_{L^2(0, t; W^{2, \infty}(\Omega)^d)}^2 
        + \| \bfeta_s (t) \|_{W^{3, \infty}(\Omega)^d}^2
        + \| p_s(t) \|_{W^{2, \infty}(\Omega)^d}^2
        + \|\bF_F\|_{L^2(0, t; H^1(\Omega)^d)}^2
        + \|g\|_{L^2(0, t; H^1(\Omega))}^2
        \\
        + \|\partial_t \bF_B\|_{L^2(0, t; H^1(\Omega)^d)} 
        + \| \bF_B (t) \|_{L^2(0, t; H^1(\Omega)^d)}
        + \|\bu_s\|_{L^2(0, t; H^3(\Omega)^d)}^2
        + \| \bfeta_s\|_{H^1(0, t; H^2(\Omega)^d)}^2
        + \|\pi_s\|_{L^2(0, t; H^2(\Omega))}^2
        \\
        + \|p_s\|_{L^2(0, t; H^3(\Omega))}^2
        + \|p_s\|_{L^2(0, t; W^{3, \infty}(\Omega))}^2
        + \| \bfeta_s \|_{H^1(0, t; W^{2, \infty}(\Omega))}^2
        + \|\bu_s (t)\|_{H^3(\Omega)^d}^2
        + \| \partial_t \bfeta_s (t)\|_{H^2(\Omega)^d}^2
        \\
        + \|\bfeta_s (t)\|_{H^3(\Omega)^d}^2
        + \|\pi_s(t)\|_{H^2(\Omega)}^2
        + \|p_s(t)\|_{H^2(\Omega)}^2
        + \|\bfeta_s\|_{H^3(0, t;W^{1, \infty}(\Omega)^d)}^2
        + \| \bfeta_s \|_{H^1(0, t; W^{3, \infty}(\Omega)^d)}^2 
        \\
        + \|p_s\|_{H^1(0, t; W^{2, \infty}(\Omega)^d)}^2 
        + \|\bu_s\|_{H^1(0, t; H^3(\Omega)^d)}^2 + \|\bfeta_s\|_{H^2(0, t; H^2(\Omega)^d)}^2
        + \|\bfeta_s\|_{H^1(0, t; H^3(\Omega)^d)}^2
        \\
        + \|\pi_s\|_{H^1(0, t; H^2(\Omega))}^2
        + \|p_s\|_{H^1(0, t; H^2(\Omega))}^2.
    \end{align*}

This can be simplified as
    \begin{align*}
        &b(t) = \|\bu_s\|_{L^2(0, t;W^{3, \infty}(\Omega)^d)}^2 
        + \|\bu_s\|_{L^{\infty}(0, t; H^{3}(\Omega)^d)}^2
        + \|\bu_s\|_{H^1(0, t; W^{1, \infty}(\Omega)^d)}^2
        \\
        &+ \|\pi_s\|_{L^2(0, t; W^{2, \infty}(\Omega)^d)}^2 
        + \|\pi_s\|_{L^{\infty}(0, T; H^2(\Omega))}^2
        + \|p_s\|_{L^2(0, t; W^{3, \infty}(\Omega)^d)}^2 
        + \|p_s\|_{L^{\infty}(0, T; W^{2, \infty}(\Omega))}^2
        \\
        &+ \| \bfeta_s \|_{L^{\infty}(0, T; W^{3, \infty}(\Omega)^d)}^2
        + \| \bfeta_s \|_{H^{3}(0, T; W^{1, \infty}(\Omega)^d)}^2
        + \| \bfeta_s \|_{H^2(0, T; H^2(\Omega)^d)}^2
        + \| \partial_t \bfeta_s \|_{L^{\infty}(0, T; H^2(\Omega)^d)}^2
        \\
        & + \|\bF_F\|_{L^2(0, t; H^1(\Omega)^d)}^2
        + \|g\|_{L^2(0, t; H^1(\Omega))}^2
        + \| \bF_B \|_{L^{\infty}(0, t; H^1(\Omega)^d)}.
    \end{align*}
\end{confidential}
\end{proof}

\subsection{Power distance weights}

Here, under the sole assumption that our weight is given by~\eqref{eq:defofPhiPsi}, we estimate the difference between the continuous solutions of the diffuse domain  and sharp interface formulations in terms of the parameter $\eps$ describing the diffuse interface width. The main differences with the previous section are, first, that we do not assume that our weight is positive, i.e.~$\delta = 0$, and that we cannot differentiate the weight. Indeed, for suitable values of $x$, we obtain
\begin{equation*}
    \nabla \Phi_F^\eps(x) = \frac1{2\eps} \cS'\left( \frac{\sdist_\Gamma(x)}\eps \right) \nabla \sdist_\Gamma(x) 
    = \frac\beta{2\eps} \left( \frac{\sdist_\Gamma(x)}\eps + 1 \right)^{\beta - 1} \nabla \sdist_\Gamma(x).
\end{equation*}
Since $\beta \in (0,1)$, this quantity is not bounded; see Figure~\ref{fig:DomainsAndLayers}. To overcome this, we will use the fact that $\Phi_i^\eps \in A_2$, for $i \in \{F,B\}$, at the expense of a slower rate of convergence in terms of $\eps$.

We now begin with the estimate. We begin by, for simplicity, assuming that $\boldsymbol\kappa$ is a constant scalar, and that we have at hand suitable extensions of the sharp interface solutions to the corresponding diffuse domains. Next, we test each problem with the difference of solutions and subtract the result. Denote
\begin{equation*}
    \be_\bu^\eps \coloneqq \bu_s - \bu_d^\eps, \qquad \be_\bfeta^\eps \coloneqq \bfeta_s - \bfeta_d^\eps, \qquad e_p^\eps \coloneqq p_s - p_d^\eps.
\end{equation*}
In much similarity to~\eqref{eq:modelerror} we obtain
\begin{multline*}
    \frac12 \frac{d}{dt} \left(
    \rho_F \| \be_\bu^\eps \|_{L^2(\Omega_F^\eps, \Phi_F^\eps)^d}^2 + \rho_B \| \partial_t \be_\bfeta^\eps \|_{L^2(\Omega_B^\eps, \Phi_B^\eps)^d}^2 + c_0 \| e_p^\eps \|_{L^2(\Omega_B^\eps, \Phi_B^\eps)}^2 
    + \| \be_\bfeta^\eps \|_{E, \eps}^2 
    \right) \\
    + 2\mu_F \| \bD(\be_\bu^\eps) \|_{L^2(\Omega_F^\eps, \Phi_F^\eps)^{d \times d}}^2 
    + \boldsymbol\kappa \|\nabla e_p^\eps \|_{L^2(\Omega_B^\eps, \Phi_B^\eps)^d}^2 = \sum_{j = 1}^4 \cR_j,
\end{multline*}
where
    \begin{align*}
        \cR_1 &=
        \rho_F \left( \int_{\ell_B^\eps} \partial_t \bu_s \cdot \be_\bu^\eps \Phi_F^\eps + 
        \int_{\ell_F^\eps} \partial_t \bu_s  \cdot \be_\bu^\eps \left( \Phi_F^\eps - \chi_F \right) \right) 
        \nonumber\\
        &+ \rho_B \left( \int_{\ell_F^\eps} \partial_{tt} \bfeta_s \cdot \partial_t \be_\bfeta^\eps \Phi_B^\eps 
        + \int_{\ell_B^\eps} \partial_{tt} \bfeta \cdot \partial_t \be_\bfeta^\eps \left(\Phi_B^\eps - \chi_B \right) \right)
        + c_0 \left( \int_{\ell_F^\eps} \partial_{t} p_s e_p^\eps \Phi_B^\eps 
        + \int_{\ell_B^\eps} \partial_t p_s e_p^\eps \left(\Phi_B^\eps - \chi_B \right) \right)
        \nonumber\\
        &+ \alpha \left( \int_{\ell_F^\eps} \nabla \cdot \partial_t \bfeta_s e_p^\eps\Phi_B^\eps 
        + \int_{\ell_B^\eps} \nabla \cdot \partial_t \bfeta_s e_p^\eps \left( \Phi_B^\eps - \chi_B \right) 
        - \int_{\ell_F^\eps} p_s \nabla \cdot \partial_t \be_\bfeta^\eps \Phi_B^\eps 
        - \int_{\ell_B^\eps} p_s \nabla \cdot \partial_t \be_\bfeta^\eps \left( \Phi_B^\eps - \chi_B \right) \right), 
        \nonumber\\
        \cR_2 &= 2\mu_F \left( \int_{\ell_B^\eps} \bD(\bu_s) : \bD(\be_\bu^\eps) \Phi_F^\eps 
        + \int_{\ell_F^\eps} \bD(\bu_s):\bD( \be_\bu^\eps ) \left( \Phi_F^\eps - \chi_F \right) \right) 
        + \int_{\ell_F^\eps} \bfsigma_E(\bfeta_s) : \bD( \partial_t \be_\bfeta^\eps ) \Phi_B^\eps 
        \nonumber\\
        &+ \int_{\ell_B^\eps} \bfsigma_E(\bfeta_s) : \bD( \partial_t \be_\bfeta^\eps ) \left( \Phi_B^\eps - \chi_B \right) 
        + \left( \int_{\ell_F^\eps} \boldsymbol\kappa \nabla p_s \cdot \nabla e_p^\eps \Phi_B^\eps 
        + \int_{\ell_B^\eps} \nabla p_s \cdot \nabla e_p^\eps \left( \Phi_B^\eps - \chi_B \right) \right),
        \\
        \cR_3 &= \int_\Gamma e_p^\eps \left( \bu_s - \partial_t \bfeta_s \right)\cdot \bn + 
        \frac1{2\eps} \int_{\ell^\eps} e_p^\eps \left( \bu_d^\eps - \partial_t \bfeta_d^\eps \right) \cdot \nabla \dist_\Gamma 
        \nonumber\\
        &- \int_\Gamma p_s \left( \be_\bu^\eps - \partial_t \be_\bfeta^\eps \right) \cdot \bn +
        \frac1{2\eps} \int_{\ell^\eps} p_d^\eps \left( \be_\bu^\eps - \partial_t \be_\bfeta^\eps \right) \cdot \nabla \dist_\Gamma 
        \nonumber\\
        &+ \alpha_{BJ} \sum_{i = 1}^{d - 1} \left( 
        \frac1{2\eps} \int_{\ell^\eps} \left( (\bu_d^\eps - \partial_t \bfeta_d^\eps) \cdot \tilde{\bftau}_i \right)
        \left( (\be_\bu^\eps - \partial_t \be_\bfeta^\eps) \cdot \tilde{\bftau}_i \right)
        - \int_\Gamma \left( (\bu_s - \partial_t \bfeta_s) \cdot \bftau_i \right)
        \left( (\be_\bu^\eps - \partial_t \be_\bfeta^\eps) \cdot \bftau_i \right)
        \right),
        \nonumber\\
        \cR_4 &= \rho_F \left( \int_{\ell_B^\eps} \bF_F \cdot \be_\bu^\eps \Phi_F^\eps 
        + \int_{\ell_F^\eps} \bF_F \cdot \be_\bu^\eps \left( \Phi_F^\eps - \chi_F \right) \right) 
        + \rho_B \left( \int_{\ell_F^\eps} \bF_B \cdot \partial_t \be_\bfeta^\eps \Phi_B^\eps 
        + \int_{\ell_B^\eps} \bF_B \cdot \partial_t \be_\bfeta^\eps \left( \Phi_B^\eps - \chi_B \right) \right) 
        \nonumber\\
        &+ \int_{\ell_F^\eps} g e_p^\eps \Phi_B^\eps + \int_{\ell_F^\eps} g e_p^\eps \left( \Phi_B^\eps - \chi_B \right),
        \nonumber
    \end{align*}
and $\chi_i$, for $i \in \{F, B\}$, denotes the characteristic function of $\Omega_i$.

The core issue at hand now is to bound the residual terms. The ideas draw heavy inspiration from the previous section,~\cite[Section 5]{Burger2017}, and~\cite[Theorem 5.5]{Bukac2023}.

\begin{theorem}[modelling error II]\label{thm:ModelErrorDistanceWeights}
Assume that the sharp interface solution and its forcing terms are sufficiently smooth, along with their extensions to the corresponding diffuse domains. Assume, in addition, that both $\{\bftau_i\}_{i = 1}^{d - 1}$ and $\{\tilde{\bftau}_i\}_{i = 1}^{d - 1}$ have suitable extensions onto $\Omega$ which belong to $C^1(\overline{\Omega_F^\eps})$. Then, for every $t \in (0,T)$, we have
\begin{multline*}
    \left(
    \rho_F \| \be_\bu^\eps \|_{L^2(\Omega_F^\eps, \Phi_F^\eps)^d}^2 + \rho_B \| \partial_t \be_\bfeta^\eps \|_{L^2(\Omega_B^\eps, \Phi_B^\eps)^d}^2 + c_0 \| e_p^\eps \|_{L^2(\Omega_B^\eps, \Phi_B^\eps)}^2 
    + \| \be_\bfeta^\eps \|_{E,\delta}^2 
    \right)(t) \\
    + \int_0^t \left( \| \bD(\be_\bu^\eps) \|_{L^2(\Omega_F^\eps, \Phi_F^\eps)^{d \times d}}^2 
    + \|\nabla e_p^\eps \|_{L^2(\Omega_B^\eps, \Phi_B^\eps)^d}^2 \right)
    \lesssim \eps^{1/2}e^{t}b(t),
\end{multline*}
where the implicit constant depends on the smoothness assumptions, but not on $\eps$. The coefficient $b$ depends on higher order norms of the sharp interface solution $(\bu_s, \bfeta_s,p_s)$ and its forcing terms $(\bF_F, \bF_B,g)$.
\end{theorem}
\begin{proof}
The proof entails estimating the terms $\{\cR_j\}_{j = 1}^4$ in a much similar way to the proof of~\cite[Theorem 5.5]{Bukac2023}. For brevity, we provide minimal details.

\framebox{Bound on $\cR_1$:} Using that, for $i \in \{F, B\}$, $0 \leq \Phi_i^\eps \leq 1$, that $\Phi_F^\eps \geq \tfrac12$ on $\ell_F^\eps$ we obtain
\begin{align*}
        \cR_1 &\leq C \left(\int_{\ell^\eps} \left( 
        |\partial_t \bu_s|^2 + |\partial_{tt} \bfeta_s|^2 + |\partial_t p_s |^2 + |\nabla\cdot \partial_t \bfeta_s|^2
        + |p_s|^2
        \right)\right) \\
        &+ \frac{\mu_F}4 \|\bD(\be_\bu^\eps)\|_{L^2(\Omega_F^\eps, \Phi_F^\eps)^{d \times d} }^2
        + \frac{\rho_B}4 \| \partial_t \be_\bfeta^\eps \|_{L^2(\Omega_B^\eps, \Phi_B^\eps)^d}^2 
        + \frac{ c_0 }4 \| e_p^\eps \|_{L^2(\Omega_B^\eps, \Phi_B^\eps)}^2 \\
        &\leq C \eps + \frac{\mu_F}4 \|\bD(\be_\bu^\eps)\|_{L^2(\Omega_F^\eps, \Phi_F^\eps)^{d \times d} }^2
        + \frac{\rho_B}4 \| \partial_t \be_\bfeta^\eps \|_{L^2(\Omega_B^\eps, \Phi_B^\eps)^d}^2 
        + \frac{ c_0 }4 \| e_p^\eps \|_{L^2(\Omega_B^\eps, \Phi_B^\eps)}^2,
    \end{align*}
where, when convenient, we used the weighted Korn's inequality \eqref{eq:weightedKorn}.

\framebox{Bound on $\cR_2$:} Similarly,
    \begin{align*}
        \cR_2 &\leq C \left(\int_{\ell^\eps} \left( 
        |\bD(\bu_s)|^2 + |\partial_{tt} \bfeta_s|^2 + |\partial_t p_s |^2 + |\nabla\cdot \partial_t \bfeta_s|^2
        + |p_s|^2 + |\nabla p_s |^2 
        \right)\right) 
        \nonumber\\
        &+ \frac{\mu_F}4 \| \bD(\be_\bu^\eps) \|_{L^2(\Omega_F^\eps, \Phi_F^\eps)^{d \times d}}^2
        + \frac{\boldsymbol\kappa}4 \| \nabla e_p^\eps \|_{L^2(\Omega_B^\eps, \Phi_B^\eps)^d}^2 
        \nonumber\\
        &+ \int_{\ell_F^\eps} \bfsigma_E(\bfeta_s) : \bD( \partial_t \be_\bfeta^\eps) \Phi_B^\eps 
        + \int_{\ell_B^\eps} \bfsigma_E(\bfeta_s) : \bD( \partial_t \be_\bfeta^\eps) \left( \Phi_B^\eps - \chi_B \right) 
        \\
        &\leq C\eps
        + \frac{\mu_F}4 \| \bD(\be_\bu^\eps) \|_{L^2(\Omega_F^\eps, \Phi_F^\eps)^{d \times d}}^2
        + \frac{\boldsymbol\kappa}4 \| \nabla e_p^\eps \|_{L^2(\Omega_B^\eps, \Phi_B^\eps)^d}^2 
        \nonumber\\
        &+ \int_{\ell_F^\eps} \bfsigma_E(\bfeta_s) : \bD( \partial_t \be_\bfeta^\eps) \Phi_B^\eps 
        + \int_{\ell_B^\eps} \bfsigma_E(\bfeta_s) : \bD( \partial_t \be_\bfeta^\eps) \left( \Phi_B^\eps - \chi_B \right).
        \nonumber
    \end{align*}

\framebox{Bound on $\cR_4$:} Next, we have
\begin{align*}
        \cR_4 &\leq C \left(\int_{\ell^\eps} \left( 
        |\bF_F|^2 + |\bF_B|^2 + |g|^2
        \right)\right) \\
        &+ \frac{\mu_F}4 \| \bD(\be_\bu^\eps) \|_{L^2(\Omega_F, \Phi_F^\eps)^{d \times d} }^2 
        + \frac{\rho_B}4 \| \partial_t \be_\bfeta^\eps \|_{L^2(\Omega_B^\eps, \Phi_B^\eps)^d}^2 
        + \frac{c_0}4 \| e_p^\eps \|_{L^2(\Omega_B^\eps, \Phi_B^\eps)}^2 \\
        &\leq C\eps 
        + \frac{\mu_F}4 \| \bD(\be_\bu^\eps) \|_{L^2(\Omega_F, \Phi_F^\eps)^{d \times d} }^2 
        + \frac{\rho_B}4 \| \partial_t \be_\bfeta^\eps \|_{L^2(\Omega_B^\eps, \Phi_B^\eps)^d}^2 
        + \frac{c_0}4 \| e_p^\eps \|_{L^2(\Omega_B^\eps, \Phi_B^\eps)}^2.
    \end{align*}

\framebox{Bound on $\cR_3$:} We leave this term last, as it warrants some more explanation. We first define
\begin{equation*}
    \cR_3 = \cR_{3, 1} + \cR_{3, 2},
\end{equation*}
with
\begin{equation*}
    \begin{aligned}
        \cR_{3, 1} &= \int_\Gamma e_p^\eps \left( \bu_s - \partial_t \bfeta_s \right)\cdot \bn + 
        \frac1{2\eps} \int_{\ell^\eps} e_p^\eps \left( \bu_d^\eps - \partial_t \bfeta_d^\eps \right) \cdot \nabla \dist_\Gamma \\
        & - \int_\Gamma p_s \left( \be_\bu^\eps - \partial_t \be_\bfeta^\eps \right) \cdot \bn +
        \frac1{2\eps} \int_{\ell^\eps} p_d^\eps \left( \be_\bu^\eps - \partial_t \be_\bfeta^\eps \right) \cdot \nabla \dist_\Gamma
    \end{aligned}
\end{equation*}
and
\begin{equation*}
    \cR_{3,2} = \alpha_{BJ} \sum_{i = 1}^{d - 1} \left( 
    \frac1{2\eps} \int_{\ell^\eps} \left( (\bu_d^\eps - \partial_t \bfeta_d^\eps) \cdot \tilde{\bftau}_i \right)
    \left( (\be_\bu^\eps - \partial_t \be_\bfeta^\eps) \cdot \tilde{\bftau}_i \right)
    - \int_\Gamma \left( (\bu_s - \partial_t \bfeta_s) \cdot \bftau_i \right)
    \left( (\be_\bu^\eps - \partial_t \be_\bfeta^\eps) \cdot \bftau_i \right)
    \right).
\end{equation*}
Next, we observe that
\begin{equation*}
    \begin{aligned}
        \cR_{3, 1} &= \left(
        \int_\Gamma e_p^\eps (\be_\bu^\eps - \partial_t \be_\bfeta^\eps )\cdot \bn 
        - \frac1{2\eps} \int_{\ell^\eps} e_p^\eps (\be_\bu^\eps - \partial_t \be_\bfeta^\eps )\cdot \nabla \dist_\Gamma
        \right) \\
        &+ \left(
        \frac1{2\eps} \int_{\ell^\eps} p_s (\be_\bu^\eps - \partial_t \be_\bfeta^\eps )\cdot \nabla \dist_\Gamma - \int_\Gamma p_s (\be_\bu^\eps - \partial_t \be_\bfeta^\eps )\cdot \bn
        \right).
    \end{aligned}
\end{equation*}
We argue as in the proof of~\cite[Theorem 5.5]{Bukac2023} to then conclude that
\begin{equation*}
    \cR_{3, 1} \lesssim \eps^{1/2}.
\end{equation*}

Finally,
    \begin{align*}
        \cR_{3, 2} &= - \frac{\alpha_{BJ}}{2\eps} \sum_{i = 1}^{d - 1} \int_{\ell^\eps} \left( (\be_\bu^\eps - \partial_t \be_\bfeta^\eps) \cdot \tilde{\bftau}_i \right)^2 
        \nonumber\\
        &+ \alpha_{BJ} \sum_{i = 1}^{d - 1} \left(
        \frac1{2\eps} \int_{\ell^\eps} \left( (\bu_s - \partial_t \bfeta_s) \cdot \tilde{\bftau}_i \right)
        \left( (\be_\bu^\eps - \partial_t \be_\bfeta^\eps) \cdot \tilde{\bftau}_i \right)
        - \int_\Gamma \left( (\bu_s - \partial_t \bfeta_s) \cdot \bftau_i \right)
        \left( (\be_\bu^\eps - \partial_t \be_\bfeta^\eps) \cdot \bftau_i \right)
        \right) 
        \nonumber\\
        &\lesssim \eps^{1/2},
    \end{align*}
where we dropped the negative term, and each one of the terms inside the sum were estimated, under the assumption that tangents have $C^1$ extensions, as in~\cite[Theorem 5.5]{Bukac2023}.

Having estimated all the terms, it is sufficient to apply Gr\"onwall's inequality. The result expediently follows.
\end{proof}

\section{Numerical results}\label{sec:numerics}

In this section, we illustrate the accuracy of our diffuse interface approach and further explore its capabilities with a series of numerical illustrations. All of our computations were carried out using the finite element 
library 
FreeFem++~\cite{Hecht2012}. For implementation, all integrals in the numerical scheme are formulated over the entire domain $\Omega$. Furthermore, as commonly done in practical applications of the diffuse interface method (see, e.g.~\cite{Stoter2017}), the phase field functions $\Phi_F^\eps$ and $\Phi_B^\eps$ are regularised as described in~\eqref{eq:regularizacija}. Otherwise, the ensuing system matrices become singular, as all the degrees of freedom that belong to one diffuse subdomain but not the other would have a zero row in these matrices. The value of the regularisation parameter $\delta$ is indicated in each example. 
The numerical method used in this section is summarised as follows: For $n = 0, \ldots, N - 1$, we seek $(\bu_h^{n + 1}, \pi_h^{n + 1}, \bfxi_h^{n + 1}, p_h^{n + 1} )$ such that, for every $(\bv_{h}, \zeta_h, \bfphi_h, q_h)$, we have:
\begin{equation}
    \begin{aligned}\label{eq:comp_method}
        &\rho_F \int_{\Omega} \dtee \bu_{h}^{n + 1} \cdot\bv_{h} \Phi_F^{\eps, \delta} 
        + 2\mu_F\int_{\Omega}\bD(\bu_h^{n + 1} ):\bD(\bv_h)\Phi_F^{\eps, \delta} 
        - \int_{\Omega }(\nabla\cdot\bv_h)\pi_h^{n + 1}\Phi_F^{\eps, \delta} 
        + \int_{\Omega}(\nabla\cdot\bu_h^{n + 1})\zeta_h \Phi_F^{\eps, \delta}  
        \\
        &+ \rho_B\int_{\Omega} \dtee \bfxi^{n + 1}_h \cdot \bfphi_h \Phi_B^{\eps, \delta}  
        + 2\mu_B \Delta t \int_{\Omega}\bD(\bfxi_h^{n + 1} ):\bD(\bfphi_h)\Phi_B^{\eps, \delta} 
        + 2\mu_B  \int_{\Omega}\bD(\bfeta_h^{n} ):\bD(\bfphi_h)\Phi_B^{\eps, \delta} 
        \\
        &  
        + \lambda_B \Delta t \int_{\Omega}\nabla \cdot \bfxi_h^{n + 1}  \nabla \cdot \bfphi_h \Phi_B^{\eps, \delta} 
        + \lambda_B\int_{\Omega}\nabla \cdot \bfeta_h^{n} \nabla \cdot \bfphi_h \Phi_B^{\eps, \delta} 
        + c_0\int_{\Omega} \dtee p_{h}^{n + 1} q_{h} \Phi_B^{\eps, \delta} 
        \\
        & 
        + \int_{\Omega} \boldsymbol\kappa_\eps\nabla p_h^{n + 1} \cdot \nabla q_h \Phi_B^{\eps, \delta} 
        - \alpha \int_{\Omega} \nabla \cdot \bfphi_h p_h^{n + 1} \Phi_B^{\eps, \delta} 
        + \alpha \int_{\Omega} \nabla \cdot \bfxi_h^{n + 1} q_h \Phi_B^{\eps, \delta} 
        \\
        &
        + \int_{\Omega} p_h^{n + 1} \bfphi_h \cdot \nabla  \Phi_B^{\eps, \delta}  
        - \int_{\Omega} q_h \bfxi_h^{n + 1} \cdot \nabla \Phi_B^{\eps, \delta} 
        + \int_{\Omega} q_h \bu_h^{n + 1} \cdot \nabla \Phi_F^{\eps, \delta}  
        - \int_{\Omega} p_h^{n + 1} \bv_h \cdot \nabla \Phi_F^{\eps, \delta}  
        \\ &
        + \alpha_{BJ}  \sum_{i = 1}^{d - 1} \int_{\Omega}((\bu_h^{n + 1} - \bfxi_h^{n + 1})\cdot \tilde{\bftau}_i) ((\bv_h - \bfphi_h)\cdot \tilde{\bftau}_i)|\nabla  \Phi_F^{\eps, \delta}  |
        \\ 
        &
        = \int_{\Omega}\bF_F^{n + 1}\cdot\bv_h\Phi_F^{\eps, \delta} 
        + \int_{\Omega}\bF_B^{n + 1}\cdot\bfphi_h\Phi_B^{\eps, \delta} 
        + \int_{\Omega}g^{n + 1} q_h \Phi_B^{\eps, \delta} .
    \end{aligned}
\end{equation}
We note that, in this case, the normal vectors are defined as $\bn_i = -\displaystyle\frac{\nabla \Phi_i^{\eps}}{|\nabla \Phi_i^{\eps}|}, i = F, B$, and the tangent vectors are constructed as described in~\cite{Stoter2017}. For simplicity, the structure problem is written in terms of the structure velocity, $\bfxi_h$. After the problem is solved, the displacement is found using
\begin{equation*}
    \bfeta^{n + 1}_h = \Delta t \bfxi_h^{n + 1} + \bfeta^n_h.
\end{equation*}
Note that it indeed appears natural to solve for the discrete velocity $\dtee\bfeta_h$, as reflected in the division-by-$\dt$ rescaling in~\eqref{eq:rescaled-discrete-form} and the error estimate~\eqref{eq:overall-error-estimate}.

\subsection{Rates of convergence}\label{sub:rates}

The first example is based on the method of manufactured solutions, which is used to compute rates of convergence. The computational domain is defined as $\Omega = (0, 1) \times (-1, 1)$, where the top half corresponds to the fluid domain and the bottom half to the poroelastic domain. 
The exact solutions are given by:
\begin{align*}
        &\bu_{ref} = \pi \cos(\pi t)
        \begin{bmatrix}
            - 3x + \cos(y)\\
            y + 1
        \end{bmatrix}, 
        \\
        &{\pi_{ref}} = e^t\sin(\pi x)\cos\left(\frac{\pi y}{2}\right) + 2\pi \cos(\pi t), 
        \\
        &\bfeta_{ref} = \sin(\pi t)
        \begin{bmatrix}
            -3x + \cos(y)\\
            y + 1
        \end{bmatrix}, 
        \\
        &{p_{ref}} = e^t\sin(\pi x)\cos\left(\frac{\pi y}{2}\right)
    \end{align*}
The forcing terms $\bF_F, \bF_B$, and $g$ are computed using the exact solutions and $\bfxi_{ref} = \partial_t \bfeta_{ref}$. Furthermore, to account for the fact that  the exact fluid velocity is not divergence-free, an additional forcing term is added to the mass conservation equation:
\begin{equation*}
    \nabla \cdot \bu = h.
\end{equation*}
Neumann boundary conditions are applied on the top boundary for the fluid problem. Dirichlet conditions are used on all other boundaries. Dirichlet boundary conditions are also used for the displacement and Darcy pressure on all boundaries. 

The following physical parameters are used: $\rho_B = \mu_B = \lambda_B = \alpha = c_0 = \gamma = \rho_F = \mu_F = 1$, and $\boldsymbol{\kappa} = \pmb{\mathbf{I}}$. The final time is $T = 0.8$s. We consider both a Lipschitz phase field function, defined as
\begin{equation*}
    \Phi_F^{\epsilon} = \frac12 \left(1 + \tanh \left(\frac{y}{\eps} \right) \right),
\end{equation*}
and a power distance weight function, defined by~\eqref{eq:defofPhiPsi} with $\beta = 0.9$. 
Both of them are regularised as in~\eqref{eq:regularizacija}. Notice that the Lipschitz weight is smooth near the interface, unlike the power distance weight.

We use $\mathbb{P}_2 - \mathbb{P}_1$ elements for the Stokes velocity and pressure, and $\mathbb{P}_2$ elements for the Biot pressure, the structure displacement, and velocity. We initially set $\Delta t = 0.1, h = 0.2$, $\eps = h$, and $\delta = 10^{-3}$. These parameters, including $\delta$, are then refined by halving them. 
The relative errors are defined as
\begin{equation*}
    \begin{aligned}
        \be_\bu &= \frac{\|\bu_{ref} - \bu_h^N \|_{L^2(\Omega, \Phi_F^{\eps, \delta})^2}}{\| \bu_{ref} \|_{L^2(\Omega, \Phi_F^{\eps, \delta})^2}},
        &\be_{\bfeta} &= \frac{\| \bfeta_{ref} - \bfeta_h^N \|_{E, \eps, \delta}}{\| \bfeta_{ref} \|_{E, \eps, \delta}},
        \\
        \be_{\partial_t \bfeta} &= \frac{\|\partial_t \bfeta_{ref} - \dtee \bfeta_h^{N} \|_{L^2(\Omega, \Phi_B^{\eps, \delta})^2}}{\| \partial _t \bfeta_{ref} \|_{L^2(\Omega,\Phi_B^{\eps, \delta})^2}},
        & e_p &= \frac{\|p_{ref} - p_h^N \|_{L^2(\Omega,\Phi_B^{\eps, \delta})}}{\| p_{ref} \|_{L^2(\Omega,\Phi_B^{\eps, \delta})}},
    \end{aligned}
\end{equation*}
evaluated at the final time $t^N = T$.

In addition to the method based on the backward Euler time discretisation described in~\eqref{eq:comp_method}, we also consider a method based on the midpoint scheme. The solution to the midpoint scheme is obtained by solving~\eqref{eq:comp_method} over a half time interval, using $\Delta t/2$, resulting in $(\bu_h^{n + \frac12}, \pi_h^{n + \frac12}, \bfxi_h^{n + \frac12}, \bfeta_h^{n + \frac12}, p_h^{n + \frac12} )$, and then extrapolating the solution as:
\begin{equation*}
    \bu^{n + 1}_h = 2 \bu^{n + \frac12}_h - \bu^n_h, \quad
    \bfxi^{n + 1}_h = 2 \bfxi^{n + \frac12}_h - \bfxi^n_h, \quad
    \bfeta^{n + 1}_h = 2 \bfeta^{n + \frac12}_h - \bfeta^n_h, \quad
    p^{n + 1}_h = 2 p^{n + \frac12}_h - p^n_h, 
\end{equation*}
at each time step, as described in~\cite{Burkardt2020}.

\begin{table}[ht]
    \begin{tabular}{c | c c | c c | c c | c c} 
        $h$ & $\be_\bu$ & rate & $e_p$ & rate & $\be_{\partial_t \bfeta}$ & rate & $\be_{\bfeta}$ & rate
        \\
        \hline  
        $\frac15$ & $8.3 \cdot 10^{-3} $ & -- & $1.1 \cdot 10^{-1}$  & -- & $7.3 \cdot 10^{-2}$ & -- &  $9.9 \cdot 10^{-1}$ & --
        \\
        $\frac1{10}$ & $7.7\cdot 10^{-3}$ & 0.11 & $8.1 \cdot 10^{-2}$ & 0.40 & $4.3 \cdot 10^{-2}$ & 0.77&  $3.3  \cdot 10^{-1}$ & 1.59
        \\
        $\frac1{20}$ & $4.0\cdot 10^{-3}$ & 0.96 & $5.3  \cdot 10^{-2}$ & 0.59  & $2.3  \cdot 10^{-2}$  & 0.86 & $1.4  \cdot 10^{-1}$ & 1.23
        \\
        $\frac1{40}$ & $2.0 \cdot 10^{-3}$ & 0.98 & $3.2 \cdot 10^{-2}$ & 0.76 & $1.2  \cdot 10^{-2}$  & 0.91 &  $ 6.5  \cdot 10^{-2}$  & 1.1
        \\
        $\frac1{80}$ & $1.0 \cdot 10^{-3}$ & 0.98 & $1.7  \cdot 10^{-2}$ & 0.86 & $6.5  \cdot 10^{-3}$  & 0.95 & $3.1  \cdot 10^{-2}$ & 1.05
        \\
    \end{tabular}       
    \caption{Rates of convergence for the fluid velocity, Biot pressure, structure velocity, and the displacement obtained using a Lipschitz phase field function and the backward Euler time discretisation.}\label{tab:ratesBE}
\end{table}

\begin{table}[ht]
    \begin{tabular}{c | c c | c c | c c | c c} 
        $h$ & $\mathbf{e}_\bu$ & rate & $e_p$ & rate & $\mathbf{e}_{\partial_t \bfeta}$ & rate & $\be_{\bfeta}$ & rate
        \\
        \hline  
        $\frac15$ & $2.7 \cdot 10^{-2} $ & -- & $7.5 \cdot 10^{-2}$  & -- & $7.1 \cdot 10^{-2}$ & -- &  $9.9 \cdot 10^{-1}$ & --
        \\
        $\frac1{10}$ & $1.4\cdot 10^{-2}$ & 1.0 & $6.8 \cdot 10^{-2}$ & 0.15 & $4.3 \cdot 10^{-2}$ & 0.75&  $3.3  \cdot 10^{-1}$ & 1.59
        \\
        $\frac1{20}$ & $6.9\cdot 10^{-3}$ & 1.0 & $4.7  \cdot 10^{-2}$ & 0.54  & $2.4  \cdot 10^{-2}$  & 0.85 & $1.4  \cdot 10^{-1}$ & 1.23
        \\
        $\frac1{40}$ & $3.4 \cdot 10^{-3}$ & 0.99 & $2.8 \cdot 10^{-2}$ & 0.73 & $1.3  \cdot 10^{-2}$  & 0.91 &  $ 6.5  \cdot 10^{-2}$  & 1.1
        \\
        $\frac1{80}$ & $1.7 \cdot 10^{-3}$ & 0.99 & $1.6  \cdot 10^{-2}$ & 0.84 & $6.6  \cdot 10^{-3}$  & 0.94 & $3.1  \cdot 10^{-2}$ & 1.05
        \\
    \end{tabular}       
\caption{Rates of convergence for the fluid velocity, Biot pressure, structure velocity, and the displacement obtained using a power distance weight function and the backward Euler time discretisation.}\label{tab:ratesBE2}
\end{table}

Table~\ref{tab:ratesBE} shows the rates of convergence obtained using a Lipschitz phase field function, and Table~\ref{tab:ratesBE2} shows the rates of convergence obtained using the power distance phase field function, both obtained using the backward Euler time discretisation.  We note that for the Lipschitz weight, the leading order error  is due to the time discretisation, which is only first order accurate, while for the distance weight, the leading order error is $\mathcal{O}(\eps^{1/4})$. However, in both cases, the first order convergence is obtained computationally for all variables except for the Biot pressure, whose rate seems to be  increasing  slowly. 

\begin{table}[ht]
    \begin{tabular}{c | c c | c c | c c | c c} 
        $h$ & $\mathbf{e}_\bu$ & rate & $e_p$ & rate & $\mathbf{e}_{\partial_t \bfeta}$ & rate & $\be_{\bfeta}$ & rate
        \\
        \hline
        $\frac15$ & $9.9\cdot 10^{-3} $ & -- & $3.0\cdot 10^{-2}$  & -- & $1.5 \cdot 10^{-2}$ & -- &  $4.6 \cdot 10^{-2}$ & --
        \\
        $\frac1{10}$ & $2.8\cdot 10^{-3}$ &1.81 & $1.2 \cdot 10^{-2}$ & 1.27 & $4.6 \cdot 10^{-3}$ & 1.8 &  $1.4  \cdot 10^{-2}$ & 1.7
        \\
        $\frac1{20}$ & $7.8\cdot 10^{-4}$ & 1.87 & $3.5  \cdot 10^{-3}$ & 1.80  & $1.2  \cdot 10^{-3}$  & 1.95 & $4.8  \cdot 10^{-3}$ & 1.56
        \\
        $\frac1{40}$ & $1.9 \cdot 10^{-4}$ & 2.01 & $8.9 \cdot 10^{-4}$ & 1.99 & $2.9  \cdot 10^{-4}$  & 2.03 &  $1.6  \cdot 10^{-3}$  & 1.58
        \\
        $\frac1{80}$ & $4.6 \cdot 10^{-5}$ & 2.08 & $2.2  \cdot 10^{-4}$ & 2.01 & $7.1  \cdot 10^{-5}$  & 2.04 & $5.4  \cdot 10^{-4}$ & 1.56
        \\
    \end{tabular}       
\caption{Rates of convergence for the fluid velocity, Biot pressure, structure velocity, and the displacement obtained using a Lipschitz phase field function and the midpoint time discretisation.}\label{tab:ratesMP}
\end{table}

\begin{table}[ht]
    \begin{tabular}{c | c c | c c | c c | c c} 
        $h$ & $\mathbf{e}_\bu$ & rate & $e_p$ & rate & $\mathbf{e}_{\partial_t \bfeta}$ & rate & $\be_{\bfeta}$ & rate
        \\
        \hline  
        $\frac15$ & $9.3\cdot 10^{-3} $ & -- & $2.3\cdot 10^{-2}$  & -- & $1.3 \cdot 10^{-2}$ & -- &  $4.3 \cdot 10^{-2}$ & --
        \\
        $\frac1{10}$ & $2.4\cdot 10^{-3}$ &1.95 & $6.8 \cdot 10^{-3}$ & 1.75 & $3.3 \cdot 10^{-3}$ & 1.97 &  $1.1  \cdot 10^{-2}$ & 1.94
        \\
        $\frac1{20}$ & $6.1\cdot 10^{-4}$ & 1.99 & $1.8  \cdot 10^{-3}$ & 1.89  & $8.8  \cdot 10^{-4}$  & 1.89 & $3.4  \cdot 10^{-3}$ & 1.69
        \\
        $\frac1{40}$ & $1.5 \cdot 10^{-4}$ & 2.02 & $5.0 \cdot 10^{-4}$ & 1.88 & $2.2  \cdot 10^{-4}$  & 1.98 &  $1.1  \cdot 10^{-3}$  & 1.63
        \\
        $\frac1{80}$ & $3.7 \cdot 10^{-5}$ & 2.03 & $1.3  \cdot 10^{-4}$ & 1.89 & $5.7  \cdot 10^{-5}$  & 1.98 & $3.6  \cdot 10^{-4}$ & 1.60
        \\
    \end{tabular}       
\caption{Rates of convergence for the fluid velocity, Biot pressure, structure velocity, and the displacement obtained using a power distance weight function and the midpoint time discretisation.}\label{tab:ratesMP2}
\end{table}

Tables~\ref{tab:ratesMP} and~\ref{tab:ratesMP2} show the errors and the rates of convergence obtained using the Lipschitz phase field function and the power distance phase field function, respectively, both obtained with the midpoint time discretisation. In case of the Lipschitz weight, the modelling error is the leading order error, with the rate of $\mathcal{O}(\eps^{3/2})$. The leading order error for the power distance weight is still $\mathcal{O}(\eps^{1/4})$. We observe the second order convergence for the fluid velocity, Biot pressure, and the structure velocity. However, the rates for displacement, measured in the energy norm, seem to be approaching 1.5. This indicates that both Lipschitz and power distance weights seem to be converging with the same order, and while for some variables the rates exceed theoretical predictions, the modelling error still dominates the rates of convergence for the displacement. 


\subsection{Comparison with a sharp interface model}\label{sec:ex2a}

The focus of this section is on the comparison of the results obtained using a sharp interface model and a diffuse interface model in a simplified 3D domain describing the flow in a channel surrounded by a poroelastic medium. The geometry consists of a cylinder with radius 0.5 and length 2 embedded in a poroelastic cube.
\begin{figure}[ht]
    \begin{center}
        \includegraphics[scale = 0.9]{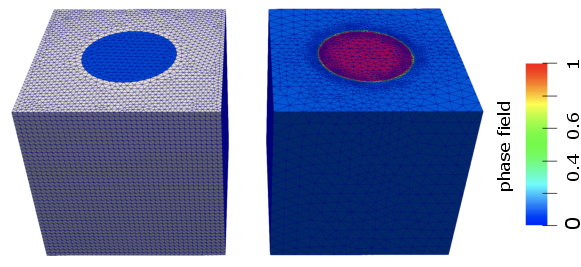}    
    \end{center}
    \caption{Geometry and the computational mesh used in Example~\ref{sec:ex2a} for the sharp (left) and the diffuse (right) interface model. The right panel also shows the phase field function used for the diffuse interface model.}\label{fig:geometry_simple}
\end{figure}

Figure~\ref{fig:geometry_simple} shows the sharp interface domain and the computational mesh (left), and the phase field function and the computational mesh on the diffuse interface domain geometry (right). The phase field function is defined as 
\begin{equation*}
    \Phi_F^{\epsilon} = \frac12 \left(1 + \tanh\left(- \frac{ x^2 + y^2-0.5^2}{\eps}\right) \right)
\end{equation*}
and regularised as in~\eqref{eq:regularizacija}. We use $\eps = 0.0095$ and $\delta = 0.001$. The flow is driven by a pressure drop imposed at the inlet and outlet sections of the cylinder: 
\begin{equation*}
    \bfsigma_F(\bu, \pi) \bn  = 
    \begin{dcases}
        - 10 \bn, & \quad \textrm{at} \; z = 0 \;\; \textrm{(bottom boundary)},
        \\
        0, & \quad \textrm{at} \; z = 2 \;\; \textrm{(top boundary)}.
    \end{dcases}
\end{equation*}
On the same boundaries (top and bottom), the structure is assumed to be fixed ($\boldsymbol \eta = 0$), and zero flux conditions are imposed for the Biot pressure ($\boldsymbol \kappa \nabla p \cdot \bn = 0$). On all the sides we impose  zero Neumann conditions for the displacement and zero Biot pressure. 
The problem is solved using parameters specified in Table~\ref{tab:example2Par} with the time step of $\Delta t = 5 \cdot 10^{-2}$ until a steady state is reached. 

\begin{table}[ht]
    \begin{tabular}{l l | l l }
        \textbf{Parameters} & \textbf{Values} & \textbf{Parameters} & \textbf{Values}  \\
        \hline
        \hline
        \textbf{Fluid density} $\rho_F$ (g/cm$^3$)& $1$ &\textbf{Dynamic viscosity} $\mu$ (poise) & $0.035$    \\
        \textbf{Structure density} $\rho_B$ (g/cm$^3$)& $1$ & \textbf{Young's modulus} $E$ (dyne/cm$^2$)  &      $5\cdot 10^5$  \\
        \textbf{Poisson's ratio} $\nu$ & $0.49$ & \textbf{Slip rate} $\alpha_{BJ} $ (g/cm$^2$ s)& $10$ \\
        \textbf{Storativity coeff.} $c_0$ (cm$^2$/dyne) & $10^{-3}$ & \textbf{Hydraulic conductivity} $\boldsymbol \kappa$ (cm$^3$ s/g) & $10^{-5} \mathbf{I}$     \\
        \textbf{Biot--Willis constant} $\alpha$ & 1 & \\
    \end{tabular}
    \caption{The parameters used in Examples~\ref{sec:ex2a} and~\ref{sec:ex2}.}\label{tab:example2Par}
\end{table}

For both sharp interface and diffuse interface problems, we use $\mathbb{P}_1$ elements for all variables, with the following stabilisation added to the fluid problem:
\begin{equation*}
    \gamma_{\rm stab} h^2 \int_\Omega \nabla \pi \cdot \nabla \zeta,
\end{equation*}
where $\gamma_{\rm stab} = 2 \cdot 10^{-3}$.

\begin{figure}[ht]
    \begin{center}
        \includegraphics[scale = 0.8]{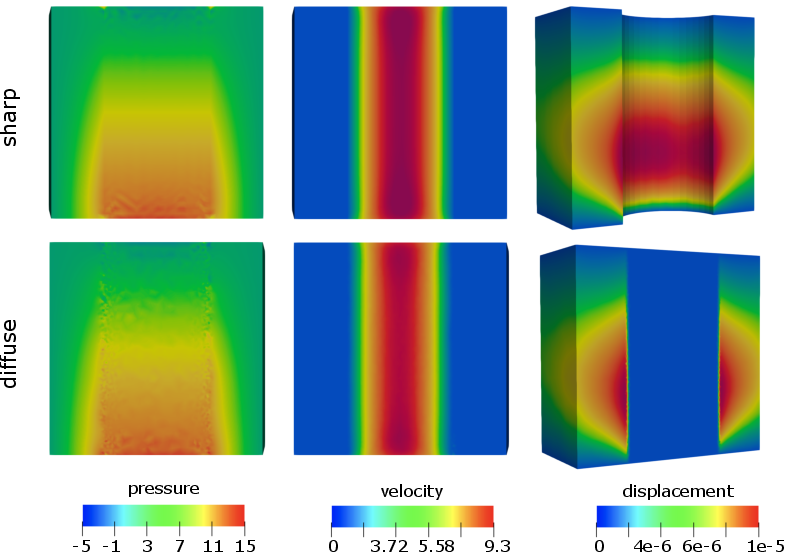}    
    \end{center}
    \caption{Comparison of the Stokes and Biot pressure (left), the velocity (middle), and the displacement (right) obtained using the sharp interface model (top) and the diffuse interface model (bottom), shown on a cross section of the domain.}\label{fig:cyl_comp}
\end{figure}

Figure~\ref{fig:cyl_comp} shows a comparison of the Stokes and Biot pressure, velocity and displacement obtained using the sharp interface model and the diffuse interface model on a cross section of the domain. Overall, an excellent agreement is observed. 

\subsection{Flow in a complex network}\label{sec:ex2}

\begin{figure}[ht]
    \begin{center}
        \includegraphics[scale = 0.55]{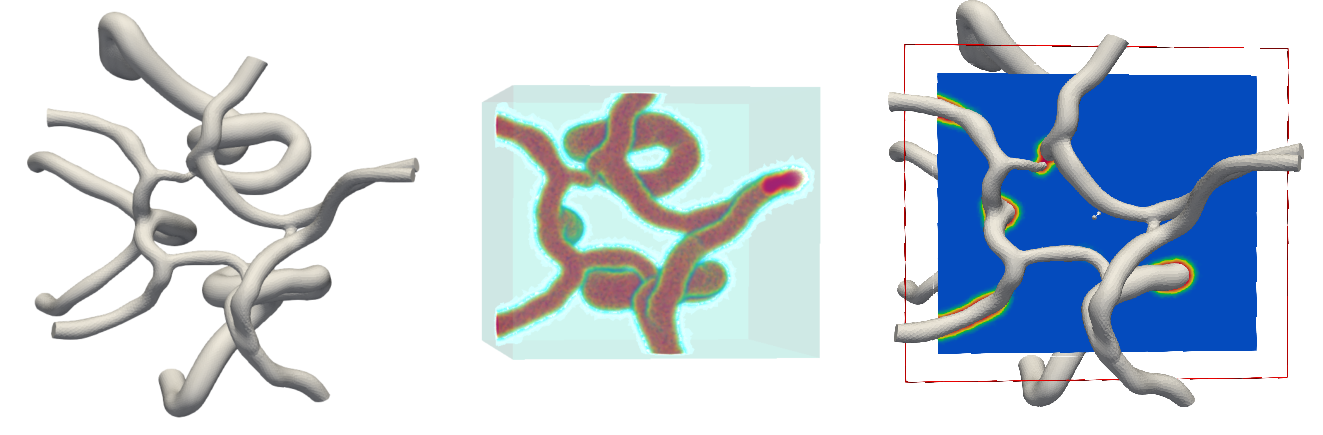}    
    \end{center}
    \caption{Left: A patient-specific vascular network containing the circle of Willis. Middle: A phase field function for the section of the network considered in our study. Right: A phase field function on a cross section of our computational domain, superimposed with the full vascular network.}\label{fig:geometry}
\end{figure}

In this example, we model fluid flow in a complex network, corresponding to a patient specific model of vasculature containing the circle of Willis, obtained from~\cite{Wilson2013}. We consider a box-like section of tissue around the network, which we assume is poroelastic. The patient specific geometry and our constructed phase field function are shown in Figure~\ref{fig:geometry}.  

The left panel shows the full vascular network containing the circle of Willis, and the middle panel shows the phase field function for the section of the network considered in our study. The right panel shows the phase field function on a cross section of our computational domain (blue box), superimposed with the full vascular network. The phase field function was obtained  by identifying nodes in our box geometry which are in the vascular network, and setting their values equal to 1 (otherwise, they were set to zero). The function was then smoothed out by solving the Allen--Cahn equation as in~\cite{Stoter2017}.

\begin{figure}[ht]
    \begin{center}
        \includegraphics[scale = 0.85]{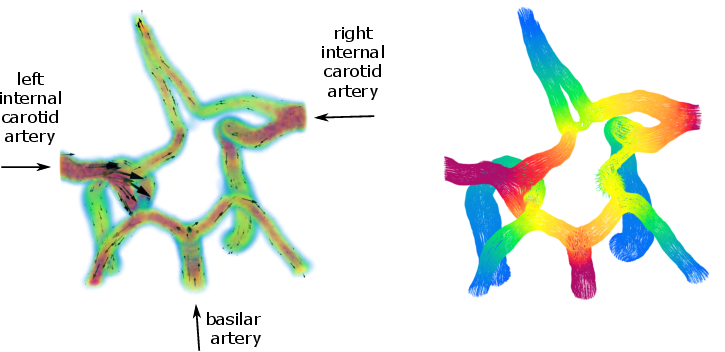}    
    \end{center}
    \caption{Left: Velocity magnitude superimposed with vectors indicating the flow direction. The marked arteries represent the inlet sections. Right: Flow streamlines coloured by the Stokes pressure.}\label{fig:microVel}
\end{figure}

To solve this problem, we impose boundary conditions which mimic the flow pattern in the circle of Willis. In particular, we impose the velocity of 35 cm/s~\cite{Cantu1992} at the inflow of the basilar artery, and the velocity of 50 cm/s at the inflow of the right and left internal carotid arteries~\cite{Devault2008, Alimi2017}, see the left panel of Figure~\ref{fig:microVel}. At all other outlets, we impose zero normal stress for the fluid problem. At the inflow boundaries, we impose zero displacement, while zero normal poroelastic stress is imposed at the remaining boundaries. Finally, zero Biot pressure is imposed at all the boundaries. The problem is solved with the parameters specified in Table~\ref{tab:example2Par}, using the same finite elements as for the cylindrical geometry case studied in Example~\ref{sec:ex2a}. As before, we use $\delta = 0.001$ and $\Delta t = 5\cdot 10^{-2}$. The problem is solved until a steady state is reached.

The left panel of Figure~\ref{fig:microVel} shows the magnitude of the velocity in the network, superimposed with velocity vectors, which show the flow direction. Marked arteries indicate the inflow sections. The flow streamlines are shown in the right panel, coloured by the fluid pressure. 

\begin{figure}[ht]
    \begin{center}
        \includegraphics[scale = 0.9]{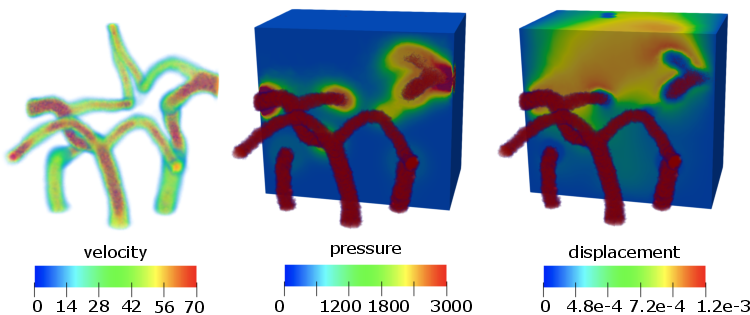}
    \end{center}
    \caption{The left panel shows the magnitude of the velocity. The middle and the right panels show the Stokes and the Darcy pressure, and the displacement, respectively, both superimposed with the outline of the vascular network.}\label{fig:microPEta}
\end{figure}
The velocity magnitude is again shown in the left panel of Figure~\ref{fig:microPEta} from a different perspective. In Figure~\ref{fig:microPEta} we also show the Stokes and Biot pressure (middle panel) and the displacement (right panel), both displayed on a cross section of the domain, superimposed with the outline of the vascular network. We can observe larger pressure values in the vicinity of the network, as well as larger displacement. We note that from the perspective shown in Figure~\ref{fig:microPEta} that the displacement is fixed on the bottom and on the left and right side. Hence, the poroelastic region surrounding the network displaces more in the top region of the domain.

\section{Conclusions}

In this work, we analysed the diffuse interface method for fluid-poroelastic structure interaction. We proved that the problem is well-posed, and derived rates of convergence of the discrete diffuse interface formulation to the continuous sharp interface formulation. In particular, for piecewise quadratic velocities, we show that in case of a power distance weight, the error is 
\begin{equation*}
    \mathcal{O}(\Delta t) + \mathcal{O}(h^2) + \mathcal{O}(\eps^{\frac14}).
\end{equation*}
With the assumption that the weight is Lipschitz, the estimate of the modelling error improves, and gives
\begin{equation*}
    \mathcal{O}(\eps^{\frac32}) + \mathcal{O}(\delta^{\frac12}), 
\end{equation*}
where $\delta$ is the regularisation parameter. The latter estimate agrees with what we observe in numerical experiments when either regularised Lipschitz or power distance weights are used. In particular, numerical results indicate that the rate of convergence obtained computationally for all variables exceeds the theoretical prediction in the case of power distance weights, and is more similar to the results obtained for Lipschitz weights. 

While we combine the approximation and the modelling errors in order to obtain the total error, we note that one needs to be careful in which order the limits are taken. In particular, we first let $(h, \dt)\to (0, 0)$, and then $\eps\to 0$. This is due to the fact that, as expected, all the implicit constants in our error estimates depend on higher order norms of the diffuse interface solution, which 
depend on the interface parameter $\eps$.
 We also note that in case when weights are regularised, 
in principle 
we can also take $\delta \rightarrow 0$. However in our implementation, this scenario is problematic since in that case the system matrix becomes singular.

While in this paper we assume that the interface is fixed and that the phase field function does not change in time, this is a first step towards analysing more complex models where the interface is evolving, such as problems with large deformations, which we will study in future work. However, even in the case considered here, the diffuse interface method is appealing in many applications, such as simulations of biomedical flows in complex geometries.

\section*{Acknowledgements}

FA is supported by a Society of Science Postdoctoral fellowship from the College of Science at the University of Notre Dame.
MB is partially supported by NSF grants DMS-2208219 and DMS-2205695. AJS is partially supported by NSF grant DMS-2111228. BM is supported by by Croatia-USA bilateral grant ``The mathematical framework for the diffuse interface method applied to coupled problems in fluid dynamics'' and by the Croatian Science Foundation, project
number IP-2022-10-2962.
Part of this work was completed while AJS was in residence at the Institute for Computational and Experimental Research in Mathematics (ICERM) in Providence, RI, during the Numerical PDEs: Analysis, Algorithms, and Data Challenges semester program. ICERM is supported by NSF grant DMS-1929284.

\appendix

\section{Consistency and projection estimates}\label{sec:taylor-expansions}

Here we collect some consistency estimates, as well as some properties of the weighted projections defined in Lemma~\ref{lem:weighted-ritz}, that are useful in the course of our error analysis.
We begin with a somewhat standard consistency estimate.

\begin{lemma}[time-consistency]\label{lem:ConsistencyDerivatives}
For $n = 0, \ldots, N - 1$ define $\cR_\eps^{n + 1} \in (\mathcal{V}^\eps_F \times \mathcal{V}^\eps_B \times \mathcal{X}^\eps)^*$ via
    \begin{align*}
        \langle \cR_\eps^{n + 1}&, (\bv, \bfphi, q) \rangle 
        = \rho_F \int_{\Omega_F^\eps} (\dtee - \partial_t) (\bu^\eps)^{n + 1} \cdot \bv \Phi_F^\eps 
        + \rho_B\int_{\Omega^\eps_B}(\dtt - \partial_{tt})(\bfeta^\eps)^{n + 1}\cdot\bfphi\Phi_B^\eps 
        \nonumber\\
        &+ c_0 \int_{\Omega_B^\eps} (\dtee - \partial_t)(p^\eps)^{n + 1} q \Phi_B^\eps
        - \frac{1}{2\eps}\int_{\ell^\eps}q(\dtee - \partial_t)(\bfeta^\eps)^{n + 1}\cdot\nabla\dist_\Gamma
        \\
        &+ \frac{\alpha_{BJ}}{2\eps}\sum_{i = 1}^{d - 1}\int_{\ell^\eps}(-(\dtee - \partial_t)(\bfeta^\eps)^{n + 1}\cdot\tilde{\bftau}_i)((\bv - \bfphi)\cdot\tilde{\bftau}_i)\vert\nabla\dist_\Gamma\vert
        + \alpha\int_{\Omega_B^\eps}\nabla\cdot(\dtee - \partial_t)(\bfeta^\eps)^{n + 1}q\Phi_B^\eps.
        \nonumber
    \end{align*}
Then, assuming~\eqref{eq:SolIsSmooth}, 
for every $1\leq K \leq N$ 
we have
    \begin{align*}
        &\dt\sum_{n = 0}^{K - 1}\langle\cR^{n + 1}_{\eps}, (\bE^{n + 1}_{\bu, h}, \dtee\bE^{n + 1}_{\bfeta, h}, E^{n + 1}_{p, h})\rangle
        \leq
        C(\dt)^2\left(
        \|\partial_{tt}\bu^\eps\|^2_{L^2(0, T; L^2(\Omega^\eps_F, \Phi^\eps_F)^d)}
        + \|\partial_{ttt}\bfeta^\eps\|^2_{L^\infty(0, T; L^2(\Omega^\eps_B, \Phi^\eps_B)^d)}
        \right.
        \nonumber\\ &
        \left.
        + \|\partial_t^4\bfeta^\eps\|^2_{L^2(0, T; L^2(\Omega^\eps_B, \Phi^\eps_B)^d)}
        + \| \partial_t \bfeta^\eps \|_{L^2(0, T; L^2(\Omega_B^\eps, \Phi_B^\eps)^d)}
        + \|\partial_{tt}p^\eps\|^2_{L^2(0, T; L^2(\Omega^\eps_B, \Phi^\eps_B))}
        + \|\partial_{tt}\bfeta^\eps\|^2_{L^2(0, T; \mathcal{V}^\eps_B)}
        \right)
        \nonumber\\ &
        + \frac{\mu_F}{4C_K}\dt\sum_{n = 0}^{K - 1}\|\bE^{n + 1}_{\bu, h}\|^2_{L^2(\Omega^\eps_F, \Phi^\eps_F)^d}
        + \frac{\mu_B}{10C_K}\|\bE^K_{\bfeta, h}\|^2_{L^2(\Omega^\eps_B, \Phi^\eps_B)^d}
        + \frac{\mu_B}{10C_K}\dt\sum_{n = 0}^{K - 1}\|\bE^{n + 1}_{\bfeta, h}\|^2_{L^2(\Omega^\eps_B, \Phi^\eps_B)^d}
        \\ &
        + \frac{3k_*}{8}\dt\sum_{n = 0}^{K - 1}\|E^{n + 1}_{p, h}\|^2_{\mathcal{X}^\eps}
        + \frac{\alpha_{BJ}}{2}\sum_{i = 1}^{d - 1}\dt\sum_{n = 0}^{K - 1}\|(\bE^{n + 1}_{\bu, h} - \dtee\bE^{n + 1}_{\bfeta, h})\cdot\tilde{\bftau}_i\|^2_{L^2\left(\ell^\eps, \frac{1}{2\eps}|\nabla\dist_\Gamma|\right)},
        \nonumber
    \end{align*}
where 
$C > 0$
depends only on
physical parameters, and the projection errors $(\bE^{n + 1}_{\bu, h}, \bE^{n + 1}_{\bfeta, h}, E^{n + 1}_{p, h})$ were defined in~\eqref{eq:truncation-errors}.
\end{lemma}
\begin{proof}
    Applying discrete integration by parts,
    Young's inequality, 
    and 
    Lemma~\ref{lem:nablaPhi},
    we have
    \begin{multline*}
        \dt\sum_{n = 0}^{K - 1}\langle\cR^{n + 1}_{\eps}, (\bE^{n + 1}_{\bu, h}, \dtee\bE^{n + 1}_{\bfeta, h}, E^{n + 1}_{p, h})\rangle
        \leq 
        \frac{\rho_F^2C_K}{\mu_F}\dt\sum_{n = 0}^{K - 1}\|(\dtee - \partial_t)(\bu^\eps)^{n + 1}\|^2_{L^2(\Omega^\eps_F, \Phi^\eps_F)^d}
        \\
        + \frac{\mu_F}{4C_K}\dt\sum_{n = 0}^{K - 1}\|\bE^{n + 1}_{\bu, h}\|^2_{L^2(\Omega_F^\eps, \Phi_F^\eps)^d}
        + 
        \frac{5C_K\rho_B^2}{2\mu_B}\|(\dtt - \partial_{tt})(\bfeta^\eps)^K\|^2_{L^2(\Omega_B^\eps, \Phi_B^\eps)^d}
        + \frac{\mu_B}{10C_K}\|\bE^K_{\bfeta}\|^2_{L^2(\Omega_B^\eps, \Phi_B^\eps)^d}
        \\
        + 
        \frac{5C_K\rho_B^2}{2\mu_B}\dt\sum_{n = 0}^{K - 1}\|(\mathrm{d}_{ttt} - \partial_{tt}\dtee)(\bfeta^\eps)^{n + 1}\|^2_{L^2(\Omega_B^\eps, \Phi_B^\eps)^d}
        + \frac{\mu_B}{10C_K}\dt\sum_{n = 0}^{K - 1}\|\bE^n_\bfeta\|^2_{L^2(\Omega_B^\eps, \Phi_B^\eps)^d}
        \\
        + \frac{2c_0^2}{k_*}\dt\sum_{n = 0}^{K - 1}\|(\dtee - \partial_t)(p^\eps)^{n + 1}\|^2_{L^2(\Omega_B^\eps, \Phi_B^\eps)}
        + \frac{k_*}{8}\dt\sum_{n = 0}^{K - 1}\|E^{n + 1}_p\|^2_{L^2(\Omega_B^\eps, \Phi_B^\eps)}
        + \frac{4C_{\rm tr}^2}{k_*}\dt\sum_{n = 0}^{K - 1}\|(\dtee - \partial_t)(\bfeta^\eps)^{n + 1}\|^2_{\mathcal{V}^\eps_B}
        \\
        + \frac{k_*}{8}\dt\sum_{n = 0}^{K - 1}\|E^{n + 1}_p\|^2_{\mathcal{X}^\eps}
        + \frac{4\alpha^2}{k_*}\dt\sum_{n = 0}^{K - 1}\|\nabla\cdot(\dtee - \partial_t)(\bfeta^\eps)^{n + 1}\|^2_{L^2(\Omega_B^\eps, \Phi_B^\eps)^d}
        + \frac{k_*}{8}\dt\sum_{n = 0}^{K - 1}\|E^{n + 1}_p\|^2_{L^2(\Omega^\eps_B, \Phi^\eps_B)}
        \\
        + 
        \frac{\alpha_{BJ}C_{\rm tr}(d - 1)}{2}\dt\sum_{n = 0}^{K - 1}\|(\dtee - \partial_t)(\bfeta^\eps)^{n + 1}\|^2_{\mathcal{V}^\eps_B}
        + \frac{\alpha_{BJ}}{2}\sum_{i = 1}^{d - 1}\dt\sum_{n = 0}^{K - 1}\|(\bE^{n + 1}_{\bu} - \dtee\bE^{n + 1}_{\bfeta})\cdot\tilde{\bftau}_i\|^2_{L^2\left(\ell^\eps, \frac{1}{2\eps}|\nabla\dist_\Gamma|\right)}
        .
    \end{multline*}
    Exactly as in the first bound in~\cite[Lemma 4]{Bukac2015} but with $\Phi_F^\eps$ incorporated into the spatial integrands, we obtain
    \begin{equation*}
        \dt\sum_{n = 0}^{K - 1}\|(\dtee - \partial_t) (\bu^\eps)^{n + 1}\|^2_{L^2(\Omega^\eps_F, \Phi_F^\eps)^d} \leq (\dt)^2\|\partial_{tt}\bu^\eps\|^2_{L^2(0, T; L^2(\Omega^\eps_F, \Phi^\eps_F)^d)},
    \end{equation*}
    and precisely the analogous bound 
    follows also for $\dt\sum_{n = 0}^{K - 1}\|(\dtee - \partial_t)(p^\eps)^{n + 1}\|^2_{L^2(\Omega_B^\eps, \Phi_B^\eps)}$ and $\dt\sum_{n = 0}^{K - 1}\|(\dtee - \partial_t)(\bfeta^\eps)^{n + 1}\|^2_{\mathcal{V}^\eps_B}$.
    By a similar calculation,
    \begin{multline*}
        \|(\dtt - \partial_{tt})(\bfeta^\eps)^K\|^2_{L^2(\Omega^\eps_B; \Phi^\eps_B)^d}
        = 
        \int_{\Omega^\eps_B}\left|\frac{-1}{(\dt)^2}\left(\int_{t^{K - 2}}^{t^{K - 1}} - \int_{t^{K - 1}}^{t^K}\right)(t - t^K)\frac{\partial_{ttt}\bfeta^\eps}{2}\Dt\right|^2\Phi^\eps_B
        \\
        \leq
        \frac{1}{2(\dt)^4}\int_{\Omega^\eps_B}\left(\int_{t^{K - 1}}^{t^{K - 1}}|\partial_{ttt}\bfeta^\eps|^2\Dt\int_{t^{K - 2}}^{t^{K - 1}}|t - t^K|^4\Dt 
        + \int_{t^{K - 1}}^{t^K}|\partial_{ttt}\bfeta^\eps|^2\Dt\int_{t^{K - 1}}^{t^K}|t - t^K|^4\Dt\right)\Phi^\eps_B
        \\
        \leq
        16(\dt)^2\|\partial_{ttt}\bfeta^\eps\|^2_{L^\infty(0, T; L^2(\Omega^\eps_B; \Phi^\eps_B)^d)}.
    \end{multline*}
    
Finally,
\begin{multline*}
    \dt\sum_{n = 0}^{K - 1}\|(\mathrm{d}_{ttt} - \partial_{tt}\dtee)(\bfeta^\eps)^{n + 1}\|^2_{L^2(\Omega_B^\eps, \Phi_B^\eps)^d}
    \\
    = \dt\sum_{n = 0}^{K - 1}\int_{\Omega^\eps_B}\left|\frac{-1}{6(\dt)^3}\left(-\int_{t^{n - 2}}^{t^{n - 1}} + 2\int_{t^{n - 1}}^{t^n} - \int_{t^n}^{t^{n + 1}}\right) 
    (t - t^{n + 1})^3\partial_t^4\bfeta^\eps\Dt 
    + \frac{1}{\dt}\int_{t^n}^{t^{n + 1}}(t - t^{n + 1})\partial_t^4\bfeta^\eps\Dt\right|^2\Phi^\eps_B
    \\
    \leq
    4\sum_{n = 0}^{K - 1}\int_{\Omega^\eps_B}\left(\frac{1}{6(\dt)^2}\int_{t^{n - 2}}^{t^{n - 1}}|\partial_t^4\bfeta^\eps|^2\Dt
    \int_{t^{n - 2}}^{t^{n - 1}}|t - t^{n + 1}|^6\Dt 
    + \frac{1}{6(\dt)^2}\int_{t^{n - 1}}^{t^n}|\partial_t^4\bfeta^\eps|^2\Dt\int_{t^{n - 1}}^{t^n}|t - t^{n + 1}|^6\Dt 
        \right. \\ \left.
    + \frac{1}{6(\dt)^2}\int_{t^n}^{t^{n + 1}}|\partial_t^4\bfeta^\eps|^2\Dt\int_{t^n}^{t^{n + 1}}|t - t^{n + 1}|^6\Dt 
    + \int_{t^n}^{t^{n + 1}}|\partial_t^4\bfeta^\eps|^2\Dt\int_{t^n}^{t^{n + 1}}|t - t^{n + 1}|^6\Dt\right)\Phi^\eps_B
    \\
    \leq
    (486(\dt)^5 + 4(\dt)^7)\|\partial_t^4\bfeta^\eps\|^2_{L^2(0, T; L^2(\Omega^\eps_B, \Phi^\eps_B)^d)}.
\end{multline*}
The estimate now follows.
\end{proof}

Finally, we record some error estimates on the weighted projections defined in Lemma~\ref{lem:weighted-ritz}.

\begin{lemma}[projection errors]\label{lem:more-projection-errors}
Recall the projection errors defined in~\eqref{eq:truncation-errors}. Under the smoothness assumption~\eqref{eq:SolIsSmooth} on the exact solution, we have
    \begin{align*}
        &\| Y^K_{p}\|^2_{\mathcal{X}^\eps}
        + \|\bY^K_{\bu}\|^2_{\mathcal{V}^\eps_F}
        + \|\dtee\bY^K_{\bfeta}\|^2_{\mathcal{V}^\eps_B}
        + \dt\sum_{n = 0}^{K - 1}\|\dtee\bY^{n + 1}_{\bu}\|^2_{\mathcal{V}^\eps_F}
        + \dt\sum_{n = 0}^{K - 1}\|\dtt\bY^{n + 1}_{\bfeta}\|^2_{\mathcal{V}^\eps_B}
        \nonumber\\ &
        + \dt\sum_{n = 0}^{K - 1}\|\dtee Y^{n + 1}_{p}\|^2_{\mathcal{X}^\eps}
        + \dt\sum_{n = 0}^{K - 1}\|\dtee\bY^{n + 1}_{\bfeta}\|^2_{\mathcal{V}^\eps_B}
        + \dt\sum_{n = 0}^{K - 1}\|\bY^{n + 1}_{\bu}\|^2_{\mathcal{V}^\eps_F}
        \\ 
        &\lesssim 
        h^{2k}\left(
        \|p^\eps\|^2_{L^\infty(0, T; \mathcal{X}^{k + 1, \eps})}
        + \|\bu^\eps\|^2_{L^\infty(0, T; \mathcal{V}^{k + 1, \eps}_F)}
        + \|\theta^\eps\|^2_{L^\infty(0, T; \mathcal{Q}^{k + 1, \eps})}
        + \|\partial_t\bfeta^\eps\|^2_{L^\infty(0, T; \mathcal{V}^{k + 1, \eps}_B)}
        \right. \nonumber\\ & \left.
        + \dt\|\partial_{tt}\bfeta^\eps\|^2_{L^\infty(0, T; \mathcal{V}^{k + 1, \eps}_B)}
        + \|\partial_t\bu^\eps\|^2_{L^2(0, T; \mathcal{V}^{k + 1, \eps}_F)}
        + \|\partial_{tt}\bfeta^\eps\|^2_{L^2(0, T; \mathcal{V}^{k + 1, \eps}_B)}
        + \|\partial_t p^\eps\|^2_{L^2(0, T; \mathcal{X}^{k + 1, \eps})}
        \right. \nonumber\\ & \left.
        + \|\partial_t\theta^\eps\|^2_{L^2(0, T; \mathcal{Q}^{k, \eps})}
        \right)
        + (\dt)^2\left(
        \|\partial_{tt}\bu^\eps\|^2_{L^2(0, T; \mathcal{V}^\eps_F)}
        + \|\partial_{tt}\bfeta^\eps\|^2_{L^2(0, T; \mathcal{V}^\eps_B)}
        + \|\partial_{ttt}\bfeta^\eps\|^2_{L^2(0, T; \mathcal{V}^\eps_B)}
        \right. \nonumber\\ & \left.
        + \|\partial_{tt}p^\eps\|^2_{L^2(0, T; \mathcal{X}^\eps)}
        \right).
        \nonumber
    \end{align*}
\end{lemma}
\begin{proof}
The proof is standard; for the sake of brevity we only demonstrate a bound for the 3\textsuperscript{rd} term on the left hand side. Extending the interpolation errors to be defined for all time, a Taylor expansion then provides the existence of $z_K\in (t^{K - 1}, t^K)$ for which
\begin{equation*}
    \dtee\bY^K_{\bfeta} = \partial_t\bY^K_{\bfeta} + \frac{\dt}{2}\partial_{tt}\bY_{\bfeta}(z_K).
\end{equation*}
Consequently,
    \begin{align*}
        \|\dtee\bY^K_{\bfeta}\|^2_{\mathcal{V}^\eps_B} 
        &\leq 2\|\partial_t\bY^K_{\bfeta}\|^2_{\mathcal{V}^\eps_B} 
        + \dt\|\partial_{tt}\bY_{\bfeta}(z_K)\|^2_{\mathcal{V}^\eps_B}
        \nonumber\\
        &\lesssim 
        h^{2k}\|\partial_t(\bfeta^\eps)^K\|^2_{\mathcal{V}^{k + 1, \eps}_B}
        + \dt h^{2k}\|\partial_{tt}\bfeta^\eps(z_K)\|^2_{\mathcal{V}^{k + 1, \eps}_B}
        \\
        &\leq
        h^{2k}\left(
            \|\partial_t\bfeta^\eps\|^2_{L^\infty_\dt(0, T; \mathcal{V}^{k + 1, \eps}_B)} 
            + \dt\|\partial_{tt}\bfeta^\eps\|^2_{L^\infty_\dt(0, T; \mathcal{V}^{k + 1, \eps}_B)}
        \right).
        \nonumber
    \end{align*}

The remaining terms follow by 
Taylor expansion with integral remainder, as in the previous Lemma,
combined with Lemma~\ref{lem:weighted-ritz}, and the observation that the regularities~\eqref{eq:SolIsSmooth} allow us to bound, for instance,
\begin{equation*}
    \|p^\eps\|_{L^\infty_\dt(0, T; \mathcal{X}^{k + 1, \eps})} \lesssim \|p^\eps\|_{L^\infty(0, T; \mathcal{X}^{k + 1, \eps})},
\end{equation*}
with constants that are independent of $\dt$ or $h$.
\end{proof}

\begin{remark}[error decoupling]
    Due to the application of the 
    Stokes-elastic
    projection in Lemma~\ref{lem:weighted-ritz}, there is dependence on the norms of $\theta^\eps$ on the right-hand-side of our main error estimate~\eqref{eq:overall-error-estimate} --- including in the contribution from the bound in Lemma~\ref{lem:more-projection-errors} --- even in the absence of an error estimate for $\theta^\eps$ (or $\pi^\eps$) in~\eqref{eq:overall-error-estimate}.
    We conjecture that this dependence may be removed by decoupling the errors of the three separate fields $(\mathfrak{u}, \vartheta, \mathfrak{p})$ in the first line of the projection error~\eqref{eq:projection-error-estimate}
    using finite element families which possess some generalisation of pressure-robustness to the (weighted) 3-field Stokes-elastic problem~\eqref{eq:weighted-stokes-proj}, but we do not pursue this.
\end{remark}

\bibliographystyle{plain}
\bibliography{bibfile}

\end{document}